\magnification=\magstep1
\input amstex.tex
\long\def\comment#1\endcomment{}
\newcount\sectioncount
\newcount\commoncount\commoncount=1
\newcount\lemmacount\lemmacount=1
\newcount\defcount\defcount=1
\newcount\eqcount\eqcount=1
\newcount\theoremcount\theoremcount=1
\newcount\Ccount\Ccount=1

\def\firstsection#1{\sectioncount=#1 \advance\sectioncount by -1}
\firstsection1
\def\renewcount{\commoncount=1\eqcount=1
}
\renewcount

\def\newsection#1#2\par{\global\advance \sectioncount by 1%
\renewcount%
\specialhead\bigbf{\the\sectioncount}.\lable{#1}\bigbf\ #2
\baselineskip=14pt
\expandafter\xdef\csname Section #1 \endcsname{%
\the\sectioncount}%
\endspecialhead
\relax}
\def\section#1{{\sl Section \csname Section #1 \endcsname}}

\def\newgrphno#1{
\the\sectioncount.\the\commoncount
\expandafter\xdef\csname Grph #1 \endcsname{%
\the\sectioncount.\the\commoncount}%
\global\advance\commoncount by1\relax}
\def\grph#1{{\csname Grph #1 \endcsname}}

\def\newtheorem#1{
\proclaim{\the\sectioncount.\the\commoncount. Theorem}\lable{#1}
\expandafter\xdef\csname Theorem #1 \endcsname{\the\sectioncount%
\the\commoncount}%
\global\advance\commoncount by1\it\relax}
\def\theorem#1{{\sl Theorem \csname Theorem #1 \endcsname}}

%

%

\def\newlemma#1{
\proclaim{\the\sectioncount.\the\commoncount.\lable{#1} Lemma}
\expandafter\xdef\csname Lemma #1 \endcsname{%
\the\sectioncount.\the\commoncount}%
\global\advance\commoncount by1\it\relax}
\def\lemma#1{{\sl Lemma \csname Lemma #1 \endcsname}}
\def\lemmano#1{\csname Lemma #1 \endcsname}

\def\newproposition#1{
\proclaim{\the\sectioncount.\the\commoncount.\lable{#1} Proposition}
\expandafter\xdef\csname Proposition #1 \endcsname{%
\the\sectioncount.\the\commoncount}%
\global\advance\commoncount by1\it\relax}
\def\proposition#1{{\sl Proposition \csname Proposition #1 \endcsname}}
\def\propositionno#1{\csname Proposition #1 \endcsname}

\def\newcorollary#1{
\proclaim{\the\sectioncount.\the\commoncount.\lable{#1} Corollary}
\expandafter\xdef\csname Corollary #1 \endcsname{%
\the\sectioncount.\the\commoncount}%
\global\advance\commoncount by1\it\relax}
\def\corollary#1{{\sl Corollary \csname Corollary #1 \endcsname}}
\def\corollaryno#1{\csname Corollary #1 \endcsname}

\def\newdefinition#1{
{\smallskip\noindent\bf\the\sectioncount.\the\commoncount.\lable{#1}
Definition.}
\expandafter\xdef\csname Definition #1 \endcsname{%
\the\sectioncount.\the\commoncount}%
\global\advance\commoncount by1\relax}
\def\definition#1{{\sl Definition \csname Definition #1 \endcsname}}
\def\definitionno#1{\csname Definition #1 \endcsname}
\let\dfntn=\definition

\def\makeeq#1{
\expandafter\xdef\csname Equation #1 \endcsname{%
\the\sectioncount.\the\eqcount}
\expandafter\xdef\csname Equationlable #1 \endcsname{{#1}}
\global\advance\eqcount by1\relax}
\def\eqqno#1{\csname Equation #1 \endcsname}

\def\eqlable#1{{\rom{[\csname Equationlable #1 \endcsname]}}}

\documentstyle{amsppt}

\def\example#1{{\sl Example \csname Example #1 \endcsname}}
\def\endremark{\smallskip}
\let\definition=\dfntn

\input epsf.tex
\newdimen\xsize
\newdimen\oldbaselineskip
\newdimen\oldlineskiplimit
\xsize=.7\hsize
\def\nolineskip{\oldbaselineskip=\baselineskip\baselineskip=0pt%
\oldlineskiplimit=\lineskiplimit\lineskiplimit=0pt}
\def\restorelineskip{\baselineskip=\oldbaselineskip%
\lineskiplimit=\oldlineskiplimit}
\def\putm[#1][#2]#3{
\vbox to 0pt{\noindent\hskip#1\xsize\lower#2\xsize\hbox{$#3$}\vss}}
%
\def\putt[#1][#2]#3{
\vbox to 0pt{\noindent\hskip#1\xsize\lower#2\xsize%
\vtop{\restorelineskip#3}\vss}}

\def\rom{\ifmmode \fam\rmfam \else\rm \fi}
\let\xpar=\par

\def\Month{\ifcase\month \or January\or February\or March\or April\or May\or
June\or July\or August\or September\or October\or November\or December\fi}

\font\bigbf=cmbx10 scaled 1200

\thinmuskip = 2mu
\medmuskip = 2.5mu plus 1.5mu minus 2.1mu  
\thickmuskip = 4mu plus 6mu
\font\teneusm=eusm10
\font\seveneusm=eusm7
\font\fiveeusm=eusm5
\newfam\eusmfam
\textfont\eusmfam=\teneusm
\scriptfont\eusmfam=\seveneusm
\scriptscriptfont\eusmfam=\fiveeusm
\def\scr#1{{\fam\eusmfam\relax#1}}
\font\tenmib=cmmib10
\font\sevenmib=cmmib7
\font\fivemib=cmmib5
\newfam\mibfam
\textfont\mibfam=\tenmib
\scriptfont\mibfam=\sevenmib
\scriptscriptfont\mibfam=\fivemib
%
\font\tensf=cmss10
\font\sevensf=cmss8 
\font\fivesf=cmr5
\newfam\sffam
\textfont\sffam=\tensf
\scriptfont\sffam=\sevensf
\scriptscriptfont\sffam=\fivesf
\def\sf{\fam\sffam}
\font\mathnine=cmmi9
\font\rmnine=cmr9
\font\cmsynine=cmsy9
\font\cmexnine=cmex10 scaled 913
\def\msmall#1{\hbox{$\displaystyle \font\ninesl=cmsl9
\textfont0=\rmnine \textfont1=\mathnine \textfont2=\cmsynine
\textfont3=\cmexnine \textfont\slfam=\ninesl
{#1}$}}
\hyphenation{Lip-schit-zian Lip-schitz}
\def\cc{{\Bbb C}}

\def\pp{{\Bbb P}}

\def\loc{{\rom loc}}

\def\coker{{\sf Coker}\,}
\def\dim{{\sf dim}\,}

\def\endo{{\sf End}}

\def\exp{{\sf exp}}
\def\ev{{\sf ev}}
\def\sfh{{\sf H}}

\def\id{{\sf Id}}
\def\ind{{\sf ind}}

\def\ker{{\sf Ker}\,}
\def\lim{\mathop{\sf lim}}

\def\ord{{\sf ord}}

\def\pr{{\sf pr}}
\def\rank{{\sf rk}}

\def\epsi{\varepsilon}

\def\barr#1{\overline{#1}}

\def\ddef{\mathrel{{=}\raise0.23pt\hbox{\rm:}}}
\def\deff{\mathrel{\raise0.23pt\hbox{\rm:}{=}}}
\def\ge{\geqslant}

\def\<{\langle}
\def\>{\rangle}
\def\fraction#1/#2{\mathchoice{{\msmall{ #1\over#2}}}%
{{ #1\over #2 }}{{#1/#2}}{{#1/#2}}}

\def\le{\leqslant}

\def\maath{\mathsurround=0pt}
\def\lowminus{\hbox{\vbox to 0.3em{\vss
\hbox{$\maath\mkern-2mu \mathord- \mkern-2mu$}}}}

\def\deform #1#2#3#4#5{#1
\setbox1=\hbox{$\maath\scriptstyle\,\,#2\,$}
\buildrel{\,\,#2}\over{
\leftarrow\mkern-9mu
\hbox to\wd1{\cleaders\lowminus\hfill}}
#3
\setbox1=\hbox{$\maath\scriptstyle\,\,#4\,$}
\buildrel{#4\,\,}\over{
\hbox to\wd1{\cleaders\lowminus\hfill}
\mkern-9mu\rightarrow}
#5}

\def\mapdown#1|#2{\llap{$\vcenter{\hbox{$\scriptstyle #1$}}$}
 {{ \big\downarrow}}
  \rlap{$\vcenter{\hbox{$\scriptstyle #2$}}$}}

\let\hook=\hookrightarrow
\def\trans{\pitchfork}
\let\wh=\widehat
\let\wt=\widetilde
\let\ti=\tilde

\def\state#1. {\medskip\noindent{\bf#1. }}

\def\Chi{\raise 2pt\hbox{$\chi$}}
\let\phI=\phi\let\phi=\varphi\let\varphi=\phI

%
%
%

\def\.{\thinspace}
\def\3{\ss}
\def\isl{\text{\sl i}}

\def\cal#1{{\scr{#1}}}
\def\cala{{\cal A}}
\def\alg(#1){\cala(\barr{#1})}
\def\calc{{\cal C}}

\def\cale{{\cal E}}
\def\calf{{\cal F}}

\def\calh{{\cal H}}
\def\calj{{\cal J}}
\def\call{{\cal L}}
\def\calm{{\cal M}}
\def\caln{{\cal N}}
\def\calo{{\cal O}}
\def\calp{{\cal P}}
\def\calq{{\cal Q}}
\def\cals{{\cal S}}
\def\calt{{\cal T}}

\def\calw{{\cal W}}
\def\calx{{\cal X}}
\def\caly{{\cal Y}}
\def\calz{{\cal Z}}
\def\Hom{{\cal H}{\mkern-3.5mu}{\italic o}{\mkern-3.2mu}{\italic m}}

\def\1{1\!{\rom l}}

\voffset -1.5 true cm

\centerline{\bf DEFORMATIONS OF NONCOMPACT COMPLEX CURVES}
\medskip
\centerline{\bf AND MEROMORPHIC ENVELOPES OF SPHERES }


\def\version{}
\leftheadtext{\hss\vtop{%
\line{\hfil S.\.Ivashkovich\ \ V.\.Shevchishin \hfil
\llap{\version}}%
\vskip 4pt \hrule }\hss}

\rightheadtext{\hss\vtop{%
\line{\rlap{\version}%
\hfil Envelopes of spheres\hfil}%
\vskip 4pt \hrule }\hss}



\topmatter
\date
23/JAN/98
\enddate
\address
Universit\'e de Lille-I, U.F.R\.\newline
des   Math\'ematiques, Villeneuve d'Ascq,
France\newline and\newline
Ruhr-Universit\"at, Bochum, Germany
\endaddress
\email
ivachkov\@gat.univ-lille1.fr \
sewa\@cplx.ruhr-uni-bochum.de
\endemail
\author
S.\,M.~Ivashkovich and V.\,V.~Shevchishin
\endauthor

\abstract
The paper is devoted to the properties of the envelopes of meromorphy
of neighborhoods of symplectically immersed  two-spheres
 in complex K\"ahler surfaces.
The method used to study the envelopes of meromorphy is based
on Gromov's theory of pseudoholomorphic curves.
The exposition includes a construction of a complete family
of holomorphic deformations of a non-compact complex curve
in a complex manifold, parametrized by a finite codimension analytic subset
of  a Banach ball.
The existence of this family is used    to prove a generalization
of Levi's continuity principle, which is
applied  to describe envelopes of meromorphy.

Bibliography: 15 titles.
\endabstract
\endtopmatter

\document

\footnote""{This research was carried out
with the financial support of the `RiP' program
of the Mathematical Institute in Oberwolfach.}

\bigskip\noindent
\bigbf 0. Introduction
\smallskip\rm

In the present paper we study the envelopes of meromorphy
of neighborhoods of  two-spheres immersed
in complex surfaces.

Throughout the paper a {\it complex surface} means
a (Hausdorff) connected complex two-dimensional manifold $X$
countable at infinity.
Let $U$ be a domain in~$X$.
Its envelope of meromorphy $(\widehat U,\pi)$
is the maximal domain over~$X$
satisfying the following conditions:
\roster
\item"(i)"
there exists a holomorphic embedding $i\:U\to\widehat U$
with $\pi\circ i=\id_U$;

\item"(ii)"
each meromorphic function $f$ on $U$
extends 
to a meromorphic function $\widehat f$ on~$\widehat U$,
that is, $\widehat f\circ i=f$.
\endroster

The envelope of meromorphy exists for each domain $U$.
This can be proved, for example,
by applying the Cartan--Thullen method
to the sheaf of meromorphic functions on $X$, see~\cite{1}.

In the sequel we shall restrict ourselves to K\"ahler
complex surfaces, that is, we assume that  $X$
carries a strictly positive closed $(1,1)$-form~$\omega$.

Let $S$ be an oriented real surface.

\smallskip\noindent\bf Definition. \sl A $C^1$-smooth immersion 
$u\:S\to(X,\omega)$
is called {\it symplectic} if $u^*\omega$
does not vanish anywhere on~$S$.

\rm The aim  of the present paper is the following result.

\proclaim{Main theorem}
Let $u\:S^2\to X$ be a symplectic immersion of the two-sphere $S^2$
in a disc-convex K\"ahler surface~$X$
such that $M:=u(S)$ has only positive double points.
Assume that $c_1(X)[M]>0$.
Then the envelope of meromorphy $(\widehat U,\pi)$
of an arbitrary neighborhood $U$ of $M$
contains a rational curve $C$ with $\pi^*c_1(X)[C]>0$.
\endproclaim

The definition of disc-convexity 
is given in \S 4.
At this point, we only observe that all compact manifolds are disc-convex.
As usual, let $c_1(X)$ be the first Chern class of~$X$.

This result is a considerable improvement of 
Theorem 1 in~\cite{2}.
First of all, the `positivity' condition for~$X$ is removed.
Moreover, another condition from~\cite{2} is weakened,
namely, the condition $c_1(X)[M]\ge\delta$.
In the case when $X=\Bbb{CP}^2$ this means that $[M]$
may have arbitrary degree in~$H_2(\Bbb{CP}^2,\Bbb Z)$,
whereas the condition $c_1(X)[M]\ge\delta$
imposes the restriction $\operatorname{deg}[M]\le8$
(cf\. \cite{2}; Corollary~1).

It may be helpful to observe (cf., for example,~\cite{3})
that a compact complex algebraic surface $X$
contains a smooth rational curve~$C$ with~$c_1(X)[C]>0$
only in the following three cases:

\roster
\item"(1)" $c_1(X)[C]=1$\ \ $\Rightarrow$\ \ $[C]^2=-1$\ \
$\Rightarrow$\ \ $C$ is an exceptional curve of the first kind\rom;

\item"(2)" $c_1(X)[C]=2$\ \ $\Rightarrow$\ \ $[C]^2=0$\ \
$\Rightarrow$\ \ $X$ can be blown-down to a ruled surface\rom;

\item"(3)" $c_1(X)[C]\ge3$\ \ $\Rightarrow$\ \ $[C]^2\ge1$\ \
$\Rightarrow$\ \ $X$ is either $\Bbb{CP}^2$,
or a Hirzebruch surface $\Sigma_n$, or a modification of the latter.
\endroster

In the last case $X\backslash C$ is pseudoconcave
and therefore, by a theorem of Grauert~\cite{4}
the envelope $\widehat U\supset C$ coincides with~$X$.

Note also that for 
{\it embedded\/} 
real surfaces (not necessarily spheres) in $\Bbb{CP}^2$
the results of~\cite{2} have been  recently improved
by Nemirovski~\cite{5}, who used a different approach
based on the Seiberg--Witten theory.

The method used in the present paper
to construct the envelopes of meromorphy
was proposed in~\cite{2}.
It is based on Gromov's theory of pseudoholomorphic curves~\cite{6}.
A considerable improvement of 
the results of~\cite{2}
is achieved by the development of the transversality
theory in the moduli space of pseudoholomorphic curves, see~\S\,2.

It should be pointed out that the proof of the continuity
principle in~\cite{2}; Theorem~5.1.3
relies heavily on the existence of a complete holomorphic
family of deformations of a non-compact complex curve
in a complex manifold (\cite{2}; Theorem~6.3.1).
In the present paper we give a complete proof
of this result, see  Theorem~3.4.

\bigskip
\bigbf 1. Moduli space of pseudoholomorphic curves
and the first variation of the $\overline\partial$-equation
\smallskip\rm

In this section we  recall briefly, in the form
convenient for  this paper,
and also complete the results of~\cite{2} concerning certain  
basic ideas of the theory of pseudoholomorphic curves.

Consider the Teichm\"uller space $\Bbb T_g$
of complex structures on a closed real oriented surface
$S$ of genus~$g$.
This is a complex manifold of dimension
$$
\dim_{\Bbb C}\Bbb T_g=\cases
0,& \text{if $g=0$;}\\
1,& \text{if $g=1$;}\\
3g-3,& \text{if $g\ge2$,}
\endcases
$$
uniquely characterized by the following property:
the product $S\times\Bbb T_g$ admits a complex structure
$J_{S\times\Bbb T_g}$
such that
\roster
\item"(i)" the natural projection $\pi|_{\Bbb T}\:S\times\Bbb T_g\to\Bbb T_g$
is holomorphic, and therefore for each $\tau\in\Bbb T_g$
the identification $S\cong S\times\{\tau\}$
defines a complex structure
$J_S(\tau):=J_{S\times\Bbb T}\big|_{S\times\{\tau\}}$
on~$S$;

\item"(ii)"
for each complex structure $J_S$ on $S$
there exists a unique $\tau\in\nomathbreak \Bbb T_g$
and a diffeomorphism $f\:S\to S$ such that $J_S=f^*J_S(\tau)$
(that is, the map $f\:(S,J_S)\to(S,J_S(\tau))$ is holomorphic)
and~$f$ is isotopic to the identity map~$\id_S\:S\to S$.
\endroster

Let $\bold G$ denote the automorphism group of $S\times\Bbb T_g$.
Then
$$
\bold G=\cases
\bold{PGl}(2,\Bbb C)& \text{for $g=0$,}\\
\bold{Sl}(2,\Bbb Z)\ltimes T^2& \text{for $g=1$,}\\
\text{is discrete}& \text{for $g\ge2$.}
\endcases
$$

We shall use the following results on  about $\Bbb T_g$.

If $g=0$, then the surface $S$ is the Riemann sphere $S^2$
and all complex structures on $S^2$ are equivalent to the standard
structure~$S^2\cong\Bbb{CP}^1$.
Hence $\Bbb T_0$ consists of one point
and  $\bold G=\bold{PGl}(2,\Bbb C)$
is the group of biholomorphisms of~$\Bbb{CP}^1$.

If $g=1$, then the surface $S$ is the torus $T^2$
and $\Bbb T_1$ is the upper half-plane
$\Bbb C_+=\{\tau\in\Bbb C:\Im(\tau)>0\}$.
Here the product $S\times\Bbb T_1$
can be identified with the quotient $(\Bbb C\times\Bbb C_+)/\Bbb Z^2$
under the action
$$
\bigl((m,n),\,(z,\tau)\bigr)\in\Bbb Z^2\times\Bbb C\times\Bbb C_+\longmapsto
(m,n)\cdot(z,\tau):=(z+m+n\,\tau,\tau)\in\Bbb C\times\Bbb C_+.
$$
In this case $T^2\equiv\Bbb R^2/\Bbb Z^2$ is the identity component
of ${\roman e}\in\bold G$; in particular, $T^2$~is normal.
The group $\bold G=\bold{Sl}(2,\Bbb Z)\ltimes T^2$
is a semi-direct product
and the action of $T^2$ on $(\Bbb C\times\Bbb C_+)/\Bbb Z^2$
is given by the formula
$$
\multline\qquad
\bigl([t_1,t_2],\,([z],\tau)\bigr)\in T^2\times(\Bbb C\times\Bbb C_+)/\Bbb Z^2
\\
\longmapsto[t_1,t_2]\cdot([z],\tau):=([z+t_1+t_2\tau],\tau)
\in(\Bbb C\times\Bbb C_+)/\Bbb Z^2.
\qquad\endmultline
$$

The action of $\bold G$ on $S\times\Bbb T_g$
is effective for all $g\ge0$
and preserves the fibers of the projection
$\pi\big|_{\Bbb T}\:S\times\Bbb T_g\to\Bbb T_g$.
This gives us an action of $\bold G$ on~$\Bbb T_g$.
Furthermore, for each $\tau\in\Bbb T_g$ and $f\in\bold G$
there exists a unique diffeomorphism
$\widehat f_\tau\:S\to S$ such that
$$
f\cdot(x,\tau)=\bigl(\widehat f_\tau(x),f\cdot\tau\bigr).
\tag 1.1
$$

For each $\tau\in\Bbb T_{\!g}$
we have the natural isomorphisms 
$T_{\tau}\Bbb T_{g}\cong\sfh^1(S,\,\calo_{\tau(TS)})$ and
$T_{\roman e}\bold G\cong\roman H^0(S,\calo_{\tau (TS)})$,
where $\calo_{\tau(TS)}$ denotes the sheaf of sections of $TS$ that are
holomorphic with respect to $J_S(\tau)$.

Below we denote elements of $\Bbb T_g$ by $J_S$
and regard them as the corresponding complex structures on~$S$.

\smallskip
Consider a symplectic manifold $(X,\omega)$
with a fixed almost complex structure $J_{\text{\rom{st}}}$.
Recall that $J_{\text{\rom{st}}}$ is called {\it $\omega$-tamed\/}
(see~\cite{6})
if for each non-zero vector $v\in T_xX$ we have
$\omega(v,J_{\text{\rom{st}}} v)>\nomathbreak 0$.
This  is equivalent to the condition
that the $(1,1)$-component of the form~$\omega$ be positive
with respect to $J_{\text{\rom{st}}}$
:
$\omega^{(1,1)}>0$.
In our applications $\omega$ and $J_{\text{\rom{st}}}$
will define a K\"ahler structure on~$X$.

Let $U$ be a relatively compact subdomain of~$X$
that may coincide with~$X$ if   $X$ is compact.
Let $S$ be a (fixed) compact oriented real surface
of genus $g\ge0$ and let $u_0\:S\to X$ be a non-constant
$C^1$-smooth map with~$u_0(S)\subset U$.

We fix~$p$, $2<p<\infty$,
and consider the Banach manifold $L^{1,p}(S,X)$
of all (continuous) maps~$u\:S\to X$
in the Sobolev class~$L^{1,p}$.
This is a smooth manifold, and its tangent space
$T_uL^{1,p}(S,X)$ at the point $u$
is the Banach space $L^{1,p}(S,u^*TX)$
of $L^{1,p}$-smooth sections of the pulled-back
tangent bundle~$TX$.
Let $\cals_U$ be     the set of  maps~$u$ in $L^{1,p}(S,X)$
that are homotopic to $u_0$ and satisfy the 
condition 
$u(S)\cap U\ne\varnothing$.
We fix $k\ge1$ and denote by $\calj^k_U$
the set of $C^k$-smooth almost complex structures $J$ on~$X$
satisfying the following two conditions:
\roster
\item"(i)" $\{x\in X:J(x)\ne J_{\text{\rom{st}}}(x)\}\Subset U$;

\item"(ii)" $J$ is $\omega$-tamed.
\endroster

The map $\ev\:S\times\cals_U\times\Bbb T\times\calj^k_U\to X$ 
given 
by the formula $\ev(x,u,J_S,J):=u(x)$
defines a bundle $E:=\ev^*(TX)$ over
$S\times\cals_U\times\Bbb T\times\calj^k_U$.
We equip $E$ with the natural complex structure, which is
equal to $J(u(x))$ on each fiber~$E_{(x,u,J_S,J)}\cong T_{u(x)}X$.
Let  $(E_u,J)$ be the restriction of $E$
to $S\times\{(u,J_S,J)\}$; it is isomorphic to~$u^*TX$.

The bundle $E$ with complex structure $J$
defines two complex Banach bundles,
$\wh{\cale }$ and $\wh{\cale }'$,
over the product $\cals_U\times\Bbb T\times\calj^k_U$. They have
the  fibers
$$
\wh{\cale }_{(u,J_S,J)} :=L^{1,p}(S,E_u)  
 \text{\ \ and\ \ } 
\wh{\cale }'_{(u,J_S,J)} :=L^p(S,E_u\otimes\Lambda^{(0,1)}S).
$$
Here $S$ carries the complex structure $J_S$,
$\Lambda^{(0,1)}S$ is the complex line bundle
of $(0,1)$-forms on $S$,
and $\otimes$ is  the complex tensor product
of corresponding complex vector bundles.
Note that $\wh{\cale }$ is the pull-back
of the tangent bundle $TL^{1,p}(S,X)$ with respect to the projection
$(u,J_S,J)\in\cals_U\times\Bbb T\times\calj^k_U
\mapsto u\in\nomathbreak L^{1,p}(S,X)$.

In its turn, the bundle $\wh{\cale }'$ is the range
of the $\overline\partial$-operator on the manifold

\noindent $L^{1,p}(S,X)\times\Bbb T\times\calj^k_U$.
Namely,  the $\overline\partial$-operator
defines a section $\sigma_{\overline\partial}$
of the bundle $\wh{\cale }'$ by the formula
$$
\sigma_{\overline\partial}(u,J_S,J):=\overline\partial_{J_S,J}u:=
\frac12\bigl(du+J\circ du\circ J_S\bigr).
\tag 1.2
$$
If $f$ is another section of $\wh{\cale}'$
(defined, for example, by an explicit geometric construction),
then we can consider the non-homogeneous
$\overline\partial$-equation
$$
\overline\partial_{J_S,J}u=f(u,J_S,J).
$$
In the present paper we consider only the homogeneous case
$f(u,J_S,J)\equiv0$.
We denote the corresponding set of solutions by
$$
\calp :=\bigl\{(u,J_S,J)\in\cals_U\times\Bbb T\times\calj^k_U:
\overline\partial_{J_S,J}u=0\bigr\}.
\tag 1.3
$$
If $u\in L^{1,p}(S,X)$ satisfies the equation
$$
\overline\partial_{J_S,J}u=0
\tag 1.4
$$
with appropriate $J_S\in\Bbb T_g$ and $J\in\calj^k_U$,
then we say that the map $u$ is pseudoholomor\-phic,
or $J$-holomorphic, or $(J_S,J)$-holomorphic,
and we call the image $M:=u(S)$
a~pseudoholomorphic (or  $J$-holomorphic) curve.

The operator $\overline\partial$ is elliptic with  Cauchy--Riemann symbol.
The theory of elliptic partial differential equations
(see, for example,~\cite{7})
shows that $\calp $ is closed in the space
$$
\calx :=\bigl\{(u,J_S,J)\in\cals_U\times\Bbb T\times\calj^k_U:
u\in C^1(S,X)\bigr\}.
$$
Note that $\calx $ is a closed Banach manifold. We set
$$
\multline\qquad
\calx^*:=\bigl\{(u,J_S,J)\in\calx :
\text{there exists a non-empty open subset $V$ of $S$}
\\
\text{such that
$u|_V$ is an embedding and $u(V)\cap u(S\setminus V)=\varnothing$}\bigr\},
\qquad\endmultline
$$
and let $\calp^*:=\calp\cap\calx^*$.
Then $\calx^*$ is open in $\calx$
and $\calp^*$ is open in $\calp$.
We consider the following natural action of the group $\bold G$ on
$\cals_U\times\Bbb T\times\calj^k_U$
by  means of  compositions:
$$
(f,u,J_S,J)\in\bold G\times\cals_U\times\Bbb T\times\calj^k_U
\longmapsto f\cdot(u,J_S,J):=(u\circ\widehat f^{-1},f\cdot J_S,J),
$$
where $\widehat f\:S\to S$ is the diffeomorphism induced
by $f\in\bold G$ and $J_S$, see~(1.1).
The inverse map $(\widehat f)^{-1}$
is introduced to make this action associative:
$(g_1\cdot g_2)\cdot(u,J_S,J)=g_1\cdot\nomathbreak (g_2\cdot(u,J_S,J))$.

The sets $\calx^*$, $\calp$ and $\calp^*$
are invariant with respect to the action of~$\bold G$.
Let $\calm:=\calp^*/\bold G$ be the quotient space
for the action of~$\bold G$
and $\pi_{\calp}\:\calp^*\to\calm$ 
the corresponding projection.
Since $\bold G$ does not act on $\calj^k_U$,
the projection
$\pi_{\calj}\:\calm\to\calj^k_U$
is well-defined.

\proclaim{Lemma 1.1}
The projections $\calx^*\to\calx^*/\bold G$ and
$\pi_{\!\calp\!}\:\calp^*\to\calm$
are principal $\bold G$-bundles.
\endproclaim

\smallskip\noindent\bf Proof. \rm 
First, we consider the case $g\ge2$.
It is known that $\bold G$
is discrete in this case and acts properly discontinuously
on~$\Bbb T_g$ (see, for example,~\cite{8}).
Hence the action of $\bold G$ on $\calx$ is also proper.
Together with the definition of $\calx^*$ this means
that $\bold G$ acts freely on~$\calx^*$.
Thus, $\calx^*\to\calx^*/\bold G$
is an unramified covering.

Now, we turn to the case $g=0$. In this
case $S=S^2$ and $\Bbb T_0=\{J_{\text{\rom{st}}}\}$.
We fix some $(u^0,J_{\text{\rom{st}}},J^0)\in\calx^*$.
Let $y_1$, $y_2$, and $y_3$ be three distinct points
on $S^2$ such that $u^0$ is an embedding
in a neighborhood of each~$y_i$
and, in particular, $du^0$ does not vanish at~$y_i$.
We consider  smooth submanifolds~$Z_i$
of~$X$ of codimension~$2$
that intersect $u^0(S^2)$ transversally
at the points $u^0(x_i)$, respectively.

Let $V\ni(u^0,J_{\text{\rom{st}}},J^0)$
be an open subset of $\calx^*$, 
$W$ its projection onto $\calx^*/\bold G$,
and $\bold G\cdot V:=\{f\cdot(u,J_{\text{\rom{st}}},J):
f\in\bold G,(u,J_{\text{\rom{st}}},J)\in V\}$ 
its $\bold G$-saturation.
We consider the set
$$
{\calz}:=\bigl\{(u,J_{\text{\rom{st}}},J)\in\bold G\cdot V:
u(y_i)\in Z_i\bigr\}.
$$
If $V$ is sufficiently small, then ${\calz}$
is a smooth submanifold of $\bold G\cdot V$, which
 intersects each orbit $\bold G\cdot(u,J_{\text{\rom{st}}},J)$
transversally at a single point.
This defines a $\bold G$-invariant
diffeomorphism $\bold G\cdot V\cong\bold G\times{\calz}$,
and therefore ${\calz}$ is a local slice of
the $\bold G$-action at ~$(u^0,J_{\text{\rom{st}}},J^0)$.

The projection of these local slices of ${\calz}$
onto $\calx^*/\bold G$
defines the structure of a smooth Banach manifold on~$\calx^*/\bold G$
and the structure of a principal $\bold G$-bundle
on the projection $\calx^*\to\calx^*/\bold G$.

\smallskip
The case $g=1$ can be regarded   as a combination
of the previous cases  because for  $g=1$
the group $\bold G$ has the `continuous' part and the `discrete' part,
$T^2$ and  $\bold{Sl}(2,\Bbb Z)$, respectively.
We fix $(u^0,J_S^0,J^0)\in\calx^*$ first
and then a point $y$ on $S$  such that
$u^0$ is an embedding in a neighborhood of~$y$
and, in particular, $du^0$ does not vanish at~$y$.
Let $Z$ be a smooth submanifold of $X$ of codimension~2
that intersects $u^0(S)$ transversally at~$u^0(y)$.
As in the case~$g=0$,
we fix a neighborhood~$V\ni(u^0,J_S^0,J^0)$ in~$X$
and consider the  set
$$
{\calz}:=\bigl\{(u,J_S,J)\in\bold G\cdot V:u(y)\in Z\bigr\},
$$
where
$\bold G\cdot V:=\{\,f\cdot(u,J_{\text{\rom{st}}},J):
f\in\bold G,\ (u,J_{\text{\rom{st}}},J)\in V\}$
is the $\bold G$-saturation of~$V$.
If $V$ is sufficiently small, then ${\calz}$
is a slice of the action of the subgroup~$T^2$ of $\bold G$.
Hence $\calx^*\to\calx^*/T^2$ is a principle $T^2$-bundle.

We consider now the action of
$\bold{Sl}(2,\Bbb Z)=\bold G/T^2$ on~$\calx^*/T^2$.
In the same way as in the case~$g\ge2$,
one can show that $\calx^*/T^2\to\calx^*/\bold G$ is a covering.
Consequently, $\calx^*\to\calx^*/\bold G$
is a principal $\bold G$-bundle.

\smallskip
Considering  the projection $\calp^*\to\calm\equiv\calp^*/\bold G$
we  observe that the natural inclusion
$\calm\hookrightarrow\calx^*/\bold G$
is continuous and closed.
This gives us a $\bold G$-invariant homeomorphism
$\calp^*\cong\calm\,\times_{\calx^*/\bold G}\calx^*$
and therefore a principal $\bold G$-bundle structure on $\calp^*$
with  base~$\calm$.

\smallskip\noindent\bf Remark. \rm 
It can be shown  that if $(u,J_S,J)\in\calp^*$,
then $u$ is an embedding in a neighborhood of each point of $S$,
except for a finitely many    points.
In particular, $J_S$ is uniquely determined by $u$ and~$J$.
In a  similar way, each class $\bold G\cdot(u,J_S,J)\in\calm$
is uniquely determined by $J\in\calj^k_U$
and the pseudoholomorphic curve~$M:=u(S)$.
Using this observation  we denote elements of~$\calm$ by~$(M,J)$.
The motivation for this notation is that
the objects that we shall  study and use
are  pseudoholomorphic curves themselves,
rather than their particular parametrizations.
We hope that this formal inaccuracy
will not lead to a misunderstanding.
\endremark

This result and the construction in the proof of the lemma
enable us to push down $\bold G$-invariant objects from $\calp^*$
to~$\calm$.

For example, there exists a (trivial) bundle over $\calp$
with fiber $S$, total space $\calp\times S$,
natural projection on $\calp$,
complex structure $J_S$ in the fiber over $(u,J_S,J)$,
and map $\ev\:\calp\times S\to X$
given by the formula $\ev(u,J_S,J;y):=u(y)$.

These structures are invariant under the action of $\bold G$
on $\calp\times S$
 extending the initial action of $\bold G$ on $S\times\Bbb T$.
This gives us a $\bold G$-bundle
$\pi_{\!\calc\!}\:{\calc}\to\calm$
with total space ${\calc}:=\calp^*\times_{\bold G}S$
and fiber $S$ and also a map~$\ev\:{\calc}\to X$.
We shall regard 
$\calm$ as the moduli space of all pseudoholomorphic
curves in~$X$ with appropriate topological properties,
$\pi_{\!\calc\!}\:{\calc}\to\calm$
as the corresponding `universal curve' family,  
and the map $\ev\:{\calc}\to X$
as the realization of this `universal curve' in~$X$.
In particular, each fiber
${\calc}_{(M,J)}:=\pi_{\!\calc\!}^{-1}(M,J)$ over $(M,J)\in\calm$
carries a   natural complex structure
$J_{{\calc}_{(M,J)}}=J_S$.

We can define Banach bundles $\cale$ and $\cale'$
over $\calx^*/\bold G\supset\calm$ in a similar way.
To this end we observe that the inclusion
$\calx^*\hook\cals_U\times\calj^k_U$
is continuous and $\bold G$-invariant.
For $(u,J_S,J)\in\cals_U\times {\Bbb T} 
\times \calj^k_U$ and $f\in\bold G$ we
define a   diffeomorphism $\widehat f\:S\to S$ by~(1.1),
so that $f\cdot(u,J_S,J)=(u\circ\widehat f^{-1},f\cdot J_S,J)$.
We define the operator $f_*$ by setting
$$
f_*\:s\in\wh{\cale}_{(u,J_S,J)}\longmapsto
(\widehat f^{-1})^*s\in\wh{\cale}_{f\cdot(u,J_S,J)}.
$$
This gives us 
a natural lifting of the $\bold G$-action
on $\cals_U \times {\Bbb T} 
\times\calj^k_U$
to a $\bold G$-action on~$\wh{\cale}$.
A  lifting of the $\bold G$-action to $\wh{\cale}'$
is defined in the same way.
Since the projection $\calx^*\to\calx^*/\bold G$
admits local $\bold G$-slices,
 there exist Banach bundles
$\cale$ and $\cale'$ over $\calx^*/\bold G$ such that
their liftings to $\calx^*$
are $\bold G$-equivariantly isomorphic
to the bundles $\wh{\cale}$ and $\wh{\cale}'$, respectively.
In~particular, if $(u,J_S,J)\in\calx^*$ represents 
the class~$\bold G\cdot\nomathbreak (u,J_S,J)\in\calx^*/\bold G$,
then we have the following natural isomorphisms:
$$
\align
\cale_{\bold G\cdot(u,J_S,J)}&\cong\wh{\cale}_{(u,J_S,J)}
=L^{1,p}(S,E_u),
\\
\cale'_{\bold G\cdot(u,J_S,J)}&\cong\wh{\cale}'_{(u,J_S,J)}
=L^p(S,E_u\otimes\Lambda^{(0,1)}S).
\endalign
$$

\smallskip
We intend  to study  the problem of deformation 
of a fixed (compact) $J_0$ - holomorphic curve $M_0=u_0(S)$
into a compact complex curve~$M_1$, which is
holomorphic with respect to some given (for instance, integrable)
almost complex structure $J_{\text{\rom{st}}}$ on~$X$,
using the continuation method. 
The idea is to find a suitable homotopy
$h(t)=J_t$, $t\in[0,1]$, of almost complex structures
between $J_0$ and $J_{\text{\rom{st}}}=J_1$
and to construct a continuous deformation
$u_t\:S\to X$ of the map $u_0$ into the  map~$u_1$ such that
 $u_t$ is $J_t$-holomorphic for all~$t\in[0,1]$.
To this end  we shall study the linearization
of the equation $\overline\partial_Ju=0$.

\proclaim{Lemma 1.2}
Let $\calx$ be a Banach manifold,
let $\cale\to\calx$ and $\cale'\to\calx$
be $C^1$-smooth Banach bundles over $\calx$, and let
$\nabla$ and~$\nabla'$ be connections
in $\cale$ and $\cale'$, respectively.       Let
$\sigma$ be  a~(local) $C^1$-section of $\cale$
and let $D\:\cale\to\cale'$ be a $C^1$-smooth bundle homomorphism.
\widestnumber\item{}
\roster
\item"(i)" If $\sigma(x)=0$ for some $x\in\calx$,
then the map $\nabla\sigma_x\:T_x\calx\to\cale_x$
does not depend on the choice of the connection $\nabla$ in~$\cale$.

\item"(ii)" Let $K_x:=\ker(D_x\:\cale_x\to\cale'_x)$ and
$Q_x:=\coker(D_x\:\cale_x\to\cale'_x)$. Fix
the corresponding inclusion $i_x\:K_x\to\cale_x$
and  the projection $p_x\:\cale'_x\to\nomathbreak Q_x$.
Let $\nabla^{\Hom}$ be the connection in $\Hom(\cale,\cale')$
induced by $\nabla$ and~$\nabla'$.
Then the map
$$
p_x\circ(\nabla^{\Hom}D_x)\circ i_x\:T_x\calx\to\Hom(K_x,Q_x)
$$
does not depend on the choice of $\nabla$ and~$\nabla'$.
\endroster
\endproclaim

\smallskip\noindent\bf Remark. \rm 
In view of the results  of the lemma,
we shall use the following notation.
For $\sigma\in\Gamma(\calx,\cale)$,
$D\in\Gamma(\calx,\Hom(\cale,\cale'))$ and $x\in\calx$
as in the lemma,
we shall write  $\nabla\sigma_x\:T_x\calx\to\cale_x$ and
$\overline{\nabla D}\:T_x\calx\times\ker D_x\to\coker D_x$
to denote the corresponding operators,
without specifying the  connections used in their definitions.
\endremark

\smallskip\noindent\bf Proof. \rm 
(i) Let $\wt\nabla$ be another connection in~$\cale$.
Then $\wt\nabla$ has the form $\wt\nabla=\nabla+A$
for some $A\in\Gamma(\calx,\Hom(T\calx,\endo (\cale)))$.
Hence, for $\xi\in T_x\calx$ we have
$\wt\nabla_\xi\sigma-\nabla_\xi\sigma=A(\xi,\sigma(x))=0$.

(ii) In a similar way, let $\wt\nabla'$ be another connection
in $\cale'$, and let $\wt\nabla^{\Hom}$
be the connection in $\Hom(\cale,\cale')$
induced by $\wt\nabla$ and~$\wt\nabla'$.
Then $\wt\nabla'$ also has the form
$\wt\nabla=\nabla+A'$ with some
$A'\in\Gamma(\calx,\Hom(T\calx,\endo(\cale')))$.
Thus, for $\xi\in T_x\calx$ we obtain
$\wt\nabla_\xi^{\Hom}D-\nabla_\xi^{\Hom}D=A'(\xi)\circ D_x-D_x\circ A(\xi)$.
 The result  of the lemma follows now
from the identities
$p_x\circ D_x=0$ and~$D_x\circ i_x=0$.

The space $\calp$ of pseudoholomorphic maps
has been  defined as the zero set  of the section
$\sigma_{\overline\partial}\:
\cals_U\times\Bbb T\times\calj^k_U\to\wh{\cale}'$.
Lemma~1.2 asserts  that for each $(u,J_S,J)\in\calp$
we have a well-defined operator
$$
\nabla\sigma_{\overline\partial}\:
T_uL^{1,p}(S,X)\oplus T_{J_S}\Bbb T_g\oplus T_J\calj^k_U
\longrightarrow\wh{\cale}'_{(u,J_S,J)}.
$$
This operator is the required
linearization of the equation
$\overline\partial_{J_S,J}u=0$.
Covariant differentiation of~(1.2) shows
that it has the following form:
$$
\nabla\sigma_{\overline\partial}\:(v,\dot J_S,\dot J)\longmapsto
D_{(u,J)}v+J\circ du\circ\dot J_S+\dot J\circ du\circ J_S
\in\wh{\cale}'_{(u,J_S,J)}.
\tag 1.5
$$
Here $D=D_{(u,J)}$ is the {\it Gromov operator}
defined on $T_u\cals_U\equiv\cale_{(u,J)}\equiv L^{1,p}(S,u^*TX)$.
To obtain an explicit formula for $D$ we
fix some connections $\nabla^S$ on~$TS$ and $\nabla^X$ on~$TX$.
Note that we can choose $\nabla^S$ to be symmetric
and $J_S$-complex for each fixed~$J_S$, which means
s that $\nabla^SJ_S=0$.
However, in general there exists no connection $\nabla^X$
with these two properties (that is, both symmetric and $J$-complex),
 and therefore we must  choose one of them.
In the present paper we prefer symmetric connections $\nabla^X$,
because this enables us to fix one      connection
for all~$J\in\calj^k_U$.
In the sequel we use the same symbol   $\nabla$
for connections in $TS$ and~$TX$,
and also in other bundles obtained from $TX$ and~$TS$,
such as~$E_u=u^*TX$.
Note that the connections $\nabla^S$ and $\nabla^X$
induce connections in the Banach bundles $\cale$ and~$\cale'$.
This yields the following formula for the Gromov
operator $D_{u,J}\:L^{1,p}(S,E_u)\to L^p(S,E_u\otimes\Lambda^{(0,1)}S)$:
$$
D_{u,J}(v)=\nabla v+\frac12 J\circ\nabla v\circ J_S+
\frac12(\nabla_vJ)\circ du\circ J_S.
\tag 1.6
$$
Explicit calculations have been carried out  in~\cite{2} (see also~\cite{9}
for the case of a $J$-complex connection $\nabla^X$ in~$TX$).
The operator $D=D_{u,J}$ is $\bold G$-equivariant
and therefore induces a bundle homomorphism
$D\:\cale\to\cale'$ over~$\calm$.

The following 
properties of  $D$ will be used below.
Let $(u,J_S,J)\in\calp$.

\proclaim{Lemma 1.3}
The Gromov operator
$D_{u,J}\:L^{1,p}(S,E_u)\to L^p(S,E_u\otimes\Lambda^{(0,1)}S)$
is an $\Bbb R$-linear differential operator
order 1 
with Cauchy--Riemann symbol.
It can be split  into the sum
$D_{u,J}=\overline\partial_{u,J}+R$,
where the $\Bbb C$-linear operator
$\overline\partial_E:=\overline\partial_{u,J}$
and the $\Bbb C$-antilinear operator~$R$
have the following properties\rom:
\roster
\item"(i)" there exists a unique  holomorphic structure
on~$E_u$ such that $\overline\partial_E=\overline\partial_{u,J}$
is the corresponding Cauchy--Riemann operator\rom;

\item"(ii)" if $J$ is $C^k$-smooth, then $R$
is a $C^{k-1}$-smooth $\Bbb C$-antilinear homomorphism
from $E_u$ to $E_u\otimes\Lambda^{(0,1)}S$,
that~is,
$R\in C^{k-1}(S,\overline{\Hom_\Bbb C}(E_u,E_u\otimes\Lambda^{(0,1)}S));$
moreover, for $x\in S$, $v\in E_x=T_{u(x)}X$ and
$\xi\in T_xS$ we have $R(v,\xi)=N_J(v,du(\xi))$,
where $N_J$ is the Nijenhuis torsion tensor 
of the almost complex structure~$J$ \rom(see~\cite{10}\rom)\rom;

\item"(iii)" the a priori continuous homomorphism
$du\:TS\to E_u$ is in fact holomorphic,
that~is, it satisfies the relation
$du\circ\overline\partial_{TS}=\overline\partial_E\circ du;$
 furthermore, $R\circ du=0$.
\endroster
\endproclaim

\smallskip\noindent\bf Proof. \rm 
Explicit formula~(1.4) shows that
the Gromov operator $D$
is a first-order differential operator
 with  Cauchy--Riemann symbol.
Since this symbol is $\Bbb C$-linear,
 the $\Bbb C$-antilinear part~$R$ of~$D$
has order~$0$, that is, $R$ is a bundle homomorphism.
The proof of the formula
$R=\nomathbreak N_J\circ\nomathbreak du$ can be found in~\cite{2}; Lemma~2.2.1.

The $\Bbb C$-linear part $\overline\partial_{u,J}$ of~$D$
is of order~$1$ and has the  Cauchy--Riemann symbol.
The fact that there exists a holomorphic structure on~$E_u$
with Cauchy--Riemann operator~$\overline\partial_{u,J}$
is standard.

Let $\xi\in L^{1,p}(S,TS)$, and let $\Phi_t\:S\to S$
be the one-parametric diffeomorphism group
generated by~$\xi$.
We set $u_t:=u\circ\Phi_{-t}$ and $J_S(t):=\Phi_{-t}^*J_S$.
Differentiating the equality
$\overline\partial_{J_S(t),J}u_t=0$ at~$t=0$
and using the identities
$\dfrac{dJ_S(t)}{dt}\bigg|_{t=0}=L_\xi J=2J\overline\partial\xi$ and
$\dfrac{du_t}{dt}\bigg|_{t=0}=du(\xi)$
we obtain the following formula:
$$
du\circ\overline\partial_{TS}=D_{u,J}\circ du.
\tag 1.7
$$
Part (iii) of the lemma follows now
by taking the $\Bbb C$-linear and the $\Bbb C$-antilinear
parts of~(1.7).

For a detailed proof of the lemma we refer to~\cite{2};~\S\,2.

\proclaim{Corollary 1.4}
\rom{(i)} The set of critical points
of a $J$-holomorphic map $u\:(S,J_S)\to(X,J)$
is discrete in~$S$, provided that~$J\in C^1(X)$\rom;

\rom{ii)}
the order of zero of $du$ at a point~$a$
is well defined for each~$a\in S$.
\endproclaim

\smallskip\noindent\bf Remark.  \rm 
Part (i) of Corollary 1.4 was proved by McDuff for~$J\in C^\infty$
and by Sikorav for~$J\in C^1$.
\endremark

Another consequence of Lemma 1.3 is the following
exact sequence of {\it coherent\/} sheaves:
$$
0\longrightarrow{\calo}(TS)\overset du\to\longrightarrow
{\calo}(E_u)\longrightarrow{\caln_u}\longrightarrow0.
\tag 1.8
$$
Here $\caln$ is  the quotient sheaf 
$\calo(E)/du(\calo(TS))$.
This sheaf splits into the direct sum
$\caln_u={\calo}(N_u)\oplus\caln_u^{\text{\rom{sing}}}$,
where $N_u$ is a holomorphic vector bundle
and $\caln_n^{\text{\rom{sing}}}=\bigoplus_{i=1}^P\Bbb C_{a_i}^{n_i}$.
By $\Bbb C_{a_i}^{n_i}$ we mean   the sheaf concentrated
at a critical point~$a_i\in S$ of the map~$du$
with stalk $\Bbb C^{n_i}$, where~$n_i=\ord_{a_i}du$
is the order of the zero of~$du$ at~$a_i$.
We call $\caln_u$ the {\it normal sheaf\/} of $M=u(S)$, 
$N_u$  the {\it normal bundle\/} of~$M$,
and $\caln_u^{\text{\rom{sing}}}$ 
the {\it ramification sheaf\/} of~$M$. 

Let $\roman Z_{du}$ be     the divisor
$\sum_{i=1}^Pn_i[a_i]$,
and let $\call(\roman Z_{du})$ be the sheaf
of meromorphic functions on~$S$
with poles of order at most~$n_i$ at~ $a_i$.
The corresponding linear bundle will be also
denoted by $\call(\roman Z_{du})$.
Then~(1.8) gives rise to the following exact sequence
of holomorphic vector bundles:
$$
0\longrightarrow TS\otimes\call(\roman Z_{du})
\overset du\to\longrightarrow E\longrightarrow N_u\longrightarrow0.
\tag 1.9
$$

Let $L^p_{(0,1)}(S,E)$ be     the space of
$L^p$-integrable $(0,1)$-forms on $S$ with coefficients in~$E$.
Then, by Lemma~1.3 and equality~(1.7)
we obtain the following commutative diagram
with exact rows:
$$
\minCDarrowwidth{5mm}
\CD 
0@>>>L^{1,p}(S,TS\otimes\call(\roman Z_{du}))@>{du}>>
L^{1,p}(S,E)@>\pr>>L^{1,p}(S,N_u)@>>>0
\\
@. @VV\overline\partial_SV @VV D_{u,J} V @VVV @.
\\
0@>>> L^p_{(0,1)}(S,TS\otimes\call(\roman Z_{du}))@>{du}>>
L^p_{(0,1)}(S,E)@>>>L^p_{(0,1)}(S,N_u)@>>>0.
\endCD
\tag 1.10
$$

This gives us the operator
$$
D_{u,J}^N\:L^{1,p}(S,N_u)\to L^p_{(0,1)}(S,N_u)
$$
of the form $D_{u,J}^N=\overline\partial_N+R_N$.
As usual, $\overline\partial_N$ is
the $\overline\partial$-operator in~$N_u$
and $R_N\in C^0(S,\Hom_\Bbb R(N_u,\Lambda^{0,1}\otimes N_u))$.

\definition{Definition 1.1}
Let $E$ be a holomorphic vector bundle
over a compact Riemann surface~$(S,J_S)$,
and let $D\:L^{1,p}(S,E)\to L^p_{(0,1)}(S,E)$
be an operator of the form $D=\overline\partial_E+R$,
where $R\in L^p\bigl(S,\Hom_\Bbb R(E,\Lambda^{0,1}S\otimes E)\bigr)$,
$2<p<\infty$.
Then we set $\roman H^0_D(S,E):=\ker D$ and
$\roman H^1_D(S,E):=\coker D$.

\smallskip\noindent\bf Remark. \rm 
It is not difficult to show that for $E$, $R\in L^p$, $2<p<\infty$, and
$D=\nomathbreak \overline\partial+\nomathbreak R$ as above,
the spaces $\roman H^i_D(S,E)$
can be defined as the kernel and the cokernel of the
 operator
$$
\overline\partial+R\:L^{1,q}(S,E)\to L^q_{(0,1)}(S,E)
$$
with arbitrary $q$, $1<q\le p$.
Consequently, the $\roman H^i_D(S,E)$ do not depend
on the choice of the function spaces.
Furthermore, the operator $D$ is Fredholm, of index
$$
\align
\ind_\Bbb R(D):=&\dim_{\Bbb R}(\ker(D))-\dim_{\Bbb R}(\coker(D))
\\
=&\ind_\Bbb R(\overline\partial)
=2(c_1(E)[S]+(1-g)\cdot\rank_\Bbb C(E)).
\endalign
$$
The proof of these facts is given in~\cite{2}.
\endremark

By the `snake lemma' from homological algebra 
the diagram~(1.10) gives us the following long
exact sequence of $D$-cohomology:
$$
\minCDarrowwidth{7mm}
\CD
0@>>>\roman H^0(S,TS\otimes\call(\roman Z_{du}))@>>>
\roman H^0_D(S,E)@>>>\roman H^0_D(S,N_u)
\\
@>\delta>>\roman H^1(S,TS\otimes\call(\roman Z_{du}))@>>>
\roman H^1_D(S,E)@>>>\roman H^1_D(S,N_u)@>>>0.
\endCD
\tag 1.11
$$

The proof of the next result  can be found in~\cite{2}; Lemma~2.3.2.

\proclaim{Lemma 1.5 \rm (Serre duality for 
$D$-cohomology)}
Let $E$ be a holomorphic vector bundle over a compact
Riemann surface~$(S,J_S)$,
and let $D\:L^{1,p}(S,E)\to L^p_{(0,1)}(S,E)$
be an operator of the form~$D=\overline\partial+R$,
where $R\in L^p\bigl(S,\Hom_\Bbb R(E,\Lambda^{0,1}S\otimes E)\bigr)$,
$2<p<\infty$.
Let $K_S:=\Lambda^{1,0}S$ be   the canonical bundle of~$(S,J_S)$.
Then there exists a naturally defined operator
$$
D^*=\overline\partial-R^*\:L^{1,p}(S,E^*\otimes K_S)
\to L^p_{(0,1)}(S,E^*\otimes K_S)
$$
with
$R^*\in L^p\bigl(S,\Hom(E^*\otimes K_S,\Lambda^{0,1}S\otimes E^*\otimes K_S)
\bigr)$
and there exist natural isomorphisms
$$
\align
\roman H^0_D(S,E)^*&\cong\roman H^1_{D^*}(S,E^*\otimes K_S),
\\
\roman H^1_D(S,E)^*&\cong\roman H^0_{D^*}(S,E^*\otimes K_S).
\endalign
$$
If, furthermore, $R$ is $\Bbb C$-antilinear,
then $R^*$ is also $\Bbb C$-antilinear.
\endproclaim

\smallskip\noindent\bf Remark. \rm 
The proof of Lemma~1.5 uses the following identity, which
that may be taken as the definition of~$D^*$.
For $\xi\in L^{1,p}(S,E)$ and
$\eta\in\nomathbreak L^{1,q}(S,E^*\otimes\nomathbreak K_S)$
we have the equalities
$$
\Re\int_S\langle\overline\partial\xi+R\xi,\eta\rangle=
\Re\int_S\overline\partial\langle\xi,\eta\rangle+
\Re\int_S\langle\xi,-(\overline\partial-R^*)\eta\rangle=
\Re\int_S\langle\xi,-D^*\eta\rangle.
\tag 1.12
$$
Hence the natural pairing
$\langle\,\cdot\,{,}\,\cdot\,\rangle\:
L^p_{(0,1)}(S,E)\times\roman H^0_{D^*}(S,E^*\otimes K_S)\to\Bbb R$,
$\langle\varphi,\eta\rangle:=\Re\dsize\int_S\langle\varphi,\eta\rangle$
vanishes on the range of $D$
and induces the pairing
$$
\langle\,\cdot\,{,}\,\cdot\,\rangle\:
\roman H^1_D(S,E)\times\roman H^0_{D^*}(S,E^*\otimes K_S)\to\Bbb R.
$$
The latter       is the required duality.
\endremark

\proclaim{Corollary 1.6 \rm \cite{6}, \cite{11}
(Vanishing theorems for $D$-cohomology)}
Let $S$ be a Riemann surface of genus~$g$
with  complex structure~$J_S$.
Let $L$ be a holomorphic {\it line\/} bundle over~$S$
with  operator $D\:L^{1,p}(S,L)\to L^p_{(0,1)}(S,L)$
of the form $D=\overline\partial+R$, where
$R\in L^p(S,\Hom_{\Bbb R}(L,\Lambda^{0,1}S\otimes L))$,
$p>2$.

If $c_1(L)<0$, then $\roman H^0_D(S,L)=0$,
and if $c_1(L)>2g-2$, then $\roman H^1_D(S,L)=\nomathbreak 0$.

In particular, if $g=0$, then either $\dim_{\Bbb R}(\roman H^1_D(S,L))=0$
or $\dim_{\Bbb R}(\roman H^0_D(S,L))=0$ and
$\dim_{\Bbb R}(\roman H^1_D(S,L))$ is positive and even.
\endproclaim

\smallskip\noindent\bf Proof. \rm 
Let $\xi$ be a non-trivial $L^{1,p}$-section of~$L$
that satisfies the condition~$D\xi=0$.
By Lemma~3.1.1 in  ~\cite{2} the section~$\xi$
has only  finitely many       zeros $a_i\in S$, which
have {\it positive\/} multiplicities~$\mu_i$,
so that~$c_1(L)=\sum\mu_i\ge0$.
Consequently, $\roman H^0_D(S,\,L)$ is trivial for ~$c_1(L)<0$.
The triviality of~$\roman H^1_D$
follows in a similar way by the Serre duality of Lemma~1.5.

In particular, if $g=0$, then either $c_1(L)<0$, or $c_1(L)>-2$.
Hence, either $\roman H^0_D(S,L)$ or $\roman H^1_D(S,\,L)$ is trivial.
The fact that $\dim_{\Bbb R}(\roman H^1_D(S,L))$ is even
follows  from the equality
$\ind_\Bbb R(D)=\ind_\Bbb R(\overline\partial)
=2\ind_\Bbb C(\overline\partial)$,
which is a consequence of the index theorem.

Note that, by Lemma~1.3
we obtain the following commutative diagram
with exact rows:
$$
\minCDarrowwidth{4mm}
\CD
0@>>>L^{1,p}(S,TS)@>{du}>>L^{1,p}(S,E)@>{\overline\pr}>>
L^{1,p}(S,E)/du\bigl(L^{1,p}(S,TS)\bigr)@>>>0
\\
@. @VV \overline\partial_S V @VV D V @VV\overline D V @.
\\
0@>>>L^p_{(0,1)}(S,TS)@>du>>L^p_{(0,1)}(S,E)@>{\overline\pr}>>
L^p_{(0,1)}(S,E)/du\bigl(L^p_{(0,1)}(S,TS)\bigr)@>>>0,
\endCD
\hskip-5mm
\tag 1.13
$$
where $\overline D$ is induced by~$D\equiv D_{u,J}$.

\proclaim{Lemma 1.7}
For the operator $\overline D$ we have the natural isomorphisms
$\ker\overline D=\roman H^0_{D_N}(S,N_{u})\oplus
\roman H^0(S,\caln_u^{\text{\rom{sing}}})$
and $\coker\overline D=\roman H^1_{D_N}(S,N_u)$.
\endproclaim

For the proof see~\cite{2}; Lemma 2.4.1.

\proclaim{Corollary 1.8}
The short exact sequence of sheaves~\thetag{1.8}
induces the following long exact sequence
of $D$-cohomology\rom:
$$
\minCDarrowwidth{5mm}
\CD
0@>>>\roman H^0(S,TS)@>du>>\roman H^0_D(S,E)@>>>
\roman H^0_D(S,N_u)\oplus\roman H^0(S,\caln_u^{\text{\rom{sing}}})
\\
@>\delta>>\roman H^1(S,TS)@>du>>\roman H^1_D(S,E)@>>>\roman H^1_D(S,N_u)@>>>0.
\endCD
\tag 1.14
$$
\endproclaim

\head
\S\,2. Transversality property of the moduli space
\endhead

To deform a pseudoholomorphic curve $M_t$
along a given path of almost complex structures~$J_t$,
it is useful to know at which points $(u,J_S,J)$
the set of pseudoholomorphic maps~$\calp$
is a Banach manifold.
Note  that, by definition, $\calp$
is the intersection of the zero section
and the section $\sigma_{\overline\partial}$
of the bundle 
$\wh{\cale}'$ over $\cals\times\Bbb T\times\calj^k_U$.
Hence the problem reduces to the question
of the transversality of these two sections.

\definition{Definition 2.1}
Let $\calx$, $\caly $, and ${\calz}$
be Banach manifolds,
and let $f\:\caly \to\calx$ and $g\:{\calz}\to\calx$
be $C^k$-smooth maps,~$k\ge1$.
We define the {\it fibered product\/} 
$\caly \times_{\!\calx\!}{\calz}$
by setting
$\caly \times_{\!\calx\!}{\calz}
:=\{(y,z)\in\caly \times{\calz}:f(y)=g(z)\}$.
The map $f$ is said to be {\it transversal\/} to~$g$
at a point $(y,z)\in\caly \times_{\!\calx\!}{\calz}$,
where $x:=f(y)=g(z)$, and $(y,z)$ is called a
{\it  transversal point} if the map
$df_y\oplus-dg_z\:T_y\caly \oplus T_z{\calz}\to T_x\calx$
is surjective and its kernel admits a closed complementary subspace.
We denote the set of transversal points in
$\caly \times_{\!\calx\!}{\calz}$
by $\caly \times^\trans_{\!\calx\!}{\calz}$,
where $\trans$ indicates  
the transversality condition 

In particular, if the map $g:{\calz}\to\calx$
is an embedding, then the fibered product 
$\caly \times_{\!\calx\!}{\calz}$
is just the inverse image $f^{-1}{\calz}$
of the set ${\calz}\subset\calx$,
and therefore each point $(y,z)\in\caly \times_{\!\calx\!}{\calz}$
is completely defined by~$y\in\caly $.
In this case we say that
$f\:\caly \to\calx$ is {\it transversal\/} to~${\calz}$
at~$y\in\caly $ if $(y,f(y))$ is a transversal point in
$\caly \times_{\!\calx\!}{\calz}\cong f^{-1}{\calz}$.

\proclaim{Lemma 2.1}
The set $\caly \times^\trans_{\!\calx\!}{\calz}$
is open in $\caly \times_{\!\calx\!}{\calz}$.
It is a $C^k$-smooth Banach manifold
with tangent space
$$
T_{(y,z)}\caly \times^\trans_{\!\calx\!}{\calz}=
\ker(df_y\oplus-dg_z\:T_y\caly \oplus T_z{\calz}\to T_x\calx).
$$
\endproclaim

\smallskip\noindent\bf Proof. \rm 
We fix some $w_0:=(y_0,z_0)
\in\caly \times^\trans_{\!\calx\!}{\calz}$
and set 
$$ 
K_0:=\ker(df_{y_0}\oplus-dg_{z_0}\:
T_{y_0}\caly \oplus T_{z_0}{\calz}\to T_x\calx).
$$
Let $Q_0$ be a closed complement of ~$K_0$.
Then the map $df_{y_0}\oplus-dg_{z_0}\:Q_0\to T_x\calx$
is an isomorphism.

Our choice of $Q_0$ ensures    that there exist a neighborhood
$V\subset\caly \times{\calz}$ of~$(y_0,z_0)$
and $C^k$-maps $w'\:V\to K_0$ and $w''\:V\to Q_0$
such that $dw'_{w_0}$ (respectively, $dw''_{w_0}$)
is the projection of $T_{y_0}\caly \oplus T_{z_0}{\calz}$ onto~$K_0$
(respectively, onto $Q_0$).
Hence $(w',w'')$ will are local variables   in a (smaller)
neighborhood $V_1\subset\caly \times{\calz}$ of~$w_0=(y_0,z_0)$.
The lemma follows now by the implicit function theorem
applied to the equation $f(y)=g(z)$
in the new coordinates~$(w',w'')$.

It is easy to see that the set $\calp$
is in fact the fibered product 
of the Banach manifold $\cals_U\times\Bbb T\times\calj^k_U$
by   itself with respect to the maps
$\sigma_0$ and $\sigma_{\overline\partial}$ into~$\cale'$.
By Lemma~2.1, $\calp$ is      a Banach manifold
at points $(u,J_S,J)\in\calp$
where $\sigma_{\overline\partial}$ is transversal to~$\sigma_0$.
However, at each point $(u,J_S,J;0)$
on the zero section $\sigma_0$ of the bundle $\wh{\cale}'$
we have the following natural decomposition:
$$
T_{(u,J_S,J;0)}\wh{\cale}'=
d\sigma_0\bigl(T_{(u,J_S,J)}
(\cals_U\times\calj_S\times\calj_U^k)\bigr)
\oplus\wh{\cale}'_{(u,J_S,J)},
$$
where the first component is the tangent space
of the zero section of $\wh{\cale}'$
and the second is the tangent space of the fiber $\wh{\cale}'_{(u,J_S,J)}$.
Let $p_2$ be     the projection onto the second component.
Then the transversality of $\sigma_{\overline\partial}$ and $\sigma_0$
in~$(u,J_S,J;0)$
is equivalent to the surjectivity of the map
$p_2\circ d\sigma_{\overline\partial}\:
T_{(u,J_S,J)}(\cals_U\times\calj_U^k)\to\cale'_{(u,J_S,J)}$.
However, by Lemma~1.2
this map is the linearization of $\sigma_{\overline\partial}$ at~$(u,J_S,J)$
and, therefore, has the form~(1.3).

Thus, the transversality of $\sigma_{\overline\partial}$ and
$\sigma_0$ at a point $(u,J_S,J)\in\calp$
is equivalent to the surjectivity of the following operator:
$$
\align
&\nabla\sigma_{\overline\partial}\:
T_uL^{1,p}(S,X)\oplus T_{J_S}\Bbb T_g\oplus T_J\calj^k_U
\longrightarrow\wh{\cale}'_{(u,J_S,J)},
\\
&\nabla\sigma_{\overline\partial}\:(v,\dot J_S,\dot J)\longmapsto
D_{(u,J)}v+J\circ du\circ\dot J_S+\dot J\circ du\circ J_S.
\endalign
$$
By Definition 1.1 the quotient of
$\wh{\cale}'_{(u,J_S,J)}=L^p_{(0,1)}(S,E_u)$
by the image $D_{u,J}$ is equal to~$\roman H^1_D(S,E_u)$.
The induced map
$\dot J_S\in T_{J_S}\Bbb T_g\mapsto
[J\circ du\circ\dot J_S]\in\roman H^1_D(S,E_u)$
is also easy to describe.
It follows from~(1.5) and from Corollary~1.8
that the image of $T_{J_S}\Bbb T_g$ under this map coincides with the range
of the homomorphism
$du\circ J_S\:\roman H^1(S,TS)\cong T_{J_S}\Bbb T_g\to\roman H^1_D(S,E_u)$,
and therefore its cokernel is~$\roman H^1_D(S,N_u)$.

Thus, it remains to understand 
the image of $T_J\calj^k_U$ in~$\roman H^1_D(S,N_u)$.
For $(u,J_S,J)\in\calp$ we define the map
$\Psi=\Psi_{(u,J)}\:T_J\calj^k_U\to\wh{\cale}'_{(u,J_S,J)}$
by the formula $\Psi_{(u,J)}(\dot J):=\dot J\circ du\circ J_S$.
Let
$\overline\Psi=\overline\Psi_{(u,J)}\:T_J\calj^k_U\to\roman H^1_D(S,N_u)$
be induced by~$\Psi$.
Recall that if $(u,J_S,J)\in\calp$,
then $J_S$ is uniquely determined by~$u$ and~$J$.

\proclaim{Lemma 2.2 \rm (infinitesimal transversality)}
The operator
$\overline\Psi\:T_J\calj^k_U\to\roman H^1_D(S,N_u)$
is surjective for each $(u,J_S,J)\in\calp^*$.
\endproclaim

\smallskip\noindent\bf Proof. \rm 
It is proved in \cite{2} (see also~\cite{9})
that if $(u,J_S,J)\in\calp^*$,
then the map $u$ is an embedding in a neighborhood
of each point $x\in S$, except for finitely many points.
Hence there exists a non-empty open subset $V\subset S$
such that $u(V)\subset U$ and $u\big|_V$ is an embedding.

By Lemma 1.5 we have the isomorphism
$\roman H^0_D(S,N_u^*\otimes K_S)\equiv\roman H^1_D(S,N_u)^*$.
The fact that an  operator of the form
$D=\overline\partial+R\:L^{1,p}(S,E)\to L^p_{(0,1)}(S,E)$
on a compact Riemann surface $(S,J_S)$ is Fredholm
shows that there exists a finite basis $\xi_1,\dots,\xi_l$
of the space~$\roman H^0_D(S,N_u^*\otimes K_S)$.
By Lemma~3.1.1 in  ~\cite{2} each  $\xi\in\roman H^0_D(S,N_u^*\otimes K_S)$
vanishes at no more than $c_1(N_u^*\otimes K_S)[S]$ points in~$S$
(cf. the proof of Corollary~1.6).
Hence  there exist sections
$\psi_i\in C^k_c(V,N\otimes\Lambda^{0,1})$, $i=1,\dots,l$,
that make up  an $\Bbb R$-basis of the space $\roman H^1_D(S,N)$.

We consider an arbitrary $\psi_i\in C^k_c(V,N\otimes\Lambda^{0,1})$.
This is a $\Bbb C$-antilinear $C^k$-smooth homomorphism
from $TS\big|_V$ into~$N\big|_V$ that vanishes outside a compact subset of~$V$.
Since $u\big|_V$ is a $C^k$-smooth embedding
and $u(V)\subset U$,  $\psi_i$ can be represented
as the composition  
$\psi_i=\pr_N\circ\dot J\circ du\circ J_S$,
where $\dot J$ is a $J$-antilinear $C^k$-smooth endomorphism
of the bundle~$TX$  vanishing outside a compact subset of~$U$.
Hence $\dot J\in T_J\calj^k_U$ and $\overline\Psi\dot J=\psi_i$.

\proclaim{Corollary 2.3} Both
$\calm$ and $\calp^*$ are $C^k$-smooth
Banach manifolds,
and the map $\pi_{\!\calj\!}\:\calm\to\calj^k_U$
is Fredholm.
For each $(M,J)$ in $\calm$ with $M=u(S)$
we have the following natural isomorphisms\rom:
$$
\align
\ker(d\pi_{\!\calj\!}\:T_{(M,J)}\calm\to T_J\calj^k_U)
&\cong\roman H^0_D(S,\caln_M),
\\
\coker(d\pi_{\!\calj\!}\:T_{(M,J)}\calm\to T_J\calj^k_U)
&\cong\roman H^1_D(S,\caln_M),
\endalign
$$
where $\caln_M=\calo(N_u)\oplus\caln^{\text{\rom{sing}}}_u$
is the normal sheaf of~$M$,
$\roman H^0_D(S,\caln_M)$ is   the sum
$\roman H^0_D(S,N_u)\oplus\roman H^0(S,\caln_u^{\text{\rom{sing}}})$,
and $\roman H^1_D(S,\caln_M)$ means  ~$\roman H^1_D(S,N_u)$.
The index of the projection $\pi_{\!\calj\!}$
is described  by the formula
$$
\ind_\Bbb R(\pi_{\!\calj\!})=\ind_\Bbb R(\caln_M)
:=\dim_{\Bbb R}\roman H^0_D(S,\caln_M)-\roman H^1_D(S,\caln_M)
$$
and is equal to $2(c_1(X)[M]+(n-3)(1-g))$, $n:=\dim_{\Bbb C}X$.
\endproclaim

\smallskip\noindent\bf Proof. \rm 
It is easy to see that the section
$\sigma_{\overline\partial}$ is $C^k$-smooth for  $J\in\calj^k_U$.
Hence the assertion about $\calp^*$ follows from Lemmas~2.1~and~2.2.
Moreover, $\calp^*$ is a $C^k$-smooth submanifold of~$\calx^*$.
It is not difficult to show that the slices
of the $\bold G$-action on $\calx^*$,
constructed in the proof of Lemma~1.1,
are $C^\infty$-smooth.
Consequently, the projection $\pi_{\!\calp\!}:\calp^*\to\calm$
induces the structure of a $C^k$-manifold in~$\calm$.

We consider now the projection  $\pi\:\calp^*\to\calj^k_U$.
The tangent space $T_{(u,J_S,J)}\calp^*$
consists of the triples $(v,\dot J_S,\dot J)$
that satisfy the condition
$$
D_{u,J}v+\frac12 J\circ du\circ\dot J_S+\frac12\dot J\circ du\circ J_S=0,
$$
and the differential $d\pi\:T_{(u,J_S,J)}\calp^*\to T_J\calj^k_U$
has the form
$(v,\dot J_S,\dot J)\in
T_{(u,J_S,J)}\calp^*\mapsto\dot J\in T_J\calj^k_U$.

Hence the kernel $\ker(d\pi)\subset T_{(u,J_S,J)}\calp^*$
is parametrized by the solutions of the equation
$$
D_{u,J}v+\frac12 J\circ du\circ\dot J_S=0,
$$
where $v\in\cale_{(u,J_S,J)}$ and $\dot J_S\in T_{J_S}\Bbb T_g$.
Since the map $\pi_{\!\calp\!}\:\calp^*\to\calm$
is a principle $\bold G$-bundle,
$\ker(d\pi_{\!\calj\!}\:T_{(M,J)}\calm\to T_J\calj^k_U)$
is the quotient of~$\ker(d\pi)$
by the tangent space of the fiber $\bold G\cdot(u,J_S,J)$.
Since $T_{\roman e}\bold G\equiv\roman H^0(S,TS)$
and $\bold G_0$ acts trivially on~$\Bbb T_g$ and on~$\calj^k_U$,
 the tangent space of the fiber $\bold G\cdot(u,J_S,J)$
in $(u,J_S,J)$ consists of the vectors of the form~$(v,0,0)$,
$v\in du(\roman H^0(S,TS))$.
>From the equalities
$$
\gather
\roman H^0(S,TS)=\ker\bigl(\overline\partial_{TS}\:L^{1,p}(S,TS)\to
L^p(S,TS\otimes\Lambda^{(0,1)}S)\bigr),
\\
T_{J_S}\Bbb T_g\cong\roman H^1(S,TS)=\coker(\overline\partial_{TS}),
\\
du\circ\overline\partial_{TS}=D_{(u,J)}\circ du
\endgather
$$
we see that
$\ker(d\pi_{\!\calj\!})$ is isomorphic to the quotient
$$
\multline\qquad
\bigl\{v\in L^{1,p}(S,E_u):Dv=du(\varphi)\text{ for some}
\\
\varphi\in L^p(S,TS\otimes\Lambda^{(0,1)}S)\bigr\}
\big/du\bigl(L^{1,p}(S,TS)\bigr).
\qquad\endmultline
$$
Hence $\ker(d\pi_{\!\calj\!}\:
T_{(M,J)}\calm\to T_J\calj^k_U)\cong\roman H^0_D(M,\caln_M)$
by Lemma~1.7.
In~particular, $\ker(d\pi_{\!\calj\!})$ is finite dimensional.

In a similar way, the image $d\pi_{\!\calj\!}$
consists of  $\dot J$ such that the equation
$$
D_{u,J}v+\frac12 J\circ du\circ\dot J_S+\frac12\dot J\circ du\circ J_S=0
$$
has a solution~$(v,\dot J_S)$.
Consequently,
$$
\Im(d\pi_{\!\calj\!})=\ker\overline\Psi \quad\text{and}\quad
\coker(d\pi_{\!\calj\!})\cong\roman H^1_D(S,N_u).
$$
Hence $d\pi_{\!\calj\!}$ is Fredholm, and therefore
 $\pi\:\calp^*\to\calj^k_U$ is also Fredholm.

Corollary 1.8 yields
the equality $\ind_\Bbb R(\caln_M)=\ind_\Bbb R(E_u)-\ind_\Bbb R(TS)$.
Using the Riemann--Roch theorem together with
the equalities $c_1(E)=c_1(X)[M]$ and $c_1(TS)$ $=2-2g$
we obtain the required formula:
$$
\ind_\Bbb R(\caln)=2\bigl(c_1(X)[M]+n(1-g)-(3-3g)\bigr)
=2\bigl(c_1(X)[M]+(n-3)(1-g)\bigr).
$$

Before stating other  results, let us introduce  further  notation.

\definition{Definition 2.2}
Let $Y$ be a $C^k$-smooth finite-dimensional manifold
that may have a non-empty $C^k$-smooth boundary~$\partial Y$,
and let $h\:Y\to\calj^k_U$ be a $C^k$-smooth map.
We define the {\it relative moduli space\/} as follows:
$$
\calm_h:=Y\times_{\calj^k_U}\calm\cong
\bigl\{(u,J_S,y)\in\cals_U\times\Bbb T_g\times Y:
(u,J_S,h(y))\in\calp^*\bigr\}/\bold G,
$$
and equip it with the natural projection $\pi_h\:\calm_h\to Y$.
In the particular case when $Y=\{J\}\hook\calj^k_U$
we obtain the moduli space of $J$-holomorphic curves
$\calm_J:=\pi_{\!\calj\!}^{-1}(J)$.
In general, the projection~$\pi_h$ has fibers $\pi_h^{-1}(y)=\calm_{h(y)}$.
We shall denote the elements of $\calm_h$ by $(M,y)$,
where $M=u(S)$ and the map $u\:S\to X$ is $h(y)$-holomorphic.

\proclaim{Lemma 2.4}
Let $Y$ be a $C^k$-smooth finite dimensional manifold,
and let $h\:Y\to\calj^k_U$ be a $C^k$-smooth map.
Assume that for some $y_0\in Y$ and $M_0=u_0(S)\in\calm_{h(y_0)}$
the map $\overline\Psi\circ dh\:T_{y_0}Y\to\roman H^1_D(S,N_{u_0})$
is surjective.
Then, in some neighborhood of $(M_0,y_0)\in\calm_h$,
the space $\calm_h$ is a $C^k$-smooth manifold
with tangent space
$$
T_{(M,y)}\calm_h=\ker\bigl(D\oplus\Psi\circ dh\:
\cale_{u,h(y)}\oplus T_yY\longrightarrow\cale'_{u,h(y)}\bigr)
\big/du\bigl(\roman H^0(S,TS)\bigr).
\tag 2.1
$$
\endproclaim

\smallskip\noindent\bf Proof. \rm 
Let $y\in Y$ and let  $(u,J_S,h(y))\in\calp^*$ and $M=u(S)$.
Then $(M,y)\in\calm_h$.
It follows from the proof of Corollary~2.3
that the range of the homomorphism
$d\pi_{\!\calj\!}\:T_{(M,h(y))}\calm\to T_{h(y)}\calj^k_U$
coincides with $\ker\bigl(\overline\Psi_{(u,h(y))}\bigr)$
and  $\overline\Psi$ maps the cokernel
$\coker\bigl(d\pi_{\!\calj\!}\bigr)$
isomorphically onto $\roman H^1_D(S,N_u)$.
Hence  the  lemma follows from Lemma~2.1.

\definition{Definition 2.3}
Let $Y$ be a compact manifold, let
$h\:Y\to\calj^k_U$ be a $C^k$-smooth map, and let
$\calm_h\subset\calm\times Y$
the corresponding moduli space
and  $(M_0,y_0)\in\calm_h$ be a point in this space. Then the
{\it  component through $(M_0,y_0)$ of the space
$\calm_h$}  is the subset $\calm_h(M_0,y_0)$ of $(M,y)\in\calm_h$
such that for each neighborhood $W$ of the image
$h(Y)\subset\calj^k_U$
there exists a continuous path
$\gamma\:[0,1]\to\calm$
with the following properties:
\roster
\item"(a)" $\gamma(0)=(M_0,h(y_0))$ and $\gamma(1)=(M,h(y))$,
that is, $\gamma$ joins $(M_0,y_0)$ with $(M,y)$ in~$\calm$;

\item"(b)" $J_t:=\pi_{\!\calj\!}(\gamma(t))\in W\subset\calj^k_U$
for all $t\in[0,1]$, that is,
the corresponding path of almost complex structures~$J_t$
lies in the  fixed neighborhood $W$ of~$h(Y)$.
\endroster

\proclaim{Lemma 2.5}
Assume that $h\:Y\to\calj^k_U$, let
$(M_0,y_0)\in\calm_h$, and  $\calm_h(M_0,y_0)$ are as in
Definition~$2.3$. Then we have the following results\rom:
\roster
\item"(i)" $\calm_h(M_0,y_0)$ is a closed subset of~$\calm_h;$

\item"(ii)"
if the component $\calm_h(M_0,y_0)$ is compact,
then there exists a subset  $\calm_h^0$
that contains $\calm_h(M_0,y_0)$
and is compact and open in~$\calm_h;$

\item"(iii)"
if the component $\calm_h(M_0,y_0)$ is non-compact,
then there exists a continuous path $\gamma\:[0,1)\to\calm_h$
with the following properties\rom:\rm
\itemitem{(a)} \it $\gamma(0)=(M_0,h(y_0))$,
that is,~$\gamma$ starts at~$(M_0,y_0);$\rm 
\itemitem{(b)} \it
there exists a sequence $t_n\nearrow1$
such that the sequence  $(M_n,J_n):=\gamma(t_n)$
lies in $\calm_h$ and is discrete there,
while the sequence $\{J_n\}$ converges in~$C^k$
to some~$J^*\in\calj^k_U.$
\endroster
\endproclaim

\smallskip\noindent\bf Proof. \rm 
(i) Let $(M',y')\in\overline{\calm_h(M_0,y_0)}\subset\calm_h$ and
let $J'=h(y')$.
Let $W$ be a neighborhood of $h(Y)\subset\calj_U^k$
and let $\{(M_n,y_n)\}$ be a sequence in~$\calm_h$ that converges to~$(M',y')$.
Then there exists a ball $B\ni(M',J')$ in~$\calm$
such that its projection to $\calj^k_U$ lies in        ~$W$.
Since $(M_n,h(y_n))$ belongs to~$B$ for $n$ sufficiently large,
there exists a path $(M_t,J_t)$ in~$\calm$
joining $(M_0,h(y_0))$ and~$(M',J')$
such that $J_t\in W$ for all~$t\in[0,1]$.
Hence $\calm_h(M_0,y_0)$ is closed.

(ii) Let $(M',y')\in\calm_h$ and let $J'=h(y')$.
We choose a finite-dimensional subspace
$F\subset T_{J'}\calj^k_U$ such that the map
$D_{u',J'}\oplus\Psi\:\cale_{u',J'}\oplus F\to\cale'_{u',J'}$
is surjective.
Let $B\ni0$ be a ball in~$F$.
Then there exists a $C^k$-smooth map $H\:Y\times B$
such that $H(y,0)\equiv h(y)$ and
$dH_{(y',0)}\:T_{(y',0)}(Y\times B)\to T_{J'}\calj^k_U$
induces an isomorphism
$T_0B\overset\cong\to\longrightarrow F\subset T_{J'}\calj^k_U$.
Thus, $\calm_H$ contains a neighborhood $V\ni(M',y',0)$
such that it is a manifold and $V\cap\calm_h$ is closed in~$V$.
It follows that $\calm_h$
is a locally compact topological space.

Since $\calm_h(M_0,y_0)$ is a compact subset of~$\calm_h$,
it has a neighborhood $V$ with compact closure
$\overline V\subset\calm_h$.
Assume that $W_i\subset\calj^k_U$ make up a fundamental
system of neighborhoods of $h(Y)$
and, in particular,  $\bigcap_iW_i=h(Y)$.
Let  $V_i$ be the set of all $(M,y)\in\calm_h$
such that $(M,h(y))$ can be joined with  $(M_0,h(y_0))$
by a path $(M_t,J_t)$ in $\calm$ with~$J_t$ in~$W_i$.
Then we have $\bigcap_i V_i=\calm_h(M_0,y_0)$.
The same argument as in part~(i)
shows that the $V_i$ are open and closed in~$\calm_h$.

We claim that there exists a positive integer $N\in\Bbb N$
such that $\bigl(\bigcap_{i=1}^NV_i\bigr)\cap\overline V\subset V$.
For otherwise,
 for each $n\in\Bbb N$ there  exists
$(M_n,y_n)\in\bigl(\bigcap_{i=1}^nV_i\bigr)\cap\overline V\setminus V$.
Then, however, a subsequence of $\{(M_n,y_n)\}$
 converges to some
$(M^*,y^*)\in\bigl(\bigcap_{i=1}^\infty V_i\bigr)\cap\overline V\setminus V$,
which is impossible because ~$\bigcap_{i=1}^\infty V_i\subset V$.

For such $N\in\Bbb N$ the set     
$$
\calm^0_h:=\biggl(\bigcap_{i=1}^NV_i\biggr)\cap\overline V
=\biggl(\bigcap_{i=1}^NV_i\biggr)\cap V
\tag 2.2
$$
satisfies the assumptions of part~(ii) of the lemma.

(iii)
Assume that $\calm_h(M_0,y_0)$ is non-compact.
Then there exists a discrete sequence
$\{(M_n,y_n)\}$ in $\calm_h(M_0,y_0)$.
Since $Y$ is compact, we may assume
that the $y_n$ converge to some~$y^*$.
For each $n\in\Bbb N$ we fix a path
$\gamma_n\:[0,1]\to\calm_h$
between $(M_{n-1},h(y_{n-1}))$ and~$(M_n,h(y_n))$.
We set $t_n:=1-\nomathbreak 2^{-n}$.
For $t\in[t_{n-1},t_n]$ we set
$\gamma(t):=\gamma_n(2^n(t-t_{n-1}))$.
Then $\gamma\:[0,1)\to\calm_h$ and $t_n\nearrow1$
are the required path and sequence.

\proclaim{Theorem 2.6}
Let $(M_0,J_0)\in\calm$ and let
$h\:[0,1]\to\calj^k_U$ be a $C^k$-smooth path
 with~$h(0)=\nomathbreak J_0$.
Assume that there exists
a compact open subset $\calm_h^0$ of $\calm_h$
 containing $(M_0,J_0)$.
Assume further that  the index
$\ind(\pi_{\!\calj\!})\allowmathbreak =2(c_1(X)[M_0]+(n-3)(1-g))$
is non-negative.
Then $h$ can be $C^k$-approximated
by smooth maps $h_n\:[0,1]\to\calj^k_U$
with the following properties\rom:
\roster
\item"(i)" each $\calm_{h_n}$ contains a component $\calm_{h_n}^0$
that is a $C^k$-smooth manifold of the expected dimension
$\dim_{\Bbb R}(\calm_{h_n}^0)=\ind(\pi_{\!\calj\!})+\nomathbreak 1;$

\item"(ii)" the sets
$\calm_{h_n(0)}^0:=\pi_{h_n}^{-1}(0)\cap\calm_{h_n}^0$ and
$\calm_{h_n(1)}^0:=\pi_{h_n}^{-1}(1)\cap\calm_{h_n}^0$
are also $C^k$-smooth manifolds of  the expected dimension
$\ind(\pi_{\!\calj\!})$ \rom($\calm^0_{h_n(1)}$ can be empty\rom{!),}
and $\calm_{h_n}$ is a $C^k$-smooth bordism 
between 
$\calm_{h_n(0)}^0$ and~$\calm_{h_n(1)}^0;$

\item"(iii)" each $\calm_{h_n(0)}^0$ can be joined with ~$(M_0,J_0)$
by a path in~$\calm$,
that~is, there exist $C^k$-smooth paths
$\gamma_n\:[0,1]\to\calm$ with $\gamma_n(0)=(M_0,J_0)$ and
$\gamma_n(1)\in\calm_{h_n(0)}^0;$
in particular, all $\calm_{h_n(0)}^0$ are not empty\rom;

\item"(iv)" for each element $(M,J)$ of $\calm_{h_n}$
we have the inequality
$\dim_{\Bbb R}\roman H^1_D(S,N_M)\le1.$
\endroster
\endproclaim

\smallskip\noindent\bf Proof. \rm 
We denote $\calm_h^0$ by~$K$.
Let $\cale_K$ and $\cale'_K$ be     the Banach bundles over~$K$
induced by the bundles $\cale\to\calm$ and $\cale'\to\calm$,
 respectively.
 Let also $\calt:=h^*T\calj^k_U$ be     the pull-back
of $T\calj^k_U$ to~$[0,1]$.

By Lemma 2.2, for each $(M,J)\in K$ with $M=u(S)$
there exists  $m_{(M,J)}\in\Bbb N$ and a $C^k$-smooth homomorphism
$P_{(M,J)}\:F_{(M,J)}\to\calt$ of the trivial vector
bundle $F_{(M,J)}$ of rank $\rank F_{(M,J)}=m_{(M,J)}$ over ~$[0,1]$
such that the operator
$$
D_{(u,J)}\oplus\Psi_{(u,J)}\circ P_{(M,J)}\:
\cale_{(M,J)}\oplus F_{(M,J)}\to\cale'_{(M,J)}
\tag 2.3
$$
is surjective.
Note  that the operator (2.3)
remains surjective for all $(M',J')$ from a neighborhood of~$(M,J)$.
Since $K$ is compact,
we may choose finitely many pairs  $(F_{\!(M,J)},P_{\!(M,J)})$ of this kind
and then take their sum to obtain
a homomorphism $P\:F\to\calt$ of   the trivial vector
bundle~$F\cong[0,1]\times\Bbb R^m$,~$m\gg0$,
such that the homomorphism
$D_{(u,J)}\oplus\Psi_{(u,J)}\circ P\:
\cale_{(M,J)}\oplus F\to\cale'_{(M,J)}$
is surjective for all~$(M,J)\in K$.

For all $\dot J\/$ in a small ball in~$T_J\calj^k_U$
we set $\exp_J(\dot J):=J\dfrac{(1-J\dot J/2)}{(1+J\dot J/2)}$.
Differentiating the identity $J^2=-1$
we obtain the equality $J\dot J=-\dot JJ$
holding for  all $\dot J\in T_J\calj^k_U$.
This gives us the following relation:
$$
\biggl(J\frac{(1-J\dot J/2)}{(1+J\dot J/2)}\biggr)^2=-1.
$$
Thus,   $\exp_J$ ranges       in~$\calj^k_U$.
Direct calculations show that the derivative of $\exp_J$
at~$0\in T_J\calj^k_U$ is the identity map of~$T_J\calj^k_U$.
Thus, $\exp_J$ is the natural exponential map for~$\calj^k_U$.

We consider  a sufficiently small ball $B=B(0,r)\subset\Bbb R^n$
and define the map $H^*\:[0,1]\times B\to\calj^k_U$
by the formula $H^*(t,y):=\exp_{h(t)}(P(t,y))$.
By our construction of $H^*$ the map
$$
D_{(u,J)}\oplus\Psi_{(u,J)}\circ dH^*(t,0)\:
\cale_{(M,J)}\oplus T_{(t,0)}([0,1]\times B)\to\cale'_{(M,J)}
$$
is surjective for all~$(M,J)\in K$ with $J=h(t)$.
After a small perturbation of~$H^*$
we obtain a $C^k$-smooth function $H\:[0,1]\times B\to\calj^k_U$
with the following properties:
\roster
\item"(i)" $H(t,0)=h(t)$, that~is, $H$ is a deformation of~$h$
with parameter space $B=B(0,r)\subset\Bbb R^n$;

\item"(ii)" $D_{(u,J)}\oplus\Psi_{(u,J)}\circ dH(t,0)\:
\cale_{(M,J)}\oplus T_0B\longrightarrow\cale'_{(M,J)}$
is surjective for all $(M,t)\in K$ with $J=h(t)$ and~$M=u(S)$;

\item"(iii)"  the map $H$ is $\calj^\infty_U$-valued and
$C^\infty$-smooth in a neighborhood of each point $(t,y)\in[0,1]\times B$
such that~$y\ne0$.
\endroster

We identify $K$ and $\calm_h$ with subsets of~$\calm_H$
using the natural embedding
$(M,t)\in\calm_h\mapsto(M,t,0)\in\calm_H$.
By Lemma~2.2 there exists a neighborhood $V$ of~$K$ in~$\calm_H$
that is a $C^k$-smooth manifold.
Shrinking this neighborhood of~$K$ if necessary
we may assume that its closure~$\overline V$ is compact
and does not intersect the other components of $\calm_h$, 
distinct  from~$\calm_h^0=K$,
that~is, $\overline V\cap\calm_h=V\cap\calm_h=K$.
Let $p\:V\to B$ be the natural projection
defined by the formula~$(u,t,y)\mapsto y$.

We claim that there exists a smaller ball
$B_1=B(0,r_1)\subset B$ such that
the set $p^{-1}(y)\subset V$ is compact for all~$y\in B_1$.
Assume the contrary.
Then there exists a sequence of~$y_n\in B$
converging to~$0\in B$ such that the  inverse images $p^{-1}(y_n)$
are non-compact.
Since $\overline V$ is compact,
 there exist $M_n$ and $t_n$ such that~$(M_n,t_n,y_n)$
lies in the closure~$\overline{p^{-1}(y_n)}\subset\overline V$,
but not in~$V$.
Passing to an appropriate subsequence of~$(M_n,t_n,y_n)$
we may assume that it converges to some point~$(M^*,t^*,0)$.
Then,  however, $(M^*,t^*,0)\in\calm_h\cap\overline V=K$
and therefore $(M^*,t^*,0)\in V$.
On the other hand $\overline V\setminus V$ is compact
and therefore  $(M^*,t^*,0)$ lies in~$\overline V\setminus V$. 
This contradiction shows that a  required ball $B_1$ exists
indeed.

We set $V_1:=p^{-1}(B_1)\subset V$.
By the construction of~$H$ the restriction of the projection
$p\:V_1\setminus K\to\check B_1:=B_1\setminus\{0\}$
is $C^\infty$-smooth.
By Sard's theorem (see, for example,~\cite{12};~\S\,3.4)
there exists a dense subset $B_1^*$ of $  B_1$
such that $p^{-1}(y)$ is a $C^\infty$-smooth compact manifold
for all~$y\in B_1^*$.
We fix a sequence $y_n\in B_1^*$ converging to~$0\in B$
and set $h_n(t):=H(t,y_n)$.
Let $\calm_{h_n}^0:=\calm_{h_n}\cap V$,
so that $\calm_{h_n}^0=p^{-1}(y_n)$.
Then each $\calm_{h_n}^0$ is a $C^\infty$-smooth non-empty manifold
that can be joined with ~$(M_0,J_0)$ by a path in~$\calm$.

By Lemma 2.4 the tangent space to
$V_1\subset\calm_H$ in~$(u,t,y)$
is canonically isomorphic to the quotient
$$
\ker\bigl(D_{u,H(t,y)}\oplus\Psi\circ dH\:
\cale_{u,H(t,y)}\oplus T_{(t,y)}([0,1]\times B_1)
\to\cale'_{u,H(t,y)}\bigr)
\big/du\bigl(\roman H^0(S,TS)\bigr).
$$
Since $p\:V_1\to B_1$ is the projection
$(u,t,y)\in V_1\mapsto y\in B_1$,
 the differential $dp_{(u,t,y)}$
maps a tangent vector of the form
$(\dot u,\dot t,\dot y)\in T_{(u,t,y)}V_1$
to~$\dot y\in T_yB_1$.
This means that $dp_{(u,t,y)}$ is the restriction to
$\ker(D_{u,H(t,y)}\oplus\Psi\circ dH)$
of the linear projection
$p_B\:\cale_{u,H(t,y)}\oplus T_{(t,y)}([0,1]\times B_1)\to T_yB_1$
defined by the formula~$p_B(\dot u,\dot t,\dot y)=\dot y$.
In particular, if $y=y_n$, then the map  $dp_{(u,t,y)}$ is surjective,
which means  the surjectivity of the map
$$
p_B\:\ker(D_{u,H(t,y)}\oplus\Psi\circ dH)\to T_{y_n}B_1.
$$
The last result    is equivalent to the surjectivity of the map
$$
D_{u,H(t,y_n)}\oplus\Psi\circ dH\oplus p_B\:
\cale_{u,H(t,y_n)}\oplus T_{(t,y_n)}([0,1]\times B_1)
\to\cale'_{u,H(t,y_n)}\oplus T_{y_n}B_1
$$
and, therefore,                     of the map
$$
D_{u,h_n(t)}\oplus\Psi\circ dh_n\:
\cale_{u,h_n(t)}\oplus T_t[0,1]\to\cale'_{u,h_n(t)}.
\tag 2.4
$$
Consequently, $\dim_{\Bbb R}\roman H^1_D(S,N_M)\le1$ for all
$(M,t)\in\calm_{h_n}\cap V_1$.

\proclaim{Corollary 2.7}
Under the hypothesis  of Theorem~$2.6$
assume in addition that $S$ is the sphere~$S^2$.
Then for all $(M,t)\in\calm^0_{h_n}$ the associated
$D_N$-operator is surjective, that~is, $\roman H^1_{D_N}(S^2,N_M)=0$.

Moreover, $\calm_{h_n}$ is the trivial bordism\rom:
$\calm_{h_n}\cong\calm_{h_n(0)}\times[0,1]$.
In particular, for each $h_n(0)$-holomorphic sphere
$M_0\in\calm_{h_n(0)}$
there exists a continuous family of $h_n(t)$-holomorphic spheres
$M_{n,t}=u_{n,t}(S^2)$ with~$M_{n,0}=M_0$.
\endproclaim

\smallskip\noindent\bf Proof. \rm 
Assume that $\roman H^1_{D_N}(	M,N_M)\ne0$ for some
$(M,t)\in\calm^0_{h_n}$.
Then $\roman H^1_{D_N}(M,N_M)$ $=1$ by Theorem~2.6.
However, this contradicts the result  of Corollary~1.6
for  $S=S^2$ and $L:=N_M$.

Let $(M,t)\in\calm^0_{h_n}$ satisfy $M=u(S)$ and $J=h_n(t)$.
Let also~$\dot J\ne0\in dh_n(T_t[0,1])$.
Then, by Lemma~2.1 and Corollary~2.3
the tangent space $T_{(M,t)}\calm^0_{h_n}$
is canonically isomorphic to
$$
\ker\bigl(D_{u,J}\oplus\Psi\:
\cale_{u,J}\oplus\Bbb R\dot J\to\cale'_{u,J}\bigr)
\big/du\bigl(\roman H^0(S,TS)\bigr),
$$
and the differential of the projection
$d\pi_{h_n}\:T_{(M,t)}\calm^0_{h_n}\to T_t[0,1]\cong\Bbb R$
takes the form
$d\pi_{h_n}[v,a\dot J]\allowmathbreak =\nomathbreak a$. 
When $S=S^2$, the space $\roman H^1(S,TS)$ is trivial
and Corollary~1.8 shows that
$D_{u,J}\:\cale_{u,J}\to\cale'_{u,J}$ is surjective.
Hence, for $a\ne0\in\Bbb R$ there exists  $v\in\cale_{u,J}$
such that~$[v,a\dot J]\in T_{(M,t)}\calm^0_{h_n}$.
This means that for each  $(M,t)\in\calm^0_{h_n}$
the projection $d\pi_{h_n}\:T_{(M,t)}\calm^0_{h_n}\to T_t[0,1]$
is surjective.
Since $\calm^0_{h_n}$ is compact,
there exists a diffeomorphism
$\calm_{h_n}\cong\calm_{h_n(0)}\times[0,1]$.

\head
\S\,3. Gromov topology and deformations of non-compact holomorphic curves
\endhead

The techniques developed in the previous section
enables us to construct local deformations 
of $J_t$-holomorphic spheres $M_t=u_t(S^2)$
for appropriate families of almost complex structures~$J_t$.
The obstruction to
the existence of a~`global family'
is that the sphere~$M_t$
can eventually `break down' into several components. 
For our purposes here  we must    know
the exact fashion of this break,
and we  must also learn to deform the reducible
(that~is, consisting of several components)
curves produced by  such a `breakdown'.
We start with the indication of a suitable 
category of reducible curves.

\definition{Definition 3.1}
The complex analytic set
$A_0:=\{(z_1,z_2)\in\Delta^2:z_1\cdot z_2=0\}$
is called the {\it standard node\/}.
A {\it nodal curve\/} is a connected reduced
complex space~$C$ of dimension~$1$
having finitely many    irreducible components
and with singularities only at  finitely many
 {\it nodal points\/}  
that have    neighborhoods isomorphic
to the standard node.
Furthermore, we assume that the boundary~$\partial C$
of ~$C$ consists of  finitely many    smooth circles
and  $\overline C:=C\cup\partial C$ is compact.
The case of $\partial C=\varnothing$ is not excluded.

\definition{Definition 3.2}
A smooth oriented real surface~$\Sigma$ with boundary~$\partial\Sigma$
{\it parametrizes\/} a nodal complex curve~$C$
if there exists a continuous map $\sigma\:\Sigma\to C$
with the following properties:
\roster
\item"(1)" if $a\in C$ is a nodal point, then
$\gamma_a:=\sigma^{-1}(a)$ is a smooth embedded circle in~$\Sigma$;

\item"(2)" $\sigma\:\overline\Sigma\setminus
\bigcup_{i=1}^N\gamma_{a_i}\to\overline C\setminus\{a_1,\dots,a_N\}$
is a diffeomorphism, where $\{a_1,\dots,a_N\}$ is
the set of all nodal points in~$C$.
\endroster
The map $\sigma$ is called a {\it parametrization\/} of~$C$.

\medskip
\vbox{\xsize=0.53\hsize\nolineskip\rm
\putm[.13][.03]{\gamma_1}%
\putm[.24][.21]{\gamma_2}%
\putm[.487][0.215]{\gamma_3}%
\putm[.62][.23]{\gamma_4}%
\putm[.84][.17]{\gamma_5}%
\putt[1.1][0]{\advance\hsize-1.1\xsize%
\centerline{Fig. 1}\smallskip 
The projection $\sigma $ `contracts'
circles $\gamma_1, ...,\gamma_5$ into the nodal points
$a_1,..., a_5$.
}
\epsfxsize=\xsize\epsfbox{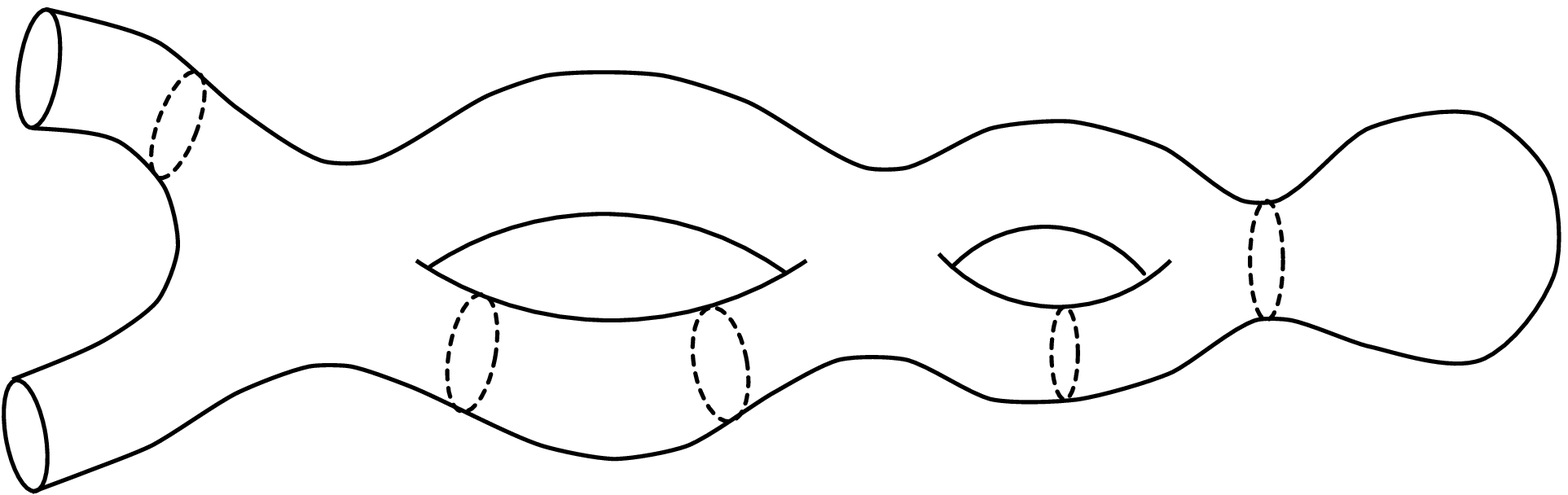} \putm[.56][.0]{\bigg\downarrow\sigma}
\vskip12pt
\putm[.14][.065]{a_1}%
\putm[.31][.27]{a_2}%
\putm[.445][.029]{a_3}%
\putm[.67][.245]{a_4}%
\putm[.805][.21]{a_5}%
\epsfxsize=\xsize\epsfbox{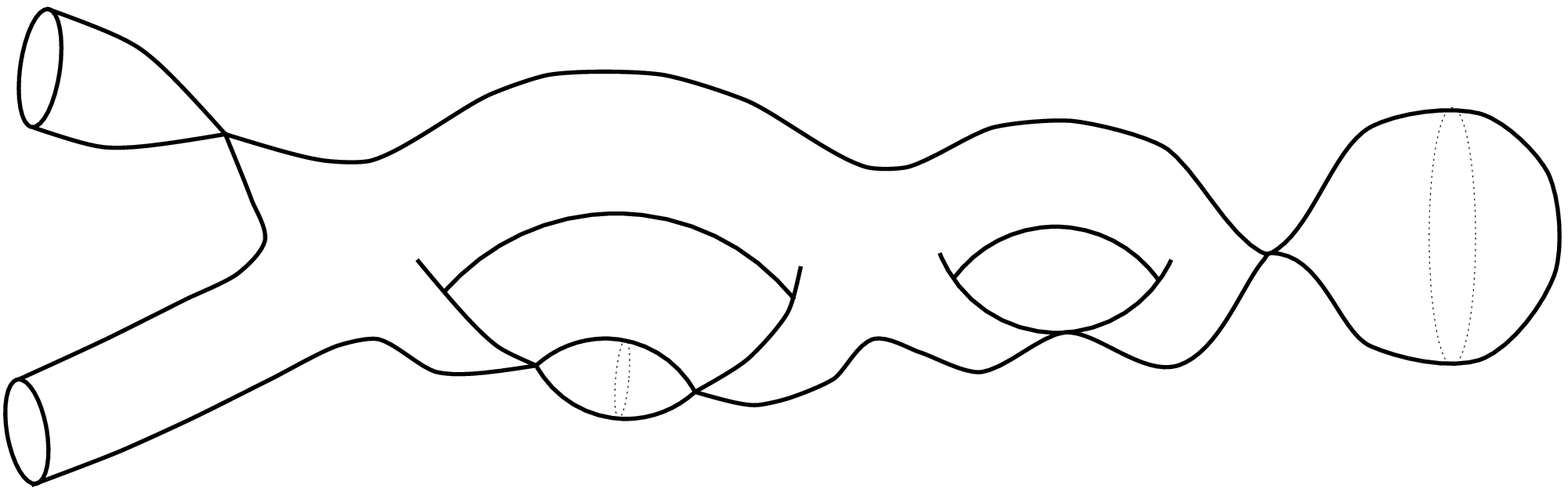}
}

\smallskip

Note that the parametrization $\sigma$ is not unique:
if $g\:\Sigma\to\Sigma$ is a diffeomorphism,
then $\sigma\circ g\:\Sigma\to C$ is also a parametrization.

\definition{Definition 3.3}
A {\it stable curve over an almost complex manifold~$(X,J)$}
is a pair $(C,u)$, where $C$ is a nodal curve with boundary~$\partial C$
and $u\:C\to X$ is a $J$-holomorphic map with the following property:
if the map $u$ is constant on a compact irreducible component~$C_j$ of
the curve~$C$,
then the group of biholomorphic automorphisms of $C_j$
preserving the nodal points of $C_j$ is finite.

It is easy to see that this condition imposes only
the following two restrictions:
\roster
\item"(1)"
if $C_j$ is rational, that~is, biholomorphic to~${\Bbb C}{\Bbb P}^1$,
then $C_j$ contains at least three nodal points of~$C$;

\item"(2)" if $C_j$ is a torus, then~$C_j$ contains at least one
nodal point of~$C$.
\endroster

This definition was given by Deligne and~Mumford
in the case of abstract algebraic curves
and was later  generalized by Kontsevich
to the case of maps into	 ~$X$.

Stable curves occur  in a natural way when
one attempts to compactify the space of
embedded or immersed curves in~$X$.
Namely, one introduces the Gromov topology
in the set of all stable curves over~$X$.
We shall describe this topology
by describing     convergent sequences.

Let $\{J_n\}$ be a sequence of $C^k$-smooth
almost complex structures on~$X$
 convergent  to some~$J$ in~$C^k$.
Let $\{(C_n,u_n)\}$ be a sequence of
$J_n$-holomorphic curves that are  stable over~$X$
and are parametrized by the same real surface~$\Sigma$.

\definition{Definition 3.4}
The sequence $\{(C_n,u_n)\}$ converges to
a stable $J$-holomorphic curve $(C_\infty,u_\infty)$ over~$X$
if the following holds:
\roster
\item"(1)" $C_\infty$ is also parametrized by~$\Sigma$;
moreover, there exist parametrizations $\sigma_n\: \Sigma 
\allowbreak \to C_n$
and $\sigma_\infty:\Sigma\to C_\infty$
such that $u_n\circ\sigma_n$ converge to
$u_\infty\circ\sigma_\infty$
in the $C^0$-topology on~$\overline\Sigma$
(that is, up to the boundary);

\item"(2)" let $\{a_1,\dots,a_N\}$ be the nodal points of~$C_\infty$
and let $\gamma_i:=\sigma_\infty^{-1}(a_i)$;
let also $K$
be a compact subset of $\overline\Sigma\setminus\bigcup_i\gamma_i$;
then for all $n\ge n^*(K)$ the set $\sigma_n(K)$
does not contain  nodal points of~$C_n$,
$u_n\circ\sigma_n$ converges to~$u_\infty\circ\sigma_\infty$
in the $L^{1,p}$-topology on~$K$,
and the inverse images $\sigma_n^*J_{C_n}$
of the complex structures $J_{C_n}$ on the curves~$C_n$
converge to the inverse image $\sigma_\infty^*J_{C_\infty}$
of the complex structure $J_{C_\infty}$ on~$C_\infty$
in the $C^\infty$-topology on~$K$.
\endroster

\definition{Definition 3.5}
An {\it annulus $A$} on a real surface or  a complex curve
is a domain  diffeomorphic (respectively,  biholomorphic)
to the standard annulus $A_{r,R}:=\{z\in\Bbb C:r<|z|<R\}$.
A subdomain
of a real surface or a complex curve
that is diffeomorphic to a   disc with two holes
is called {\it pants\/}.
In both cases we assume that the boundary of the domain
consists of smooth embedded circles.
An annulus $A$ {\it is adjacent\/} to a  circle~$\gamma$
if $\gamma$ is a  component of its boundary~$\partial A$.

\medskip
\vbox{
 \hbox{%
   \vtop{
\hsize=.46\hsize\epsfxsize=\hsize\epsfbox{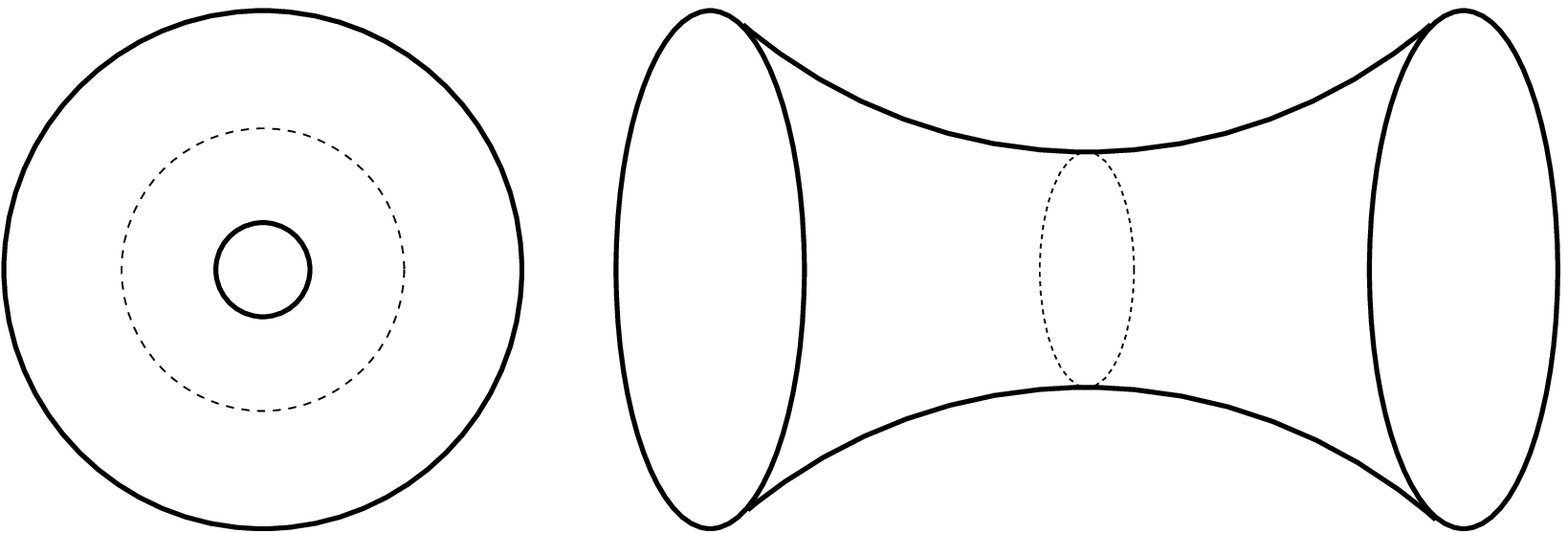}%
\smallskip
\centerline{Fig.\ 2. Annulus}
\smallskip
It is useful to represent an annulus with a complex structure 
as a cylinder or a ``tube". Contracting the circle
in the middle of the annulus we obtain the nodal point.
   }

\hskip.05\hsize
   \vtop{
\hsize=.49\hsize\epsfxsize=\hsize\epsfbox{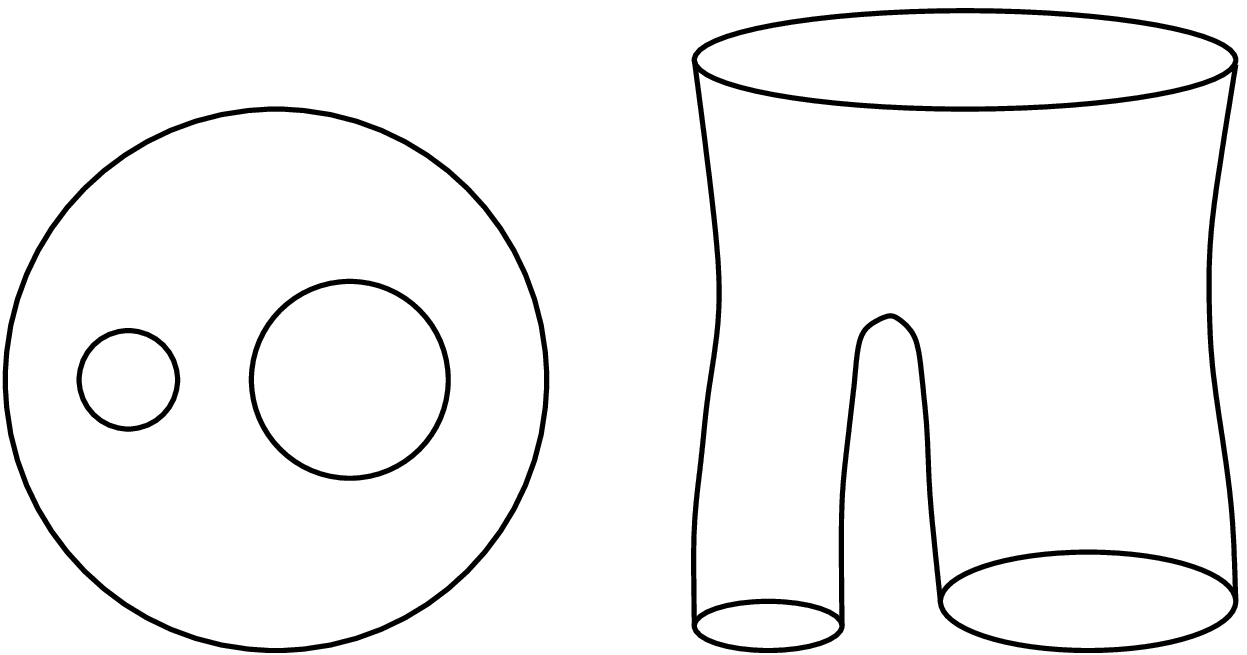}
\smallskip
\centerline{Fig.\ 3. Pants.}
\smallskip
One can think about pants as a
sphere with three holes.
   }
 }
}

\medskip
The significance of the Gromov topology is expressed
by the following result.

\proclaim{Theorem 3.1}
Let $\{(C_n,u_n)\}$ be a sequence
of stable $J_n$-holomorphic curves over $X$
that satisfy the following conditions\rom:
\roster
\item"(1)" the $J_n$ are $C^k$-smooth
and converge to $J$ in~$C^k$ for some~$k\ge1;$

\item"(2)" $\operatorname{area}[u_n(C_n)]\le A$ for all~$n;$

\item"(3)" all $C_n$ are parametrized by the same
real surface~$\Sigma;$

\item"(4)" the sequence $(C_n,u_n)$ converges near the boundary\rom;
this means that there exist parametrizations~$\sigma_n\:\Sigma\to C_n$
and annuli~$A_\alpha$ in $\Sigma$ such that \rm
\itemitem{(a)} \it each~$A_\alpha$ is adjacent
to a single component $\gamma_\alpha$ of the boundary~$\partial\Sigma$,
the images $\sigma_n(A_\alpha)$ do not contain the nodal points of~$C_n$
and the complex structures~$\sigma_n^*J_{C_n}$
are constant on each annulus~$A_\alpha;$\rm
\itemitem{(b)} \it the maps $u_n\circ\sigma_n$
converge in the~$L^{1,p}$-topology on each~$A_\alpha$.
\endroster

Then there exists a subsequence of stable curves,
also denoted by $\{(C_n,u_n)\}$,
that converges in the Gromov topology
to a curve $(C_\infty,u_\infty)$ stable over~$X$.
\endproclaim

This theorem is due to Gromov~\cite{6};
proofs can be also found in~\cite{13} and~\cite{14}.
We shall require a more precise description of Gromov convergency,
which can be derived from~\cite{14}.

\proclaim{Proposition 3.2}
Under the assumptions of Theorem~$3.1$
one can  choose para\-metri\-zations $\sigma_n\:\Sigma\to C_n$
and a finite cover of $\Sigma$ by open subsets~$\{V_\alpha\}$
and $\{V_{\alpha\beta}\}$
with the following properties \rom{(see~Fig.~6):}
\roster
\item"(i)" each $V_\alpha$ is a disc,  an annulus, or pants,
and $V_\alpha\cap V_\beta=\varnothing$ for all~$\alpha\ne\beta$;

\item"(ii)" for each component of the boundary of $\Sigma$
there exists a single annulus $V_\alpha$ adjacent to this component\rom;

\item"(iii)" each $V_{\alpha\beta}$ is an annulus
 intersecting $V_\alpha$ and $V_\beta$ by annuli
$W_{\alpha,\beta}:=V_{\alpha\beta}\cap V_\alpha$ and
$W_{\beta,\alpha}:=V_{\alpha\beta}\cap V_\beta;$
furthermore,~$V_{\alpha\beta}=V_{\beta\alpha}$
and the annuli $W_{\alpha,\beta}$ and $W_{\beta,\alpha}$
are distinct and disjoint\rom;

\item"(iv)" $\sigma_n^*J_{C_n}\bigm|_{V_\alpha}$
does not depend on~$n$,
and the sequence of $u_n\circ\sigma_n|_{V_\alpha}$ converges in~$L^{1,p}$
to~$u_\infty\circ\sigma_\infty|_{V_\alpha}$\rom;

\item"(v)" there exist $L^{1,p}$-smooth maps
$\varphi_{\alpha\beta}^n\:V_{\alpha\beta}\to\Delta^2$ and
$\varphi_{\alpha\beta}^\infty\:V_{\alpha\beta}\to\Delta^2$
that induce biholomorphisms
$\overline\varphi_{\alpha\beta}^n\:
\sigma_n(V_{\alpha\beta})\overset\cong\to\longrightarrow
\{(z_1,z_2)\in\Delta^2:z_1\cdot z_2=\lambda_{\alpha\beta}^n\}$
and
$\overline\varphi_{\alpha\beta}^\infty\:
\sigma_\infty(V_{\alpha\beta})\overset\cong\to\longrightarrow
\{(z_1,z_2)\in\Delta^2:z_1\cdot z_2=\lambda_{\alpha\beta}^\infty\};$
here $\lambda_{\alpha\beta}^n\in\Delta$,
$\lambda_{\alpha\beta}^n$ converge to
$\lambda_{\alpha\beta}^\infty\in\Delta$, and
$\varphi_{\alpha\beta}^n$ converge to
$\varphi_{\alpha\beta}^\infty$ in~$L^{1,p};$

\item"(vi)"
for the coordinate functions $z_1$ and $z_2$ on $\Delta^2$
the compositions 
$z_1\circ\varphi_{\alpha\beta}^n\bigm|_{W_{\alpha,\beta}}$ and
$z_2\circ\varphi_{\alpha\beta}^n\bigm|_{W_{\beta,\alpha}}$
do not depend on~$n$\rom;
moreover, the images
$z_1\circ\varphi_{\alpha\beta}^n(W_{\alpha,\beta})$ and
$z_2\circ\varphi_{\alpha\beta}^n(W_{\beta,\alpha})$
are some annuli $A_{r_{\alpha,\beta}1}$ and $A_{r_{\beta,\alpha}1}$,
respectively\rom;

\item"(vii)" for a fixed limiting  curve $(C_\infty,u_\infty)$
the corresponding covering $\{V_\alpha,V_{\alpha\beta}\}$
of the real surface~$\Sigma$
can be chosen to depend only on~$\{(C_\infty,u_\infty)\}$,
that is, to be the same for all sequences~$\{(C_n,u_n)\}$
converging to~$(C_\infty,u_\infty)$.
\endroster
\endproclaim

Thus, Theorem 3.1 ensures that a deformation
of a pseudoholomorphic curve can
`break down' only into
pseudoholomorphic curves,
which gives us a possibility to continue the process of local deformation.

 We shall now consider local deformations of non-compact curves.

\definition{Definition 3.6}
A {\it Banach ball\/} is a ball in some complex Banach space.
A subset $\calm$ of the Banach ball~$B$
is called a {\it Banach analytic set of finite codimension\/}
(a BASFC) 
 if there exists a holomorphic
map~$F\:B\to\Bbb C^N$,~$N<\infty$
such that~$\calm=\{x\in B:F(x)=0\}$.

This concept is important
because, unlike  general Banach analytic sets,
 BASFC's  have   properties  similar to those
of usual finite-dimensional analytic sets.
Namely, we have the following result.

\proclaim{Theorem 3.3 \rm \cite{15}}
Let $B$ be a ball in a Banach space~$\calf$, let
$\calm\subset B$ be a BASFC, and 
$x_0$ a point in~$\calm$.
Then there exists a neighborhood $U\ni x_0$ in~$B$
such that $\calm\cap U$ is a finite union of BASFC's
$\calm_j$,
each of them 
irreducible at~$x_0$.

Moreover, each $\calm_j$ can be represented
as a proper ramified covering over a domain
in a closed linear subspace~${\calf}_j\subset{\calf}$
of finite codimension.
\endproclaim

Our aim in this section is to prove
the existence of a complete family
of holomorphic deformations of a stable curve over~$X$,
parametrized by a BASFC.
Before stating the result we introduce the following definition.

\definition{Definition 3.7}
Let $C$ be a nodal curve, let $E$ be a holomorphic vector bundle over~$C$,
and  let $C=\bigcup_{i=1}^lC_i$ be the decomposition of $C$
into irreducible components.
Assume that $E$ 
extends 
sufficiently smoothly 
up to 
the boundary~$\partial C$.
We define an {\it $L^{1,p}$-section $v$} of the bundle~$E$ over~$C$
as a collection $(v_i)_{i=1}^l$ of~$v_i\in L^{1,p}(C_i,E)$
such that at each nodal point $z\in C_i\cap C_j$ we have $v_i(z)=v_j(z)$.
We also define an {\it $E$-valued $L^p$-integrable $(0,1)$-form\/}
$\xi$ on~$C$ as a collection $(\xi_i)_{i=1}^l$ of $(0,1)$-forms
$\xi_i\in L^p(C_i,E\otimes\Lambda^{(0,1)})$.
Let $L^{1,p}(C,E)$ be     the Banach space
of $L^{1,p}$-sections of~$E$ over~$C$
and let $L^p(C,E\otimes\Lambda^{(0,1)})$ be the Banach space
of $L^p$-integrable $(0,1)$-forms~$C$.
We also denote  by ${\calh}^{1,p}(C,E)$
the Banach space of holomorphic $L^{1,p}$-sections of~$E$ over~$C$.

In a similar way, for each complex manifold $X$
we shall mean by $L^{1,p}(C,X)$
the set of all collections $u=(u_i)_{i=1}^l$ of 
maps 
$u_i\in L^{1,p}(C_i,X)$
satisfying the equality $u_i(z)=u_j(z)$ at each nodal point~$z\in C_i\cap C_j$.
It is easy to see that $L^{1,p}(C,X)$ is a Banach manifold
with tangent space~$T_uL^{1,p}(C,X)=L^{1,p}(C,u^*TX)$.
Further, we denote
the space of all holomorphic $L^{1,p}$-maps from~$C$ to~$X$
 by ${\calh}^{1,p}(C,X)$.
Note  that for each ~$u\in{\calh}^{1,p}(C,X)$ we have
$u(C)\subset u(\overline C)\Subset X$.

The aim  of this section is the following result.

\proclaim{Theorem 3.4}
Let $(X,J)$ be a complex manifold and 
$(C_0,u_0)$ a stable complex curve over~$X$
parametrized by a real surface~$\Sigma$.
Then there exist BASFC's $\calm$ and ${\calc}$
and holomorphic maps $F\:{\calc}\to X$ and
$\pi\:{\calc}\to\calm$
with the following properties\rom:
\roster
\item"(a)" for each $\lambda\in\calm$ the fiber
$C_\lambda=\pi^{-1}(\lambda)$ is a nodal curve parametrized by~$\Sigma$,
and $C_{\lambda_0}\cong C_0$ for some~$\lambda_0;$

\item"(b)" for $F_\lambda:=F\big|_{C_\lambda}$
the pair $(F_\lambda,C_\lambda)$ is a stable curve over~$X$,
and $F_{\lambda_0}=u_0;$

\item"(c)" if $(C',u')$ is a stable curve over~$X$
that is sufficiently close to~$(C_0,u_0)$ in the Gromov topology,
then there exists $\lambda'\in\calm$
such that $(C',u')=(C_{\lambda'},F_{\lambda'});$

\item"(d)" for an appropriate integer $N\in\Bbb N$ and a small ball~$B$
in the Banach space ${\calh}^{1,p}(C_0,u_0^*TX)\allowbreak\oplus\Bbb C^N$
the BASFC $\calm$ can be realized
as the zero set of a holomorphic map $\Phi$ from~$B$
into the finite dimensional space~$\roman H^1(C,u_0^*TX)$.
\endroster
\endproclaim

The proof relies on the construction of local deformations
of stable curves and on the analysis of  conditions
for  patching together
local models.
The following result on the solution
of a Cousin-type problem plays an important role in this proof.

\proclaim{Lemma 3.5}
Let $C$ be a nodal curve and 
$E$ a holomorphic vector bundle over~$C$ that is
$C^1$-smooth up to the boundary.
Let $\{V_i\}_{i=1}^l$ be a finite cover of~$C$
by Stein domains with piecewise smooth boundaries.
Set $V_{ij}:=V_i\cap V_j$ and assume that
all triple intersections
$V_i\cap V_j\cap V_k$ with $i\ne j\ne k\ne i$ are empty.

Then for all $2\le p<\infty$ the {\v C}ech coboundary operator
$$
\matrix
\delta\:&\sum_{i=1}^l{\calh}^{1,p}(V_i,E)&\longrightarrow&
\sum_{i<j}{\calh}^{1,p}(V_{ij},E),
\\
\delta\:&(v_i)_{i=1}^l&\longmapsto&(v_i-v_j)
\endmatrix
\tag 3.1
$$
has the following properties\rom:
\roster
\item"(i)" the image $\Im(\delta)$ is closed and has finite codimension\rom;
moreover,
$\coker(\delta)=\roman H^1(C,E)=\roman H^1(C_{\text{\rom{comp}}},E)$,
where $C_{\text{\rom{comp}}}$
is      the union of all compact irreducible components of~$C;$

\item"(ii)" the kernel $\ker(\delta)$ is isomorphic to
${\calh}^{1,p}(C,E)$
and has  a closed complementary subspace.
\endroster
\endproclaim

\smallskip\noindent\bf Proof. \rm 
Before  the {\v C}ech complex,
let us discuss  the corresponding $\overline\partial$-problem.
We consider the following operator:
$$
\overline\partial\:L^{1,p}(C,E)\longrightarrow
L^p(C,E\otimes\Lambda^{(0,1)}_C).
\tag 3.2
$$
First, we shall prove  that the properties of
this operator are similar to those  of~$(3.1)$,
that~is, $\ker(\overline\partial)$ has a closed complementary subspace,
$\Im(\overline\partial)$ is closed and of finite codimension,
and $\coker(\overline\partial)=\roman H^1(C,E)=\roman H^1(C_{\text{\rom{comp}}},E)$.
Moreover, we shall construct explicit isomorphisms
between the (co)kernels of~(3.1) and~(3.2).

Since the boundary of $C$ is smooth,
there exist nodal curves $C^+$ and~$C^{++}$
such that $C\Subset C^+\Subset C^{++}$ and
the difference~$C^+\setminus\overline C$
(respectively, $C^{++}\setminus\overline C^+$)
consists of annuli $A^+_\alpha$ (respectively, $A^{++}_\alpha$)
adjacent to the corresponding components~$\gamma_\alpha$
(respectively,
$\gamma^+_\alpha$) of the boundary~$\partial C$
(respectively, $\partial C^+$, see~Fig.~4).
Then $E$ extends 
to a holomorphic vector bundle over~$C^{++}$, which we also denote by~$E$.

\bigskip
\vbox{\nolineskip\rm\xsize.55\hsize%
\putt[1.05][0]{\advance\hsize-1.05\xsize
\bigskip
\centerline{Fig.4.\ \ $C\Subset C^+ \Subset C^{++}$.}
\medskip
Boundaries of $C^{++}$, $C^+$ and $C$ are marked by solid,
dotted, and dashed lines, respectively.
\smallskip}%
\epsfxsize=\xsize\epsfbox{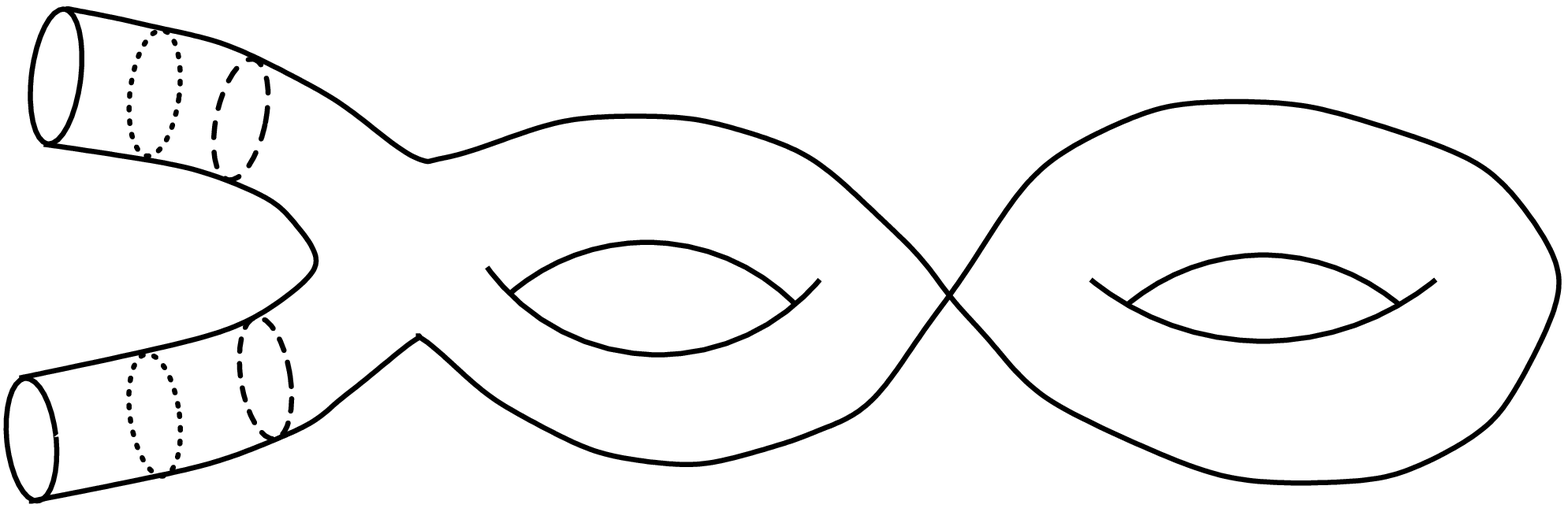}%
}
\bigskip

We consider the following sheaves on $C^{++}$:
$$
\matrix
\format\r \; \; &\c \; & \; \c \; & \; \l
\\
L^{1,p}_\loc(\,\cdot\,,E)\:&
V&\mapsto&L^{1,p}_\loc(V,E),
\\
\vspace{4pt}
L^p_\loc(\,\cdot\,,E\otimes\Lambda^{(0,1)}_{C^{++}})\:&
V&\mapsto&L^p_\loc(V,E\otimes\Lambda^{(0,1)}_{C^{++}})
\endmatrix
$$
with the sheaf homomorphism induced by the operator
$$
\overline\partial\:
L^{1,p}_\loc(V,E)\longrightarrow L^p_\loc(V,E\otimes\Lambda^{(0,1)}_{C^{++}}).
$$
The sheaves $L^{1,p}_\loc(\,\cdot\,,E)$ and
$L^p_\loc(\,\cdot\,,E\otimes 
\Lambda^{(0,1)}_{C^{++}})$, 
together with the $\overline\partial$-homomorphism,
form a fine resolution    of the (coherent) sheaf~$\calo^E$
of holomorphic sections of~$E$ over~$C^{++}$.
At smooth points of $C^{++}$ this follows from the
$L^p$-regularity of the elliptic operator~$\overline\partial$,
while 
at nodal points we use   the following argument.

Let $z\in C$ be a nodal point in the intersection
of two irreducible components, $C_i$ and~$C_j$, of~$C$.
Let $\xi_i$ (respectively, $\xi_j$) be $E$-valued
$L^p_\loc$-integrable $(0,1)$-forms,
defined in a neighborhood of~$z$ in~$C_i$ (respectively, in~$C_j$).
We find an $L^{1,p}_\loc$-solution $v_i$ (respectively,~$v_j$)
of~$\overline\partial v_i=\xi_i$
(respectively,~$\overline\partial v_j=\xi_j$).
By adding a local holomorphic section of $E$ over~$C_i$
we obtain that~$v_i(z)=v_j(z)$.
Now, the pair $(v_i,v_j)$ defines a section of~$L^{1,p}_\loc(\,\cdot\,,E)$
in a neighborhood of~$z$.
This shows that the Dolbeault lemma for the holomorphic bundle~$E$
holds also in a neighborhood of a nodal point.

This gives us the following natural isomorphisms:
$$
\aligned
\ker\Bigl(\overline\partial\:L^{1,p}_\loc(C^{++},E)\to
L^p_\loc\bigl(C^{++},E\otimes\Lambda^{(0,1)}_{C^{++}}\bigr)\Bigr)
=\roman H^0(C^{++},E),
\\
\coker\Bigl(\overline\partial\:L^{1,p}_\loc(C^{++},E)\to
L^p_\loc\bigl(C^{++},E\otimes\Lambda^{(0,1)}_{C^{++}}\bigr)\Bigr)
=\roman H^1(C^{++},E).
\endaligned
\tag 3.3
$$
We point out that similar isomorphisms exist for $C$ and~$C^+$.
Note also that we have the natural isomorphism
$$
\roman H^1(C^{++},E)=\roman H^1(C^+,E)=\roman H^1(C,E)
=\roman H^1(C_{\text{\rom{comp}}},E)
$$
induced by the restrictions
$$
\align
L^p_\loc\bigl(C^{++},E\otimes\Lambda^{(0,1)}_{C^{++}}\bigr)&\longrightarrow
L^p_\loc\bigl(C^+,E\otimes\Lambda^{(0,1)}_{C^+}\bigr)\longrightarrow
L^p_\loc\bigl(C,E\otimes\Lambda^{(0,1)}_C\bigr)
\\
&\longrightarrow
L^p_\loc\bigl(C_{\text{\rom{comp}}},E\otimes\Lambda^{(0,1)}_C\bigr).
\endalign
$$

Now, fix arbitrary  $\xi\in L^2(C^+\!\!,E\otimes\Lambda^{(0,1)}_{C^+})$
with zero cohomology class in~$\roman H^1(C^+\!\!,E)$.
We can extend $\xi$ by zero 
to an element $\ti\xi\in L^2_\loc(C^{++},E\otimes\Lambda^{(0,1)}_{C^{++}})$.
Since $[\ti\xi]_{\overline\partial}=[\xi]_{\overline\partial}=0$,
 there exists a section~$\ti v\in L^{1,2}_\loc(C^{++},E)$
such that $\overline\partial\ti v=\ti\xi$.
The restriction $v:=\ti v\big|_{C^+}$
satisfies the relations $\overline\partial v=\xi$ and $v\in L^{1,2}(C^+,E)$.
This shows that the range of the (continuous!) operator
$$
\overline\partial\:L^{1,2}(C^+,E)\longrightarrow
L^2(C^+,E\otimes\Lambda^{(0,1)}_{C^+})
\tag 3.4
$$
has finite codimension.
By Banach's open mapping theorem   this range is closed.
Further, since $L^{1,2}(C^+,E)$ is a Hilbert space,
the kernel of~(3.4) admits a direct
complement~$Q\subset L^{1,2}(C^+,E)$.
Moreover, the operator~(3.4) maps $Q$ isomorphically onto its image.
Hence the operator~(3.4) splits,
that is, there exists a continuous operator
$$
T^+\:L^2(C^+,E\otimes\Lambda^{(0,1)}_{C^+})\longrightarrow L^{1,2}(C^+,E)
\tag 3.5
$$
such that $\Im(T^+)=Q$ and for each
$\xi\in L^2(C^+,E\otimes\Lambda^{(0,1)}_{C^+})$ with
$[\xi]_{\overline\partial}=0\in\roman H^1(C^+,E)$
we have $\overline\partial(T^+\xi)=\xi$.

We define the operator
$T\:L^2(C,E\otimes\Lambda^{(0,1)}_C)\to L^{1,2}(C,E)$ as follows.
We extend  each $\xi\in L^2(C,E\otimes\Lambda^{(0,1)}_C)$
by zero 
to $\ti\xi\in L^2(C^+,E\otimes\Lambda^{(0,1)}_{C^+})$
and set $T(\xi)=T^+(\ti\xi)\big|_C$.
Then $T$ is obviously continuous and, moreover,
$$
\|T^+(\ti\xi)\|_{L^{1,2}(C^+)}\le c\cdot\|\xi\|_{L^2(C)},
$$
where the constant $c$ does not depend on $\xi$.
If $2\le p<\infty$ and $v\in L^{1,p}_\loc(C^+,E)$,
then by  the $L^p$-regularity of the elliptic $\overline\partial$-operator
(see, for example,~\cite{7}) we obtain the following interior estimate:
$$
\|v\|_{L^{1,p}(C)}\le c'\cdot
\bigl(\|v\|_{L^{1,2}(C^+)}+\|\overline\partial v\|_{L^p(C^+)}\bigr),
\tag 3.6
$$
where the constant $c'$ does not depend on $v$.
It follows that for $2\le p<\infty$
and $\xi\in L^p(C,E\otimes\Lambda^{(0,1)}_C)$
with $[\xi]_{\overline\partial}=0\in\roman H^1(C,E)$
we have the estimate
$$
\|T(\xi)\|_{L^{1,p}(C)}\le c''\cdot\|\xi\|_{L^p(C)},
$$
where the constant $c''$ does not depend on $\xi$.
This means that the operator $T$ is a splitting of~(3.2).
Hence   the operator~(3.2) has           properties
~(i) and~(ii) from Lemma~3.5.

We now return to the {\v C}ech coboundary operator~(3.1).
We fix a partition of unity
$\boldkey1=\sum_{i=1}^l\varphi_i$
subordinate to the cover $\{V_i\}_{i=1}^l$ of the curve~$C$ and
fix a cocycle
$w=(w_{ij})\in\sum_{i<j}{\calh}^{1,p}(V_{ij},E)$, where
for $i>j$  we set $w_{ij}:=-w_{ji}$.
We also set  $f_i:=\sum_j\varphi_jw_{ij}$.
Then $f_i\in L^{1,p}(V_i,E)$ and $f_i-f_j=w_{ij}$.
Consequently, $\overline\partial f_i\in L^p(V_i,E\otimes\Lambda^{(0,1)}_C)$
and $\overline\partial f_i=\overline\partial f_j$ in~$V_{ij}$.
Hence $\overline\partial f_i=\xi\big|_{V_i}$
for some well-defined section $\xi\in L^p(C,E\otimes\Lambda^{(0,1)}_C)$.
Moreover, $(w_{ij})$ and $\xi$ define the same cohomology class
$[w_{ij}]=[\xi]$ in~$\roman H^1(C,E)$.

Assume in addition that the induced cohomology class
$[w_{ij}]$ is trivial. We set $f:=T(\xi)$ and $v_i=f_i-f$.
Then $v_i\in L^{1,p}(V_i,E)$, $v_i-v_j=w_{ij}$, and
$\overline\partial v_i=\overline\partial f_i-\overline\partial f=0$.
Hence $v:=(v_i)\in\sum_{i=1}^l{\calh}^{1,p}(V_i,E)$ and
$\delta(v)=w$. It follows that the formula
$T_\delta\:w\mapsto v$ defines an operator $T_\delta$
that is a splitting of~$\delta$.
This  explicit construction shows that $T_\delta$ is continuous,
which proves Lemma~3.5.

\proclaim{Lemma 3.6}
Let $C$ be a Stein nodal curve with  piecewise smooth boundary,
and let $X$ be a complex manifold. Then
\roster
\item"(i)" ${\calh}^{1,p}(C,X)$
has a natural  complex manifold structure
with tangent space
$T_u{\calh}^{1,p}(C,X)\allowmathbreak ={\calh}^{1,p}(C,u^*TX);$

\item"(ii)"
if $C'\!\subset C\!$ is a nodal curve, then the restriction map
${\calh}^{1,p}(C,U)\to{\calh}^{1,p}(C'\!\!,U)$
is holomorphic and its differential at a point
$u\in{\calh}^{1,p}(C,U)$ is also  the restriction map
${\calh}^{1,p}(C,u^*TU)\to{\calh}^{1,p}(C',u^*TU)$,
$v\mapsto v\big|_{C'}$.
\endroster
\endproclaim

\smallskip\noindent\bf Proof. \rm  This consists of several steps.

\smallskip
{\it Step\/} 1.
Assume first that $u(C)$
lies in a coordinate chart $U\subset X$
with complex coordinates
$w=(w_1,\dots,w_n)\:U\overset\cong\to\longrightarrow U'\subset\Bbb C^n$.
Then the set ${\calh}^{1,p}(C,U)$ is an open neighborhood
of~$u$ in~${\calh}^{1,p}(C,X)$
and can be naturally identified with the set  ${\calh}^{1,p}(C,U')$,
which is an open subset of the Banach space
${\calh}^{1,p}(C,\Bbb C^n)$.
This gives us a complex Banach manifold structure
on ${\calh}^{1,p}(C,U)$
with tangent space
$T_u{\calh}^{1,p}(C,U)\cong{\calh}^{1,p}(C,\Bbb C^n)
\cong{\calh}^{1,p}(C,u^*TU)$
at~$u\in{\calh}^{1,p}(C,U)$.

Note that if $u_t$, $t\in[0,1]$, is a $C^1$-curve in
${\calh}^{1,p}(C,U)$, then the tangent vector
$v\in{\calh}^{1,p}(C,u^*TU)$ to $u_t$ in $u_0$
is given by the formula
$v(z)=\dfrac{\partial u}{\partial t}(z)\in T_{u(z)}U$.
This last formula does not depend
on the choice of complex coordinates
$w=(w_1,\dots,w_n)\:U\to\Bbb C^n$ in~$U$.
This has  the following two consequences.

Firstly, {\it the complex structure on ${\calh}^{1,p}(C,U)$
does not depend on the choice of complex coordinates
$w=(w_1,\dots,w_n)\:U\to\Bbb C^n$ in~$U$}.
Secondly, {\it $C$ has  property~(ii) from the
statement of the lemma}. 

Thus, Lemma~3.6 is proved in the case when $u(C)$ lies in a coordinate chart.

\smallskip
{\it Step\/} 2.
Assume that $u_0\in{\calh}^{1,p}(C,X)$ is fixed
and there exists a finite cover $\{V_i\}_{i=1}^l$ of~$C$
such that, firstly, the assumptions of Lemma~3.5 are satisfied
and, secondly, Lemma~3.6 holds  for each $V_i$
(for example, assume that each $u_0(V_i)$ lies in a chart~$U_i\subset X$).

We set  $V_{ij}:=V_i\cap V_j$ and choose balls
$$
B_{ij}\subset{\calh}^{1,p}(V_{ij},u_0^*TX)
\cong T_{u_0}{\calh}^{1,p}(V_{ij},X),
$$
such that there exist biholomorphisms
$\psi_{ij}\:B_{ij}\overset\cong\to\longrightarrow
B'_{ij}\subset{\calh}^{1,p}(V_{ij},X)$
with $\psi_{ij}(0)=u_0\big|_{V_{ij}}$ and
$d\psi_{ij}(0)=\id\:{\calh}^{1,p}(V_{ij},u_0^*TX)
\to{\calh}^{1,p}(V_{ij},u_0^*TX)$.
Then we  choose balls $B_i\subset{\calh}^{1,p}(V_i,X)$ such that
 $u_0\big|_{V_i}\in B_i$ and $u_i\big|_{V_{ij}}\in B'_{ij}$
for each~$u_i\in B_i$

This defines   holomorphic maps
$\varphi_{ij}\:B_i\to B_{ij}\subset{\calh}^{1,p}(V_{ij},u_0^*TX)$
such that
$\varphi_{ij}\:u_i\mapsto\psi_{ij}^{-1}(u_i|_{V_{ij}})$,
and a    holomorphic map
$$
\matrix
\Phi\:&\prod_{i=1}^lB_i&\longrightarrow&\sum_{i<j}{\calh}^{1,p}(V_{ij},u_0^*TX),
\\
\Phi\:&(u_i)_{i=1}^l&\mapsto&\varphi_{ij}(u_i)-\varphi_{ji}(u_j).
\endmatrix
$$
It is easy to see that the map~$\Phi$
gives us a  condition for  the compatibility
of local holomorphic maps
$u_i\:V_i\to X$;
namely, $(u_i)_{i=1}^l\in\prod_{i=1}^lB_i$
defines a holomorphic map $u\:C\to X$
if and only if~$\Phi(u_i)=0$.
Furthermore, the differential $d\Phi$ in $(u_0|_{V_i})$
is equal to the {\v C}ech coboundary operator~(3.1).
Since $C$ is Stein, it follows that~$\roman H^1(C,u_0^*TX)=0$.
Using Lemma~3.5 and the implicit function theorem
we conclude that parts~(i) and~(ii) of Lemma~3.6
hold in a neighborhood of the map
$u_0\in{\calh}^{1,p}(C,X)$.

\smallskip
{\it Step\/} 3.
Applying step~2 sufficiently many times
one can show that for each Stein nodal curve~$C$
and  each~$u\in{\calh}^{1,p}(C,X)$
parts~(i) and~(ii) of Lemma~3.6
hold in a neighborhood of the map~$u$.
For example, if $C$ is the annulus~$A_{r,R}$,
then we can cover it with narrow annuli~$A_{r_i,R_i}$,
$0<\nomathbreak R_i-\nomathbreak r_i\ll\nomathbreak 1$,
and then cover each $A_{r_i,R_i}$  with sectors
$V_{ij}=\{z=\rho e^{\isl\theta}\in\Bbb C:
r_i<\rho<R_i,\allowmathbreak \alpha_j<\theta<\beta_j\}$, where
 $0<\beta_j-\alpha_j\ll 1$.
We leave the details to the reader.

One of the difficulties 
in the construction of {\it holomorphic\/}
families of stable curves
is that the moduli space of holomorphic structures
on a non-compact Riemann surface~$\Sigma$
does not have a natural complex structure
and, moreover, its real dimension may be odd.
For example, if $\Sigma$ is an annulus,
then it is biholomorphic to the standard annulus~$A_{r,1}$
for some unique~$r\in(0,1)$,
and therefore the corresponding moduli space is the interval~$(0,1)$.
In general, if $\Sigma$ is of genus~$g$
and has $k$ boundary components,
then the real dimension of the moduli space
is equal to~$d=6g-6+3k$,
except for the four cases
when $\Sigma$ is either a sphere ($g=0$, $k=0$), or a torus~($g=1$, $k=0$),
or a disc~($g=0$, $k=1$), or an annulus~($g=0$, $k=2$)
(see, for example,~\cite{8}).
Note that these are the only cases
when the dimension of the group of holomorphic automorphisms
of the corresponding complex curve~$(\Sigma,J)$ is positive.

The problem  may be fixed 
by an introduction of  $k$ additional parameters,
namely, by fixing $k$ marked points,
one on each boundary component.
Let $A$ be an annulus with boundary circles~$\gamma_0$ and~$\gamma_1$,
and let $X$ be a complex manifold.

\proclaim{Theorem 3.7}
There exist complex Banach manifolds
$\calm(A,X)$ and ${\calc}(A,X)$,
a~ho\-lo\-morphic projection
$\pi_{\!\calc\!}\:\!{\calc}(A,X)\to\calm(A,X)$,
and holomorphic maps $\ev\!\:\!{\calc}(A,X)\to X$,
$z_1\:{\calc}(A,X)\to\Delta$, $z_2\:{\calc}(A,X)\to\Delta$ and
$\lambda_{\!\calm\!}\:\calm(A,X)\to\Delta$
with the following properties\rom:
\roster
\item"(i)" for each $y\in\calm(A,X)$ the fiber
$C_y:=\pi_{\!\calc\!}^{-1}(y)$ is a nodal curve
parametrized by the annulus~$A;$
moreover, the map $(z_1,z_2)\:C_y\to\Delta^2$
is a biholomorphism onto the curve
$\{(z_1,z_2)\in\Delta^2:z_1\cdot\nomathbreak z_2
=\lambda_{\!\calm\!}(y)\}$;
in~particular, $C_y$ is either a standard node
(if $\lambda_{\!\calm\!}(y)=0$),
or a holomorphic annulus $\{|\lambda_{\!\calm\!}(y)|<|z_1|<1\};$

\item"(ii)" the  diagram
$$
\CD
{\calc}(A,X)@>(\ev,z_1,z_2)>>X\times\Delta^2
\\
@VV\pi_{\!\calc\!}V @VV\lambda=z_1\cdot z_2V
\\
\calm(A,X)@>\lambda_{\!\calm\!}>>\Delta
\endCD
\tag 3.7
$$
is commutative\rom;
moreover, for each $y\in\calm(A,X)$ the restriction
$\ev\big|_{C_y}$ belongs to ${\calh}^{1,p}(A_a,X)$ with
$a=\lambda_{\!\calm\!}(y)$
and the maps
$\ev_1\:y\in\calm(A,X)\mapsto\ev\big|_{C_y}(z_1^{-1}(1))$ and
$\ev_2\:y\in\calm(A,X)\mapsto\ev\big|_{C_y}(z_2^{-1}(1))$
are holomorphic\rom;

\item"(iii)" let $C$ be an annulus or a node with smooth boundary
$\partial C=\gamma_1\sqcup\gamma_2$,
$p_i\in\gamma_i$ marked points, 
and $u\:C\to X$ 
a holomorphic map of class $L^{1,p}$\rom;
then there exists {\it  unique\/}  $y\in\calm(A,X)$
and a {\sl  unique\/} biholomorphism $\varphi\:C\to C_y$
such that $\ev\circ\varphi=u\:C\to X$ and
$z_i\circ\varphi(p_i)=1\in\overline\Delta$\rom;
in other words,
$\calm(A,X)$ parameterizes  holomorphic maps into~$X$
of annuli and nodes with marked boundary points\rom;

\item"(iv)" if the  commutative diagram
$$
\CD
{\calz}@>(\ev^{\calz},\ti z_1,\ti z_2)>>X\times\Delta^2
\\
@VV\pi_{\calz}V @VV\lambda=\ti z_1\cdot\ti z_2V
\\
{\calw}@>\lambda_{\calw}>>\Delta
\endCD
\tag 3.8
$$
of complex spaces ${\calw}$ and ${\calz}$
and holomorphic maps
has properties \rom{(i)} and \rom{(ii)},
in particular, if the fibers
${\calz}_w:=\pi_{\calz}^{-1}(w)$
are nodal curves with induced maps
$f_w:=\ev^{\calz}\big|_{{\calz}_w}\in{\calh}^{1,p}({\calz}_w,X)$,
then the two diagrams~\thetag{3.7} and~\thetag{3.8}
can be completed in a unique way
to 
the following commutative diagram\rom:
$$
\CD
{\calz}@>\wt F>>{\calc}(A,X)@>(\ev,z_1,z_2)>>X\times\Delta^2
\\
@VV\pi_{\calz}V @VV\pi_{\!\calc\!}V @VV\lambda=z_1\cdot z_2V
\\
{\calw}@>F>>\calm(A,X)@>\lambda_{\!\calm\!}>>\phantom{,}\Delta,
\endCD
\tag 3.9
$$
where $\lambda_{\!\calm\!}\circ F=\lambda_{\calw}$ and
$(\ev,z_1,z_2)\circ\wt F=(\ev^{\calz},\ti z_1,\ti z_2)$\rom;

\item"(v)" the differential
$d\lambda_{\!\calm\!}\:T_y\calm(A,X)\to
T_{\lambda_{\!\calm\!}(y)}\Delta\cong\Bbb C$
is non-degenerate at each  point~$y$ in $\calm(A,X)$,
and for each $a\in\Delta$ the fiber $\lambda_{\!\calm\!}^{-1}(a)$
is naturally isomorphic to the manifold ${\calh}^{1,p}(A_a,X)$,
where $A_a$ is the curve $\{(z_1,z_2)\in\Delta^2:z_1\cdot z_2=a\}$\rom;
in~particular, for each $y\in\calm(A,X)$
one  obtains a biholomorphism~$C_y\cong A_{\lambda_{\!\calm\!}(y)}$
and the following natural exact sequence\rom:
$$
\minCDarrowwidth{7mm}
\CD
0@>>>{\calh}^{1,p}(C_y,u^*TX)@>\iota_y>>T_y\calm(A,X)
@>d\lambda_{\!\calm\!}(y)>>\Bbb C@>>>0.
\endCD
$$
\endroster
\endproclaim

\smallskip\noindent\bf Proof. \rm 
Let $(A,p_1,p_2)$ be a smooth annulus
with a marked point on each boundary component~$\gamma_i\cong S^1$,
and let $J$ be a complex structure on~$A$.
We know that $(A,J)$ is biholomorphic to
one of the annuli~$A_{r,1}=\nomathbreak \{r<\nomathbreak |z|<\nomathbreak 1\}$.
It is easy to see that there exists only one
isomorphism $\psi\:(A,J)\to A_{r,1}$
that extends smoothly to a diffeomorphism
$\psi\:\overline A\to\overline A_{r,1}$ with~$\varphi(p_1)=1$.
Set~$a:=\varphi(p_2)$.
It is now evident that there exists a unique biholomorphism
$\varphi\:(A,J)\to A_a:=\{(z_1,z_2)\in\Delta^2:z_1\cdot z_2=a\}$
such that $\varphi(p_1)=1$ and $\varphi(p_2)=a$.

Thus, the map $\lambda\:\Delta^2\to\Delta$,
$\lambda(z_1,z_2)=z_1\cdot z_2$, with  fiber $A_a$
over $a\in\nomathbreak \Delta$
forms 
the holomorphic moduli space of annuli
with marked points on boundary components
completed by   the standard node at~$a=0$.
If $a\ne0$, then the coordinate functions $z_i$, $i=1,2$,
define 
an embedding of each $A_a$ in~$\Bbb C$ such
that $\gamma_i$ becomes the outer unit circle.
As $a\to0$, the annuli $A_a$ degenerate into the standard node,
and each $z_i$ becomes            
the standard coordinate function on the corresponding component of the node.

\smallskip\noindent\bf Remark. \rm 
In what follows, we denote by $A_a$ an annulus (or a node)
{\it with marked points on its  boundary~$\partial A_a$}
and with  coordinate functions $z_1$ and~$z_2$
defined as above.
\endremark

Fix $r$, $0<r<1$. For $|a|<r$ we define the maps
$\zeta^a_1,\zeta^a_2\:A_{r,1}\to\nomathbreak A_a$
by the formulae $\zeta^a_1(z):=z$ and~$\zeta^a_2(z):=a/z$,
so that the $\zeta^a_i$ are the reciprocals
of the coordinate functions~$z_i$.
We consider the following map:
$$
\matrix
\Psi_r\:&\coprod_{|a|<r}{\calh}^{1,p}(A_a,X)&\longrightarrow&
{\calh}^{1,p}(A_{r,1},X)\times{\calh}^{1,p}(A_{r,1},X)\times\Delta(r),
\\
\Psi_r\:&u\in{\calh}^{1,p}(A_a,X)&\mapsto&(u\circ\zeta^a_1,u\circ\zeta^a_2,a).
\endmatrix
$$
It is easy to see that $\Psi_r$ is holomorphic on each
space ${\calh}^{1,p}(A_a,X)$ and  the image of
$\coprod_{|a|<r}{\calh}^{1,p}(A_a,X)$
consists of  triples $(u_1,u_2,a)$ such
that each map $u_i\in{\calh}^{1,p}(A_{r,1},X)$
extends 
to a map $u_i\in{\calh}^{1,p}(A_{|a|,1},X)$ and~$u_2(z)=u_1(a/z)$.
Hence $\Psi_r$ is injective in
$\coprod_{|a|<r}{\calh}^{1,p}(A_a,X)$ and this image is closed.
We consider the topology induced by the maps~$\Psi_r$
on the disjoint union
$\calm(A,X):=\coprod_{a\in\Delta}{\calh}^{1,p}(A_a,X)$,
Clearly, it  is compatible with the topology on each fiber
${\calh}^{1,p}(A_a,X)$.

Our aim is to construct an appropriate holomorphic structure
on~$\calm(A,X)$
compatible with the holomorphic structures
on the fibers ${\calh}^{1,p}(A_a,X)$
and with the topology on $\calm(A,X)$
introduced above.

We start with the special case of $X=\Bbb C^n$.
It is easy to see that for $a\ne0$ each function
$f\in{\calh}^{1,p}(A_a,\Bbb C^n)$
can be uniquely expanded in the Laurent series
$f(z_1)=\sum_{i=-\infty}^\infty c_iz_1^i$.
We set $f^+(z_1):=\sum_{i=0}^\infty c_iz_1^i$ and
$f^-(z_1):=\sum_{i=-\infty}^0c_iz_1^i$.
It is also convenient to regard  ~$f^-$
as a function of the variable $z_2=a/z_1$
with~$f^-(z_2)=\sum_{i=0}^\infty c_{-i}(z_2/a)^i$.
We have
$f^+\in{\calh}^{1,p}(\{|z_1|<1\},\Bbb C^n)$,
$f^-\in{\calh}^{1,p}(\{|z_2|<1\},\Bbb C^n)$,
$f^+(0)=f^-(z_2=0)=c_0$, and
$f(z_1)=f^+(z_1)+f^-(a/z_1)-c_0$,
so that the pair $(f^+,f^-)$
defines a holomorphic function
$\widehat f\in{\calh}^{1,p}(A_0,\Bbb C^n)$.
The resulting canonical isomorphisms
${\calh}^{1,p}(A_a,\Bbb C^n)\cong{\calh}^{1,p}(A_0,\Bbb C^n)$
define
a structure of a trivial Banach bundle over~$\{|a|<1\}$,
and therefore a structure of a Banach manifold,
 on~$\coprod_{|a|<1}{\calh}^{1,p}(A_a,\Bbb C^n)$.

Now, the map
$\Psi_r\:\coprod_{|a|<r}{\calh}^{1,p}(A_a,\Bbb C^n)
\to{\calh}^{1,p}(A_{r,1},\Bbb C^n)^2\times\Delta(r)$
is holomorphic.
If $U$ is an open subset of $\Bbb C^n$,
then $\calm(A,U)$ is also open in $\calm(A,\Bbb C^n)$
and therefore inherits a   holomorphic structure.
The natural projection
$\lambda_{\!\calm\!}\:\calm(A,U)\to\Delta$
is now  holomorphic.
It is easy to see that the differential
$d\lambda_{\!\calm\!}$ is non-degenerate.
Hence  we can define
the {\it universal family of curves\/} ${\calc}(A,U)$
as the fibered product $\calm(A,U)\times_\Delta\Delta^2$
with respect to the maps
$\lambda_{\!\calm\!}\:\calm(A,U)\to\Delta$ and
$\lambda\:\Delta^2\to\Delta$.
This is a holomorphic Banach manifold
because $d\lambda_{\!\calm\!}$ is non-degenerate.

Let $\pi_{\!\calc\!}\:{\calc}(A,U)\to\calm(A,U)$ be
the natural projection. Then the fiber $C_y$ over
$y\in\calm(A,U)$ is biholomorphic to~$A_a$ with
$a=\lambda_{\!\calm\!}(y)$.
The natural projection of ${\calc}(A,U)$ onto~$\Delta^2$
induces  $\Delta$-valued holomorphic functions $z_1$~and~$z_2$
on ${\calc}(A,U)$
 that have property~(i) from  Theorem~3.7.

Assume now that $a\ne0\in\Delta$ and let
$f\in{\calh}^{1,p}(A_a,\Bbb C^n)$.
We can represent $f$ in the form $f(z_1)=f^+(z_1)+f^-(a/{z_1})-f_0$,
where  $f^\pm\in{\calh}^{1,p}(\Delta,\Bbb C^n)$ and
$f^+(0)=f^-(0)=f_0$.
In a similar way, for each
$f\in{\calh}^{1,p}(A_0,\Bbb C^n)$
we have $f=(f^+,f^-)$,
where we also have  $f^\pm\in{\calh}^{1,p}(\Delta,\Bbb C^n)$
with~$f^+(0)=f^-(0)=f_0$.
Consider 
the holomorphic function
$\ti f(z_1,z_2):=f^+(z_1)+f^-(z_2)-f_0$,
$\ti f\in{\calh}^{1,p}(\Delta^2,\Bbb C^n)$, and
define the maps
$\ti\ev_a\:{\calh}^{1,p}(A_a,\Bbb C^n)\times\Delta^2\to\Bbb C^n$
by the formula $\ti\ev(f,z_1,z_2):=\ti f(z_1,z_2)$.
It is easy to see that the $\ti\ev_a$ 
define 
a holomorphic map
$\ti\ev\:\calm(A,\Bbb C^n)\times\Delta^2\to\Bbb C^n$.
Let $\ev$ be     the restriction of this map to
${\calc}(A,\Bbb C^n)\subset\calm(A,\Bbb C^n)\times\Delta^2$.
We leave it to the reader to verify that
the assertions of Theorem~3.7   hold for
$\calm(A,\Bbb C^n)$, ${\calc}(A,\Bbb C^n)$,
$\ev\:{\calc}(A,\Bbb C^n)\to\Bbb C^n$, and
$z_{1,2}\:{\calc}(A,\Bbb C^n)\to\Delta$.

It follows now that for $U\subset\Bbb C^n$ we can define
the map $\ev\:{\calc}(A,U)\to U$
as a mere  restriction~$\ev\:{\calc}(A,\Bbb C^n)\to\Bbb C^n$.
The assertions of Theorem~3.7 hold again.
In particular, if $G\:U\to U'\subset\Bbb C^n$
is biholomorphic, then the natural bijections
$\calm(A,U)\overset\cong\to\longrightarrow\calm(A,U')$ and
${\calc}(A,U)\overset\cong\to\longrightarrow{\calc}(A,U')$
are biholomorphisms.
This means that the holomorphic structure
on~$\calm(A,U)$ does not depend on the embedding~$U\subset\Bbb C^n$.

Let $C=A_{r,1}$ be an annulus
and let $u\:C\to X$ be a holomorphic embedding.
Then $du\:TC\to u^*TX$ is an embedding of holomorphic bundles over~$C$,
which allows us  to define the holomorphic normal bundle
as the quotient bundle~$N_C:=u^*TX/TC$.
Since $C$ is Stein,  $N_C$ is holomorphically trivial.
We fix a holomorphic frame
$\sigma_1,\dots,\sigma_{n-1}\in{\calh}^{1,p}(C,N_C)$,
$n=\dim_{\Bbb C}X$,
and its lifting
 $\ti\sigma_1,\dots,\ti\sigma_{n-1}\in{\calh}^{1,p}(C,u^*TX)$.
Let $B^{n-1}(r)$ be  the ball of radius~$r$
in the space $\Bbb C^{n-1}$ with coordinate functions~$w=(w_1,\dots,w_{n-1})$.
By Lemma~3.6  there exists a holomorphic map
$\Psi\:C\times B^{n-1}(r)\to X$ such that
$\dfrac{\partial\Psi}{\partial w_i}\bigg|_{z\in C,w=0}=\sigma_i(z)$.
Hence $\Psi$ is biholomorphic in a neighborhood of~$C\equiv C\times\{0\}$.
In~particular, if $r$ is sufficiently small,
then the image $U:=\Psi(C\times B^{n-1}(r))$
is a local chart with coordinates~$(z,w_1,\dots,w_{n-1})$.
Note that the image $u(C')$ of each  smaller annulus~$C'\Subset C$
 also lies in~$U$.

To complete the proof of the theorem it  remains only
to consider the general case when $C\cong A_a$ is arbitrary
and $u\:C\to X$ is a holomorphic map.
If $a=0$, then $C$ is a node and therefore
there exists a neighborhood $V_0$ of the nodal point
such that $u(\overline V_0)$ lies in some coordinate chart in~$X$.
If $a\ne0$ and $u(\overline C)$ does not lie in a coordinate chart in~$X$,
then $u$ is not constant. Hence the map $u$
is an embedding in a neighborhood of the circle
$S^1_r:=\{|z_1|=r\}\subset C\cong\{|a|<|z_1|<1\}$ for some $|a|<r<1$.

In any case we obtain a cover $\{V_0,V_1,V_2\}$
of the curve $C=\{(z_1,z_2)\in\Delta^2:z_1\cdot z_2=a\}$
that has  the following form:
$$
\align
V_1&=\{(z_1,z_2)\in C:r_1<|z_1|<1\},
\\
V_2&=\{(z_1,z_2)\in C:r_2<|z_2|<1\},
\\
V_0&=\{(z_1,z_2)\in C:|z_1|<R_1,|z_2|<R_2\}.
\endalign
$$
Here $0<r_1<R_1<1$, $0<r_2<R_2<1$,
$r_1\cdot R_2>|a|<r_2\cdot R_1$
and $V_0$ has the following property:    ~$u(\overline V_0)$
lies in a coordinate chart~$U$ in~$X$.
\bigskip
\vbox{\xsize.35\hsize\nolineskip\rm
\putm[.15][0.]{|z_2|}%
\vskip.06\xsize
\putm[.095][.07]{1}%
\putm[.07][.32]{R_2}%
\putm[.09][0.5]{r_2}%
\putm[.345][.32]{V_2}%
\putm[.33][.68]{V_0}%
\putm[.66][.64]{V_1}%
\putm[1.01][.82]{|z_1|}
\putm[.48][.9]{r_1}%
\putm[.66][.9]{R_1}%
\putm[.92][.9]{1}%
\putt[1.2][0]{\advance\hsize-1.4\xsize%
\bigskip
\centerline{Fig.5.\ \ Covering $\{V_0, V_1, V_2\}$.}
\bigskip
\parindent=0pt
The curve $A_a$ is drawn by solid line as a piece of a hyperbola,
the elements of the covering $V_0$, $V_1$, and $V_2$ by punctured line.
}
\epsfxsize=\xsize\epsfbox{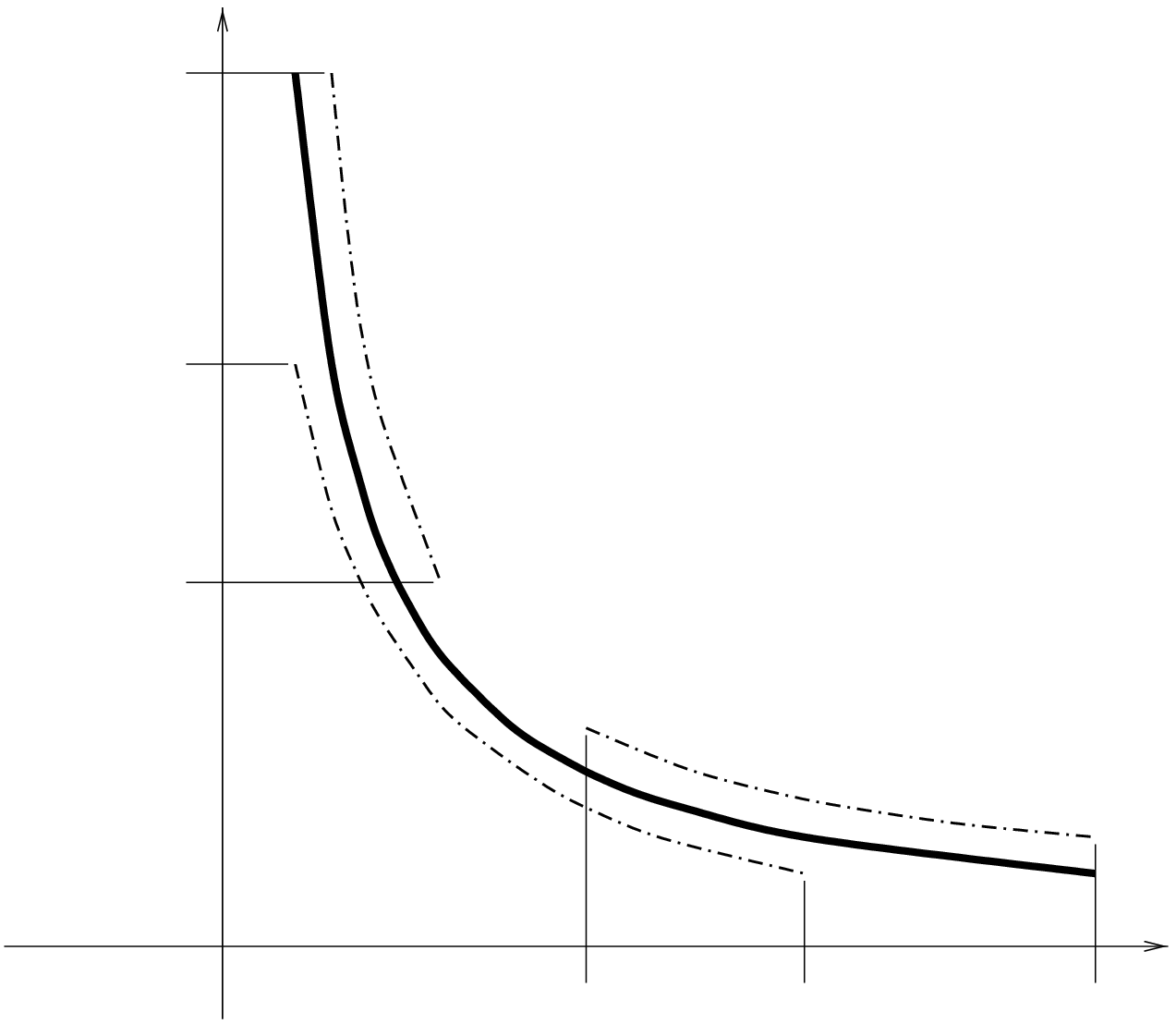}
}

\bigskip
We fix the coordinate
function $z_1$ in $V_1$ and the coordinate function $z_2$ in~$V_2$.
Next we introduce the new coordinates
$\ti z_1:=z_1/R_1$ and $\ti z_2:=z_2/R_2$ on $V_0$
and fix the marked points
$\ti p_1:=R_1$ and $\ti p_2:=a/R_2$.
We set $\ti a:=a/(R_1R_2)$. Then
$V_0\cong A_{\ti a}=\{(\ti z_1,\ti z_2)\in\Delta^2:
\ti z_1\cdot\nomathbreak \ti z_2=\ti a\}$
and $\ti z_i(\ti p_1)=1$.
It is easy to see that changing  
the complex parameter $a=\nomathbreak z_1z_2$ on~$C$, which
parameterizes the holomorphic structures on an  annulus
with marked points on the boundary,
can be reduced to changing  
 a similar parameter in~$V_0$.
Namely, let
$C'\cong A_{a'}=\{(z_1,z_2)\in\Delta^2:z_1\cdot z_2=a'\}$
be a result of a (small) deformation of~$a$.
We set $\ti a{}':=a'/(R_1R_2)$ and regard $C'$ as
a result of gluing together the complex curves $V_0'$, $V_1$ and~$V_2$
 defined in~$C'$
by the same equations
$V_0'=\{(z_1,z_2)\in C':|z_1|<r_1,\ |z_2|<r_2\}$ and
$V_i=\{r_i<|z_i|<1\}$ for $i=1,2$,
with  coordinate functions
$(z_1,z_2)$ satisfy now the new relation
$z_1\cdot z_2=a'$, so that $V_0'\cong A_{\ti a{}'}$.

Thus, in an appropriate small neighborhood
${\calw}\subset\calm(A,X)$
of the curve $(C,u)$ over $X$
we have the following holomorphic map:
$$
\matrix
\Theta\:&{\calw}&\longrightarrow&
\calm(A,U)\times{\calh}^{1,p}(V_1,X)\times{\calh}^{1,p}(V_2,X);
\\
\vspace{5pt}
\Theta\:&(C',u')&\mapsto&\bigl((A_{\ti a{}'},u|_{V_0'}),
u|_{V_1},u|_{V_2}\bigr).
\endmatrix
$$
On the other hand  the collection of maps~$u_0\:V_0'\to X$,
$u_1\:V_1\to X$ and $u_2\:V_2\to X$
defines a map $u'\:C'\to X$ if and only if $u_1$ coincides with~$u_0$
on~$W_1:=V_1\cap V_0'$ and $u_2$ coincides with~$u_0$
on~$W_2:=V_2\cap V_0'$.
We point out that the domains $W_i\subset V_i$
do not change in  the deformation of the complex
structure of~$C$.
The  construction of gluing  from step~2
in the proof of Lemma~3.6
completes the proof of the theorem.

Now we can complete the proof of Theorem~3.4.

\smallskip\noindent\bf Proof. \rm{Proof of Theorem~\rom{3.4}}
Let $(C_0,u_0)$ be a stable curve over a complex manifold~$X$
with  parametrization~$\sigma_0\:\Sigma\to C_0$.
We use  Proposition~3.2 to fix a cover
$\{V_\alpha,V_{\alpha\beta}\}$ of~$\Sigma$
with properties~(i)--(vi). In particular, there exist biholomorphisms
$\varphi^0_{\alpha\beta}\:\sigma_0(V_{\alpha\beta})\to
A_{\lambda^0_{\alpha\beta}}$.
It follows from properties~(i)--(vi)
that for each collection ${\bold\lambda}:=(\lambda_{\alpha\beta})$ that is
sufficiently close  to~${\bold\lambda}^0:=(\lambda^0_{\alpha\beta})$,
there exists a nodal curve~$C$
with  parametrization~$\sigma\:\Sigma\to C$
such that the properties~(i)--(vi) still hold
and there exist biholomorphisms
$\varphi_{\alpha\beta}\:\sigma(V_{\alpha\beta})\to A_{\lambda_{\alpha\beta}}$.
In particular,  complex structures on each
$\sigma(W_{\alpha\beta})$ do not change.
Moreover, we can choose holomorphic coordinate functions~$z_1$ and~$z_2$
in~$V_{\alpha\beta}$ such that  $z_1\cdot z_2\equiv\lambda_{\alpha\beta}$,
 and both~$z_1\big|_{W_{\alpha,\beta}}$ and $z_2\big|_{W_{\beta,\alpha}}$
do not change under a variation of ~$\lambda_{\alpha\beta}$.
This means that the disc
$\Delta_{\alpha\beta}:=\{\lambda_{\alpha\beta}:
|\lambda^0_{\alpha\beta}-\lambda_{\alpha\beta}|\le\epsi\}$
parameterizes a {\it holomorphic\/} family of curves
of the form~$\sigma(V_{\alpha\beta})$.

The deformation of  complex structure on~$V_{\alpha\beta}$ 
is shown in Fig.~6.

\medskip
\vbox{\nolineskip\xsize.548\hsize%
\putm[.01][.38]{\underbrace{\hskip.35\xsize}_{V_\alpha}}
\putm[.54][.44]{\underbrace{\hskip.46\xsize}_{V_\beta}}
\putm[.27][.16]{\overbrace{\hskip.36\xsize}^{V_{\alpha\beta}}}
\putm[.27][.3]{\underbrace{\hskip.09\xsize}_{W_{\alpha\beta}}}
\putm[.54][.3]{\underbrace{\hskip.09\xsize}_{W_{\beta\alpha}}}
  \putt[1.05][0]{\advance\hsize-1.05\xsize\parindent=0pt%
  \centerline{Fig.\ 6. }
 \smallskip
On this picture one sees $V_{\alpha\beta}$ together with adjoint
$V_\alpha$ and $V_\beta$. One can suppose that, under varying of
$\lambda _{\alpha \beta}$, the complex structure varies only in shadowed
domain.
}
\epsfxsize=\xsize\epsfbox{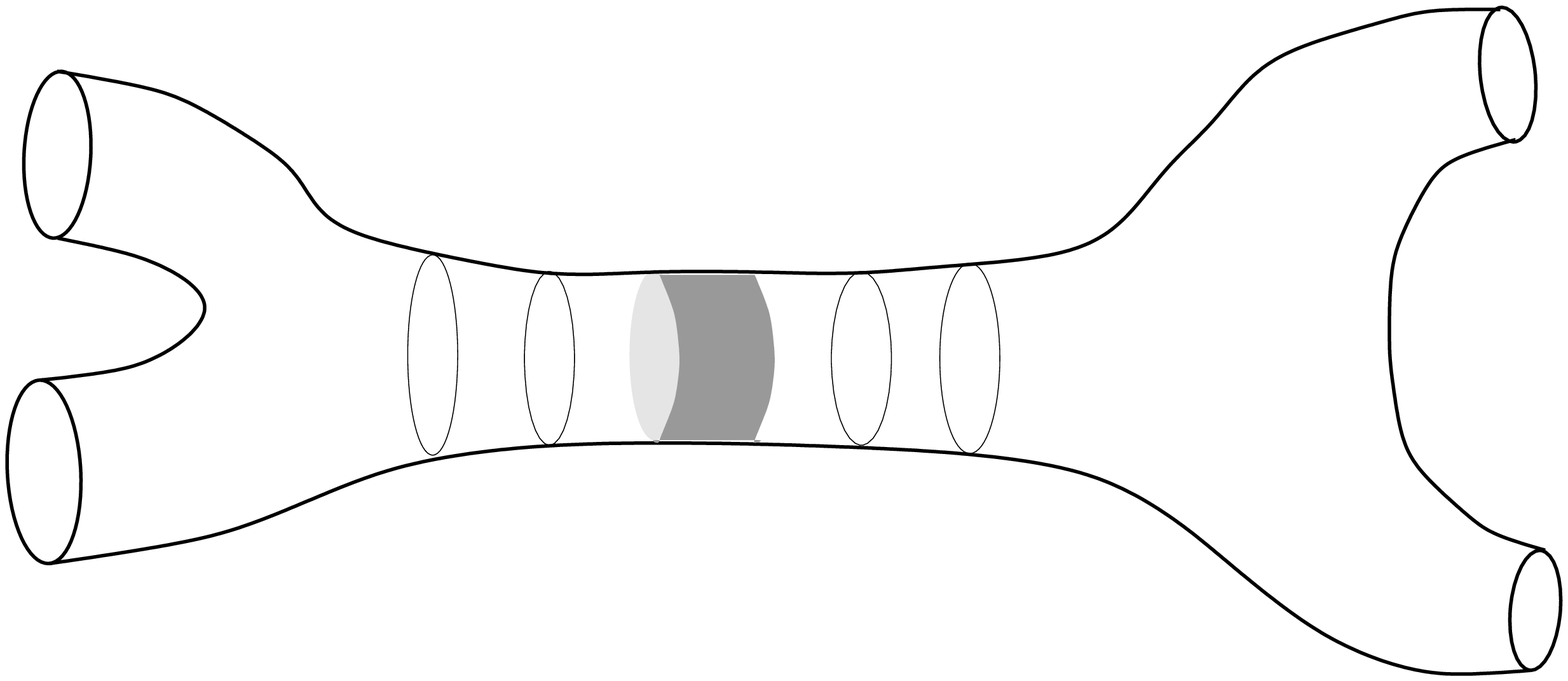}
}

\bigskip
Let $N\in\Bbb N$ be the number of elements of the covering
of the type~$V_{\alpha\beta}=V_{\beta\alpha}$.
Then, for sufficiently small $\epsi>0$ the polydisc
$$
\Delta^N_\lambda:=\{{\bold\lambda}:=(\lambda_{\alpha\beta}):
|\lambda^0_{\alpha\beta}-\lambda_{\alpha\beta}|\le\epsi\}
$$
parameterizes a holomorphic family of nodal curves
$\{C_{\bold\lambda}\}_{{\bold\lambda}\in\Delta^N_\lambda}$.
Each curve $C_{\bold\lambda}=\sigma(\Sigma)$
is obtained by patching  together the pieces~$\sigma(V_\alpha)$ 
and~$\sigma(V_{\alpha\beta})$.
We set $\calm:=\bigcup_{{\bold\lambda}\in\Delta^N_\lambda}
{\calh}^{1,p}(C_{\bold\lambda},X)$.
The pieces~$\sigma(V_\alpha)$ do~not contain nodal points,
and the complex structures on~$\sigma(V_\alpha)$
are constant and independent of ~${\bold\lambda}=(\lambda_{\alpha\beta})$.
Furthermore, we have the natural isomorphisms
$\sigma(V_{\alpha\beta})\cong A_{\lambda_{\alpha\beta}}$.
Hence the following map is well defined:
$$
\matrix
\Theta\:&\calm&\longrightarrow&\prod_\alpha{\calh}^{1,p}(V_\alpha,X)
\times\prod_{\alpha\beta}\calm(V_{\alpha\beta},X);
\\
\Theta\:&(C_{\bold\lambda},u)&\mapsto&\bigl(u|_{V_\alpha},
(\lambda_{\alpha\beta},u|_{V_{\alpha\beta}})\bigr).
\endmatrix
$$
It is easy to see that a collection
$\bigl(\!u_\alpha,(\lambda_{\alpha\beta},u_{\alpha\beta})\!\bigr)
\in\prod_\alpha{\calh}^{1,p}(V_\alpha,X)\times
\prod_{\alpha\beta}\calm(V_{\alpha\beta},X)$
belongs to     $\Theta(\calm)$
if and only if the gluing conditions
$u_\alpha\big|_{W_{\alpha,\beta}}=u_{\alpha\beta}\big|_{W_{\alpha,\beta}}$
are satisfied for all pairs~$(\alpha,\beta)$.

We shall repeat the gluing procedure of step~2 in Lemma~3.6.
To this end we choose the balls
$$
\align
B_\alpha&\subset{\calh}^{1,p}(V_\alpha,u_0^*TX)
\cong T_{u_0}{\calh}^{1,p}(V_\alpha,X),
\\
B_{\alpha\beta}&\subset{\calh}^{1,p}(V_{\alpha\beta},u_0^*TX)\oplus\Bbb C
\cong T_{u_0}\calm(V_{\alpha\beta},X),
\\
B'_{\alpha,\beta}&\subset{\calh}^{1,p}(W_{\alpha,\beta},u_0^*TX)
\cong T_{u_0}{\calh}^{1,p}(W_{\alpha,\beta},X)
\endalign
$$
such  that there exist biholomorphisms
$$
\align
\psi_\alpha&\:B_\alpha\overset\cong\to\longrightarrow\psi_\alpha(B_\alpha)
\subset{\calh}^{1,p}(V_\alpha,X),
\\
\psi_{\alpha\beta}&\:B_{\alpha\beta}\overset\cong\to\longrightarrow
\psi_{\alpha\beta}(B_{\alpha\beta})\subset\calm(V_{\alpha\beta},X),
\\
\psi'_{\alpha,\beta}&\:B'_{\alpha,\beta}\overset\cong\to\longrightarrow
\psi'_{\alpha,\beta}(B'_{\alpha,\beta})\subset{\calh}^{1,p}(W_{\alpha,\beta},X)
\endalign
$$
with the following properties:
$$
\alignat2
\psi_\alpha(0)&=u_0\big|_{V_\alpha},&\quad
d\psi_\alpha(0)&=\id\:T_{u_0}{\calh}^{1,p}(V_\alpha,X)\to
T_{u_0}{\calh}^{1,p}(V_\alpha,X),
\\
\psi_{\alpha\beta}(0)&=u_0\big|_{V_{\alpha\beta}},&\quad
d\psi_{\alpha\beta}(0)&=\id\:T_{(\lambda^0_{\alpha,\beta},u_0)}
\calm(V_\alpha,X)\to
T_{(\lambda^0_{\alpha,\beta},u_0)}\calm(V_\alpha,X),
\\
\psi'_{\alpha,\beta}(0)&=u_0\big|_{W_{\alpha,\beta}},&\quad
d\psi'_{\alpha,\beta}(0)&=\id\:T_{u_0}{\calh}^{1,p}(W_{\alpha,\beta},X)
\to T_{u_0}{\calh}^{1,p}(W_{\alpha,\beta},X).
\endalignat
$$
Shrinking the balls $B_\alpha$ and $B_{\alpha\beta}$ if necessary,
we may assume that for all $\xi_\alpha\in B_\alpha$
and all $\xi_{\alpha\beta}\in B_{\alpha\beta}$
the restrictions of their images
$\psi_\alpha(\xi_\alpha)\big|_{W_{\alpha,\beta}}$ and
$\psi_{\alpha\beta}(\xi_{\alpha\beta})\big|_{W_{\alpha,\beta}}$
belong to  the image~$\psi'_{\alpha,\beta}(B'_{\alpha,\beta})$.
We consider the following holomorphic map:
$$
\matrix
\Psi\:&\prod_\alpha B_\alpha\times\prod_{\alpha<\beta}B_{\alpha\beta}
&\longrightarrow&\sum_{\alpha,\beta}{\calh}^{1,p}(W_{\alpha,\beta},u_0^*TX);
\\
\Psi\:&(v_\alpha,v_{\alpha\beta})&\mapsto&
\psi^{\prime\,-1}_{\alpha,\beta}
\bigl(\psi_\alpha(v_\alpha)|_{W_{\alpha,\beta}}\bigr)
-\psi^{\prime\,-1}_{\alpha,\beta}
\bigl(\psi_{\alpha\beta}(v_{\alpha\beta})|_{W_{\alpha,\beta}}\bigr).
\endmatrix
$$
As in  similar situations above, the map~$\Psi$
gives us  gluing conditions for local holomorphic maps
$\psi_\alpha(v_\alpha)\:V_\alpha\to X$ and
$\psi_{\alpha\beta}(v_{\alpha\beta})\:V_{\alpha\beta}\to X$.
Hence we can identify
$\calm\cap\prod_\alpha B_\alpha\times\prod_{\alpha<\beta}B_{\alpha\beta}$
with the set~$\Psi^{-1}(0)$.

 We shall now study in greater detail the behavior of $\Psi$
at the point
$y_0\in\prod B_\alpha\times\prod B_{\alpha\beta}$ that we are interested
in: namely,
$y_0=\bigl(\psi_\alpha^{-1}(u_0|_{V_\alpha}),
\psi_{\alpha\beta}^{-1}(u_0|_{V_{\alpha\beta}})\bigr)$, so that
$\Psi(y_0)=0\in\sum_{\alpha,\beta}{\calh}^{1,p}(W_{\alpha,\beta},u_0^*TX)$.
It is easy to see that the tangent space at~$y_0$
is
$$
T_{y_0}\Bigl(\prod B_\alpha\times\prod B_{\alpha\beta}\Bigr)
=\sum_\alpha{\calh}^{1,p}(V_\alpha,u_0^*TX)\oplus
\sum_{\alpha\beta}{\calh}^{1,p}(V_{\alpha\beta},u_0^*TX)\oplus\Bbb C^N
$$
and the differential $d\Psi(y_0)$ coincides on the term
$\sum{\calh}^{1,p}(V_\alpha,u_0^*TX)\oplus\sum{\calh}^{1,p}(V_{\alpha\beta},
\allowbreak
u_0^*TX)$  with the {\v C}ech codifferential~(3.1)
with respect to the cover $\{V_\alpha,V_{\alpha\beta}\}$ of the curve~$C_0$.
By Lemma~3.5 we can decompose
$\sum{\calh}^{1,p}(W_{\alpha,\beta},u_0^*TX)$
into the direct sum ${\calw}\oplus{\calq}$,
where ${\calw}=\Im(d\Psi(y_0))$ and ${\calq}$ is isomorphic to
$\roman H^1(C_0,u_0^*TX)$ and finite-dimensional.
Let $\Psi_{\calw}$ and $\Psi_{\calq}$
be the components of
$\Psi=(\Psi_{\calw},\Psi_{\calq})$
with respect to this decomposition,
and let~$\wt{\calm}:=\Psi_{\calw}^{-1}(0)$.
It follows from Lemma~3.5 and the implicit function theorem
that $\wt{\calm}$ is a complex submanifold of
$\prod B_\alpha\times\prod B_{\alpha\beta}$
with tangent space ${\calh}^{1,p}(C_0,u_0^*TX)\oplus\Bbb C^N$
at the point~$y_0\in\wt{\calm}$,
while   $\calm$ is defined in~$\wt{\calm}$
as the zero set  of the holomorphic map
$\Phi:=\Psi_{\calq}\big|_{\wt{\calm}}\:\wt{\calm}\to{\calq}
\cong\roman H^1(C,u_0^*TX)$.
This defines on~$\calm$ a structure of a Banach analytic set
of finite codimension.

\smallskip
To complete the proof it remains to construct the corresponding
family of nodal curves $\pi\:{\calc}\to\calm$
and a holomorphic map $F\:{\calc}\to X$.
Note that with  each ball $B_\alpha$ we can  associate in a natural way
the trivial family
$\pi_\alpha$: ${\calc}_\alpha:=B_\alpha\times V_\alpha\to B_\alpha$
and with each $B_{\alpha\beta}$ we can associate the holomorphic family
$\pi_{\alpha\beta}\:{\calc}_{\alpha\beta}\to B_{\alpha\beta}$
with fiber $\pi_{\alpha\beta}^{-1}(v_{\alpha\beta})=A_{\lambda_{\alpha\beta}}$,
where $\lambda_{\alpha\beta}$ is uniquely determined by the relation
$\psi_{\alpha\beta}(v_{\alpha\beta})=(\lambda_{\alpha\beta},u_{\alpha\beta})
\in\calm(V_{\alpha\beta},X)$.

We extend these families to the families
$\ti\pi_\alpha\:\wt{\calc}_\alpha\to\prod B_\alpha\times
\prod B_{\alpha\beta}$
and
$\ti\pi_{\alpha\beta}:\wt{\calc}_{\alpha\beta}\to\prod B_\alpha\times
\prod B_{\alpha\beta}$.
Clearly, $\wt{\calc}_\alpha$ and $\wt{\calc}_{\alpha\beta}$
can be glued together canonically, producing  a global
family of nodal curves
$\ti\pi\:\wt{\calc}\to\prod B_\alpha\times\prod B_{\alpha\beta}$.
Note that ${\calc}$ is  a Banach manifold.
Furthermore, we obtain well-defined holomorphic maps
$F_\alpha\:\wt{\calc}_\alpha\to X$ and
$F_{\alpha\beta}\:\wt{\calc}_{\alpha\beta}\to X$
such that $F(v_\alpha,v_{\alpha\beta},z):=\psi_\alpha(v_\alpha)[z]$
for each $z\in V_\alpha$
and a similar relation holds for~$F_{\alpha\beta}$.

We define ${\calc}$ to be the restriction
${\calc}:=\wt{\calc}|_{\!\calm\!}$.
We observe that the restriction of the {\it trivial\/}
holomorphic family
$\wt{\calc}_\alpha=V_\alpha\times\prod B_\alpha\times\prod B_{\alpha\beta}$
to~$\calm$ is also a trivial holomorphic family.
It follows that ${\calc}$ is a BASFC in a neighborhood
of each point ~$y\in{\calc}\cap\wt{\calc}_\alpha$. In a
similar way, each holomorphic family of curves
$\wt{\calc}_{\alpha\beta}$ can be distinguished
in the trivial bundle
$\Delta^2\times\prod B_\alpha\times\prod B_{\alpha\beta}
\overset\pr\to\longrightarrow\prod B_\alpha\times\prod B_{\alpha\beta}$
by the condition
$z_1\cdot z_2-\lambda_{\alpha\beta}=0$,
where $\lambda_{\alpha\beta}\:B_{\alpha\beta}\to\Delta$
is the holomorphic parameter of the deformation
of the complex structure on~$V_{\alpha\beta}$
and $(z_1,z_2)$ are the standard coordinate functions on~$\Delta^2$.
Hence ${\calc}$ is also a BASFC      
in a neighborhood of each point
$y\in{\calc}\cap\wt{\calc}_{\alpha\beta}$.
Since $\calm$ has been in fact defined
by the condition that the local maps
$F_\alpha$ and $F_{\alpha\beta}$ coincide,
there exists a global holomorphic map
from~${\calc}$ to~$F\:{\calc}\to X$.

Properties (a), (b) and (d) from Theorem~3.4 follow now  for
$\calm$, ${\calc}$, and $F$
directly from the construction,
and ~(c) follows by Proposition~3.2.

The proof of the required variant of the continuity principle
for meromorphic maps (see  Theorem~4.2
in the next section or Theorem~5.1.3 in~\cite{2})
uses the following consequence of the main results in this section.

Assume that a sequence $(C_n,u_n)$ of irreducible stable curves over~$X$
converges to a stable curve~$(C_\infty,u_\infty)$.

\proclaim{Lemma 3.8}
For some positive integer~$N$ there exist
a complex \rom(maybe, singular\rom) surface~$Z$
and 
holomorphic maps $\pi_Z\:Z\to\Delta$ and $F\:Z\to X$
 defining  a holomorphic family of  stable nodal
curves over~$X$ joining $(C_N,u_N)$ and~$(C_\infty,u_\infty)$.
More precisely, the following results hold\rom:
\roster
\item"(1)" for each $\lambda\in\Delta$ the fiber
$C_\lambda=\pi_Z^{-1}(\lambda)$
is a connected nodal curve with boundary~$\partial C_\lambda$
and the pair $(C_\lambda,u_\lambda)$ with $u_\lambda:=F\big|_{C_\lambda}$
is a stable curve over~$X;$

\item"(2)" all  $C_\lambda$, except for  finitely many curves,
are connected and smooth\rom;

\item"(3)" $(C_0,u_0)$ is equal to~$(C_\infty,u_\infty)$
and there exists $\lambda_N\in\Delta$ such that
$(C_{\lambda_N},u_{\lambda_N})\allowmathbreak =(C_N,u_N);$

\item"(4)" there exist open sets $V_1,\dots,V_m$ in~$Z$
such that each $V_j$ is biholomorphic to~$\Delta\times A_j$
for some annulus $A_j;$
moreover, the diagram
$$
\CD
V_j@>\cong>>\Delta\times A_j
\\
@VV\pi V @VV\pi_\Delta V
\\
\Delta@=\Delta
\endCD
$$
is commutative,
each annulus $C_\lambda\cap V_j\cong\{\lambda\}\times A_j$
is adjacent to a unique boundary component of~$\partial C_\lambda$,
and the number $m$ of the domains~$V_j$
is equal to the number of boundary components of each~$C_\lambda$.
\endroster
\endproclaim

\smallskip\noindent\bf Remark. \rm 
This lemma was stated without proof
in Proposition~5.1.1 of~\cite{2}.
\endremark

\smallskip\noindent\bf Proof. \rm 
Assume that  $\pi_{\!\calc\!}\:{\calc}\to\calm$ and
$\ev\:{\calc}\to X$
define the complete family of holomorphic deformations
of the stable nodal curve $(C_\infty,u_\infty)$ over~$X$
constructed in Theorem~3.4.
Assume that $\lambda^*\in\calm$ parameterizes
 the curve~$(C_\infty,u_\infty)$.
We consider a sequence $\lambda_n\to\lambda^*$ in~$\calm$
such  that $(C_n,u_n)\cong\bigl(\pi_{\!\calc\!}^{-1}(\lambda_n),
\ev|_{\pi_{\!\calc\!}^{-1}(\lambda_n)}\bigr)$
for all sufficiently large~$n$.
Since $(C_\infty,u_\infty)=\bigl(\pi_{\!\calc\!}^{-1}(\lambda^*),\
\ev|_{\pi_{\!\calc\!}^{-1}(\lambda^*)}\bigr)$
can be lifted  in a neighborhood of the boundary of~$\widehat U$,
shrinking $\calm$
we can assume that this is true for all~$\lambda\in\calm$.

By construction, the space $\calm$ is a BASFC      
and therefore Theorem~3.3 applies.
In~particular, $\calm$ has  finitely many
irreducible components  at~$\lambda^*$.
Let $\calm_1$ be a component of~$\calm$
that contains infinitely many~$\lambda_n$.
We can represent $\calm_1$ in a neighborhood of~$\lambda^*$
as a proper ramified cover $\pi_1\:\calm_1\to B_1$
over a Banach ball~$B_1$.
If $\Delta$ is an arbitrary embedded disc in~$B_1$,
then $\pi_1^{-1}(\Delta)$ is a one-dimensional analytic subset
with irreducible components also parametrized by discs.
Hence there exists a holomorphic map
$\varphi\:\Delta\to\calm_1$
passing through $\lambda^*$ and $\lambda_N$
for some $N\gg 1$.
The inverse image of the family
$\pi_{\!\calc\!}\:{\calc}\to\calm$
under~$\varphi$ is a holomorphic family of stable nodal curves over~$X$,
which has the
 total space $\pi_Z\:Z\to\Delta$ equipped with a map $F\:Z\to X$ and
 contains $(C_\infty,u_\infty)$ and~$(C_N,u_N)$.

Since $C_N$ is smooth,  the general fiber
$C_\lambda=\pi_Z^{-1}(\lambda)$ is also smooth.
Shrinking the disc if necessary,
we may assume that the $C_\lambda$ are singular
only for finitely many ~$\lambda\in\Delta$.
The other properties from~(1)--(4) follow by   the
construction of the family~$\pi_Z\:Z\to\Delta$
and the map~$F\:Z\to X$.

\head
\S\,4. Continuity principle and the proof of the main theorem
\endhead

 We shall now formulate the continuity principle from~\cite{2}
required for the proof of the main theorem.

\definition{Definition 4.1}
A Hermitian complex manifold $(X,h)$
is said to be  {\it disc-\allowlinebreak convex}
if for each sequence $(C_n,u_n)$ of stable curves over~$X$
parametrized by the same real surface~$\Sigma$
and having the properties

\roster
\item"(1)" the curves $C_n$ are irreducible for all~$n$
 and the boundaries $\partial C_n$ are not empty;

\item"(2)" the areas $\operatorname{area}_h[u_n(C_n)]$
are uniformly bounded;

\item"(3)" the $(C_n,u_n)$ converge in a neighborhood
of their boundaries~$\partial C_n$,
\endroster
there exists a compact subset $K$ of $  X$
 containing all~$u_n(C_n)$.

In particular, all compact manifolds are disc-convex.
The property to be disc-convex remains invariant
under the  replacement of the Hermitian metric~$h$
by an equivalent metric~$h'$, $c\cdot h\le h'\le C\cdot h$.
We point out also that by Theorem~3.1,
such a sequence $(C_n,u_n)$
contains a subsequence converging in the Gromov topology.

Let $U$ be a domain in a complex manifold~$X$,
and let $Y$ be a complex space.

\definition{Definition 4.2}
The {\it envelope of meromorphy\/} of~$U$
with respect to $Y$ is the maximal domain
$(\widehat U_Y,\widehat\pi)$ over $X$
containing~$U$ (that is, there exists an embedding $i\:U\to\widehat U_Y$
with~$\widehat\pi\circ i=\id$)
such that each meromorphic map $f\:U\to Y$
extends 
to a meromorphic map~$\widehat f\:\widehat U\to Y$.

The Cartan--Thullen construction (cf.~\cite{1})
applies here to prove the existence and the uniqueness of the envelope.

\proclaim{Proposition 4.1}
For each domain $U$ in a complex manifold~$X$
and  each complex space~$Y$
there exists an  envelope of meromorphy
$(\widehat U_Y,\widehat\pi)$ of  $U$
with respect to~$Y$.
\endproclaim

The following result was proved in \cite{2}; Theorem 5.1.3.

\proclaim{Theorem 4.2}
For each domain $U\subset X$ in a disc-convex
complex Hermitian surface $(X,h)$
and  each disc-convex K\"ahler space~$Y$
the envelope of meromorphy $(\widehat U_Y,\widehat\pi)$
equipped with the Hermitian metric~$\widehat\pi^*h$
is also disc-convex.
\endproclaim

\smallskip\noindent\bf Remark. \rm 
Let us explain the meaning of this theorem using the following example.
Assume that $X$ is disc-\allowlinebreak convex
and consider a domain~$U\subset X$.
Let $f\:U\to Y$ be a meromorphic map
and let $\{(C_n,u_n)\}$ be a sequence of stable holomorphic curves over~$X$
convergent to~$(C_\infty,u_\infty)$ in the Gromov topology
such that the images of the boundaries~$u_n(\partial C_n)$
and~$u_\infty(\partial C_\infty)$ lie in~$U$.
Assume further that $f$ 
extends 
along each~$u_n(C_n)$.
This means that there exist a complex surface~$V_n$ containing~$C_n$
and a holomorphic, locally biholomorphic map
$u'_n\:V_n\to X$ such that $u'_n\big|_{C_n}=u_n$
and $f$ 
extends 
meromorphically from $u_n^{\prime\,-1} U$
to the whole of~$V_n$.

The last assumption is equivalent to the condition that
the curves~$(C_n,u_n)$ can be lifted to the envelope $\widehat U$,
that is, there exist holomorphic maps
$\widehat u_n\:C_n\to\widehat U$
such that $\widehat\pi\circ\widehat u_n=u_n$.
In other words,
we may take~$V_n$ to be a~neighborhood
of the lifting of~$C_n$ in~$\widehat U$.
It is easy to see that the lifted curves
$(C_n,\widehat u_n)$ are stable over~$\widehat U$,
have uniformly bounded area,
and converge near the boundary~$\partial C_n$.
By Theorem~4.2 and  Gromov's compactness theorem
(Theorem~3.1)
a subsequence of $(C_n,\widehat u_n)$ converges
to a $\widehat U$-stable curve~$(C_\infty,\widehat u_\infty)$ such
that $\widehat\pi\circ\widehat u_\infty=u_\infty$.
This means that $f$ 
extends 
along~$u_\infty(C_\infty)$.
Thus, Theorem~4.2 is a generalization of Levi's continuity
principle.
\endremark

We shall actually prove a stronger result  than the main theorem.
Namely, in place of meromorphic functions
(that~is, meromorphic maps into~$\Bbb{CP}^1$)
we shall consider the general case of meromorphic maps
into an arbitrary disc-convex K\"ahler space~$Y$.

\proclaim{Theorem 4.3}
Let $u\:S^2\to X$ be a symplectic immersion of the sphere~$S^2$
into a disc-convex K\"ahler surface~$X$
such that $M:=u(S)$ has only positive double points.
Assume that~$c_1(X)[M]>0$.
Then the envelope of meromorphy $(\widehat U_Y,\pi)$
of a neighborhood $U$ of~$M$
with respect to a disc-convex K\"ahler space~$Y$
contains a rational curve~$C$ with~$\pi^*c_1(X)[C]>0$.
\endproclaim

\smallskip\noindent\bf Proof. \rm 
Let $u\:S^2\to X$ be a symplectically immersed sphere
with only positive self-intersections.
Let $U$ be a relatively compact subdomain of~$X$ containing~$M:=u(S^2)$.
We denote its envelope of meromorphy with respect to~$Y$
by $(\widehat U_Y,\widehat\pi_Y)$.

\smallskip
{\it Step\/} 1. {\it
There exists an \rom($\omega$-tamed\/\rom)
almost complex structure $J_0\in\calj_U$
such that $M$ is a $J_0$-holomorphic curve}.

This has been proved in~\cite{2}; Lemma~1.1.2.
Moreover, there exists a smooth homotopy~$h\:[0,1]\to\calj_U$
between~$J_0=h(0)$ and~$J_{\text{\rom{st}}}=h(1)$.
We set~$M_0:=M$.
As in Lemma~2.5, let $\calm_h(M_0,J_0)$
be     the component of~$\calm_h$
passing through~$(M_0,J_0)$.

\smallskip
{\it Step\/} 2. {\it Assume that the component
$\calm_h(M_0,J_0)$ is non-compact}.

Then part~(iii) of Lemma~2.5 says  that there exists
a continuous curve~$\gamma\:[0,1)\to\calm_h$
starting at~$(M_0,J_0)$ and having  property~(b)
from part~(iii) of Lemma~2.5.
We consider now  $J_n$-holomorphic spheres~$M_n$ that
 form a discrete set in~$\calm_h(M_0,J_0)$,
as $J_n$ converges to~$J^*\in\calj_U$.

If $M_n\cap U=\varnothing$ for some~$n$,
then $M_n$ is the required rational curve
because $J_n=J_{\text{\rom{st}}}$ on~$\widehat U\setminus U$.

If $M_n\cap U\ne\varnothing$ for all~$n$,
then there exists a subsequence,
which we still denote by $M_n$, that
converges in the Gromov topology
to some (in general, reducible) curve~$M^{(1)}$.
If $M^{(1)}$ has an irreducible component~$M^{(1)}_0$
lying outside ~$U$ and satisfying ~$c_1(X)[M^{(1)}_0]>0$,
then $M^{(1)}_0$ is the required rational curve.

Otherwise, there exists a component $M^{(1)}_0$
of the limiting curve~$M^{(1)}$
intersecting~$U$ and such that  $c_1(X)[M^{(1)}_0]>0$.
We
repeat step~2 with $M^{(1)}_0$ in place of~$M_0$.
Since the area of pseudoholomorphic curves
is bounded from below (see~\cite{6}),
after several repetitions of step~2
we shall either find a rational curve in the envelope~$\widehat U_Y$,
or arrive at the  situation described below.

\smallskip
{\it Step\/} 3. {\it $\calm_h(M_0,J_0)$ is compact}.

By part~(iii) of Theorem~2.6 and Corollary~2.7
we obtain continuous paths $(M^n_t,J^n_t)$
such that:

\goodbreak 

\roster
\item"(1)" $M^n_0=M_0$ for all $n$;

\item"(2)" $J^n_0=J_0$ for all $n$;

\item"(3)" $J^n_1\to J_{\text{\rom{st}}}$.
\endroster

By  Gromov's compactness theorem
a subsequence of $M^n_1$ converges
to a $J_{\text{\rom{st}}}$-holomorphic nodal curve~$C^*$
over~$\widehat U_Y$.
We choose an irreducible component $C$ of~$C^*$
with~$\pi^*c_1(X)[C]>0$.
Then $C$ is the required rational curve in~$\widehat U_Y$.

\comment
\spaceskip=4pt plus3.5pt minus 1.5pt
\spaceskip=5pt plus4pt minus 2pt
\font\csc=cmcsc10

\lleftskip=3.8\parindent
\length=\hsize \advance\length-\lleftskip
\def\entry#1#2#3#4\par{\parshape=2  0pt  \hsize%
\lleftskip \length%
\noindent\hbox to \lleftskip%
{\bf[#1]\hfill}{\csc{#2 }}{\sl{#3}}#4%
\medskip
}
\ifx \twelvebf\undefined \font\twelvebf=cmbx12\fi

\bigskip\bigskip
\bigskip\bigskip
\centerline{\twelvebf References}
\bigskip

\entry{Ab}{Abikoff W.:}{Real analytic theory of Teichmuller space} Springer,
LNM 820, (1980).

\entry{Ba}{Barlet D.:}{Espace analytique reduit des cycles analytiques
complexes compacts d'un espace analytique complexe de dimension finie .}
Semiar Norguet IX, Lect. Notes Math., {\bf 482}, 1--157, (1975).

\entry{B-P-V}{Barth W., Peters C., Van de Ven A.:}{Compact complex sufaces}
Springer, Berlin, (1984).

\entry{D-G}{Docquier F., Grauert H.:}{ Levisches Problem und Rungescher
Satz f\"ur Teilgebiete Steinscher Mannigfaltikeiten.} Math. Ann. {\bf 140},
94-123, (1960).

\entry{Gau}{Gauduchon P.:}{Les metriques standard d'une surface a premier
nombre de Betti  pair.} Asterisque.\ Soc.\ Math.\ France.  {\bf126},
129--135, (1985).

\entry{Gr}{Griffiths P.:}{Two theorems on extensions of holomorphic
mappings.} Invent.\ math. {\bf14}, 27--62, (1971).

\entry{Gro}{Gromov M.:}{Pseudo holomorphic curves in symplectic manifolds.}
Invent. mat. {\bf82}, 307--347 (1985).

\entry{H-L-S}{Hofer H., Lizan V., Sikorav J.-C.}{On genericity
for holomorphic curves in 4-dimensional almost-complex manifolds.}
Preprint (1994).

\entry{H}{Hummel C.:}{Gromov's compactness theorem for pseudo - holomorphic
curves} Progress in Math., {\bf Vol.151}, Basel, Birk\"auser, 1997.

\entry{D-S}{McDuff D., Salomon D.:}{$J$-holomorphic curves and quantum
cohomology}AMS University Lecture Series, {\bf Vol.6}, 1994.

\entry{Mo}{Morrey C.~B.:}{Multiple integrals in the calculus of
variations.}Springer, New-York,(1966).

\entry{I-1}{Ivashkovich, S.:}{Extension of analytic object by the method of
Cartan - Thullen.} Proc. Conf. Compl. Anal. Math. Phys., 53--61,
Krasnojarsk, (1988).

\entry{I-2}{Ivashkovich, S.:}{The Hartogs-type extension theorem for
meromorphic maps into compact K\"ahler manifolds.} Invent. Math. {\bf109},
47--54 (1992).

\entry{I-S}{Ivashkovich S., Shevchishin V.:}{Pseudo-holomorphic curves and
envelopes of meromorphy of two-spheres in $\cc\pp^2$.} Preprint, Bochum
SFB-237 (1995) http://xxx.lanl.gov 9804014.

\entry{Le}{Levi E.:}{Studii sui punti singolari essenziali delle funzioni
analitiche di due o pi\'u variabili complesse.} Annali di Mat. pura ed appl.
 {\bf17}, n 3, 61--87 (1910).

\entry{N}{Nemirovsky S.:}{Stein domains and imbedded Riemann surfaces.}
Preprint, Moscow, (1997).

\entry{P-W}{Parker Th., Wolfson, J.:}{Pseudo-holomorphic maps and bubble
trees} J.~of Geometrical Analysis, {\bf 3}, 63--98, (1993).

\entry{Ra}{Ramis J.-P.:}{Sous-ensembles analytiques d'une veri\'et\'e
banachique complexe.}, Springer, Berlin, 1970.

\entry{Re}{Remmert R.:}{Holomorphe und meromorphe Abbildungen komplexer
R\"aume.} Math.\ Ann. {\bf133}, 328-370, (1957).

\entry{Fe}{Federer G.:}{Geometric measure theory} Springer (1989).
\endcomment

\bigskip\bigskip\noindent
\vtop{\hsize=.4\hsize
Universit\'e de Lille-I\par
U.F.R. de Math\'ematiques\par
Villeneuve d'Ascq Cedex\par
59655 France\par
ivachkov\@gat.univ-lille1.fr
}\hfill
\vtop{\hsize=.4\hsize
Ruhr-Universit\"at, Bochum\par
Fakult\"at f\"ur Mathematik \par
Universit\"atsstasse 150\par
D-44780 Bochum, \par
Germany\par
sewa\@cplx.ruhr-uni-bochum.\par
de
}

\bigskip\bigskip
\centerline{\vtop{\hsize=.4\hsize
IAPMM NAS Ukraine\xpar
Naukova str. 3b,\xpar
290053 Lviv,\xpar
Ukraine
}}

\end